\documentclass[final]{amsart}
\usepackage{amsfonts,amssymb,amsmath,amsgen,amsthm,datetime,float}
\usepackage[top=0.6in, bottom=0.8in, scale=0.75, centering, headheight=10pt]{geometry}
\usepackage{hyperref}
\usepackage{cleveref}
\usepackage{xcolor}

\definecolor{mypetrol}{RGB}{0,125,145}

\hypersetup{colorlinks=true,linkcolor=blue,citecolor=mypetrol}
\usepackage{enumitem}
\usepackage{showkeys}
\usepackage{latexsym,amssymb}
\usepackage{graphicx}
\usepackage{tikz}      
\usetikzlibrary{positioning}  
\usetikzlibrary{arrows.meta}
\usetikzlibrary{calc}
\usepackage{etoolbox}
\theoremstyle{plain}
\newtheorem{theorem}{Theorem}[section]
\newtheorem{defi}[theorem]{Definition}
\newtheorem{assumption}[theorem]{Assumption}
\newtheorem{lemma}[theorem]{Lemma}
\newtheorem{corollary}[theorem]{Corollary}
\newtheorem{prop}[theorem]{Proposition}

\theoremstyle{remark}
\newtheorem{remark}[theorem]{Remark}
\newtheorem{notation}[theorem]{Notation}

\newcommand{\ds}{\displaystyle}
\newcommand{\N}{\mathbb{N}}
\newcommand{\R}{\mathbb{R}}

\newcommand{\C}{\mathbb{C}}
\newcommand{\Lp}[3]{L^{#1}(#2,#3)}

\newcommand{\siegel}[2]{\mathrm{Sieg}_{#1}(#2)}
\newcommand{\symetrique}[2]{\mathcal{S}_{#1}(#2)}

\newcommand{\normLp}[3]{\| #1 \|_{L^{#2}(#3)}}
\newcommand{\normtriple}[1]{||| #1 |||}
\newcommand{\normfrob}[1]{||| #1|||_{\mathrm{F}}}
\newcommand{\normspec}[1]{||| #1|||_{2}}
\newcommand{\transpose}[1]{#1^{\mathsf{T}}}

\newcommand{\diff}{\mathrm{d}}

\newcommand{\Ima}{\mathrm{Im}}
\newcommand{\Rea}{\mathrm{Re}}
\newcommand{\eps}{\varepsilon}
\renewcommand{\d}{\partial}
\newcommand{\psiapp}{\psi^{\varepsilon}_{\mathrm{app}}}
\newcommand{\psiWP}{\psi_{_{\mathrm{WP}}}^{\varepsilon}}
\newcommand{\nl}[1]{\lambda \, \eps \, #1 \, \log \, \lvert #1 \lvert^2}
\newcommand{\matinv}[2]{\mathrm{GL}_{#1}(#2)}
\newcommand{\valeurmoyenne}[4]{\ds\int_{\R^d} #1(x) \exp \left( - \frac{1}{\varepsilon} \Big( \transpose{(x - #2)} #3 (x - #2) + 2 \Ima (#4) \Big) \right) \diff x}

\newcommand{\smalleigenvalue}[1]{\lambda_{\mathrm{min}}\Big[#1\Big]}
\newcommand{\largeeigenvalue}[1]{\lambda_{\mathrm{max}}\Big[#1\Big]}

\newcommand{\deter}[1]{\mathrm{det} \Big[ #1 \Big]}

\numberwithin{equation}{section}

\begin{document}

\title{Variational approximation by complex Gaussians for the nonlinear logarithmic equation}

\author[M. Curely]{Marianne Curely}
\address[M. Curely]{Univ Angers, CNRS, LAREMA, SFR MATHSTIC, F-49000 Angers, France}
\email{marianne.curely@univ-angers.fr}

\author[L. Hari]{Lysianne Hari}
\address[L. Hari]{Université Marie et Louis Pasteur, CNRS, LmB (UMR 6623), F-25000 Besançon, France}
\email{lysianne.hari@umlp.fr}

\begin{abstract} We consider the semiclassical logarithmic Schrödinger equation, with an at most quadratic potential and with Gaussian initial states. For a critical nonlinear term in the semiclassical regime, we analyze an approximation based on the variational evolution of Gaussian wave packets. We also prove a nonlinear superposition principle for solutions issued from the sum of two Gaussian wave packets initially localized at distinct points in phase space. For both a single Gaussian wave packet and a superposition of two Gaussian wave packets, we compare the variational approximation with the standard coherent-state approximation and highlight their respective features. The required assumptions are weakened in this setting without becoming technically involved.
\end{abstract}

\maketitle

\vspace{-0.3cm}



\tableofcontents
\addtocontents{toc}{\protect\setcounter{tocdepth}{1}}

\vspace{-0.5cm}

\section{Introduction}

\medskip

\noindent Let \(d \geqslant 1\) be an integer. We consider the following logarithmic Schrödinger equation in the semiclassical regime~$0 < \varepsilon \ll 1$
\begin{equation}
\label{eq:logequation}
\left\{
        \begin{array}{r c l l}
        i \varepsilon \partial_t \psi^{\varepsilon} + \ds\frac{\varepsilon^2}{2} \Delta \psi^{\varepsilon} & = & V \psi^{\varepsilon} + \nl{\psi^{\varepsilon}} \\
        \psi^{\varepsilon}_{|{t=0}} & = & \psi_{0}^{\varepsilon}
        \end{array}
    \right.
\end{equation}
\noindent where $\psi^{\varepsilon}:=\psi^{\varepsilon}(t,x)$ is a complex-valued \textit{wave function}, $\lambda \in \R$, the function $V:=V(x)$ is a real-valued \textit{potential} and the initial data $\psi_{0}^{\varepsilon}=\psi_{0}^{\varepsilon}(x)$ is a bounded family of $\Lp{2}{\R^d}{\C}$. Throughout this work, the potential $V$ satisfies the following assumption.

\smallskip

\begin{assumption}
\label{hyp:potential}
The function \( V \) is smooth and at most quadratic, in the sense that
\begin{equation*}
    \forall \, \beta \in \mathbb{N}^d, \,\, |\beta| \geqslant 2, \quad \partial^{\beta} V \in \Lp{\infty}{\R^d}{\R}.
\end{equation*}
\end{assumption}

\noindent Under our assumptions on the potential \(V\), the results of \cite{CarlesFerriere21} apply, allowing us to establish the existence of global solutions to~\eqref{eq:logequation} with an initial data chosen in a suitable subset of $L^2(\R^d)$. The solution enjoys conservation properties as both the mass and the energy are conserved for all time. More details are provided in Appendix~\ref{appendix:existence}.

\smallskip

\noindent Since its introduction in \cite{BBM76}, the logarithmic Schrödinger equation has been used to model a wide range of physical phenomena (see for example~\cite{PhDFerriere} and \cite{carles2025log} for a detailed list of references). From a mathematical point of view, the first results for $V=0$ are due to Cazenave and Haraux in \cite{cazenave80, cazenave83}. They overcame the difficulties arising from the logarithmic nonlinearity, which is singular for small values of the modulus and does not have a definite sign, by considering an appropriate energy space to solve the Cauchy problem. Their analysis relies on the following key estimate. For all~$(z_1, z_2) \in \C^2$, we have
    \begin{equation}
    \label{lem:majlog}
        \Big| \mathrm{Im} \Big( \Big( z_2 \log | z_2 |^2 \, - \, z_1 \log | z_1 |^2 \Big) (\overline{z_2} - \overline{z_1}) \Big) \Big| \leqslant 2 |z_2 - z_1 |^2.
    \end{equation}

\noindent The Cauchy problem was recently revisited by Hayashi and Ozawa in \cite{HayashiOzawa25} with a constructive method, which was a pioneering work to handle more general cases, such as low regularity solutions and their uniqueness (see \cite{CHO24,Hayashi26} and also \cite{Bensaid26} for competition between logarithmic and power type nonlinearities).

\smallskip

\noindent A long interval separates the early works of Cazenave and Haraux from the more recent developments on the equation in the 2010s that extend over a broad spectrum of research topics (see for instance \cite{GLN10, AMS14log, Ardila16, CarlesGallagher18}).
In addition to the references mentioned above, we refer the reader to more recent works as in \cite{Ferriere20,Ferriere21multi,Bensaid26-2} for deeper discussions about special solutions like Gaussons, multi-solitons and orbital stability. The survey \cite{CarlesEMS} provides a recent account of the results around the equation, before the latest works. The equation was also actively studied from a numerical point of view as in \cite{BCST19, BCST19-2, BCST22, CC22, CC23} and with higher order methods in more recent papers \cite{WYZ24, ZW25, LWZ26, WY26}.

\smallskip

\noindent Among the contributions, the logarithmic Schrödinger equation with a semiclassical parameter has also been studied (see  \cite{CarlesNouri17,Ferriere21,Ferriere22} for results involving Wigner measures, Wigner transforms and WKB states). More recently, the propagation of coherent states was investigated in~\cite{carles2025log}, which is particularly relevant to our work. The propagation of wave packets allows to construct approximate solutions and have been extensively studied. The case of the scalar linear Schrödinger equation is covered in~\cite{Hagedorn1980, livreRobertComb2012}. For nonlinear scalar Schrödinger equations, we refer, among others, to~\cite{carlesfermanian2011, APPP11, CaoCarles, Carles12interaction}. In~\cite{carles2025log}, the authors indeed construct a suitable approximate solution, at leading order, using coherent states for the semiclassical logarithmic Schrödinger equation with a more general nonlinear term than the one considered here. They assume that the potential $V$ satisfies Assumption~\ref{hyp:potential} and they consider the following Cauchy problem
\begin{equation}
\label{eq:logequationalpha}
i \varepsilon \partial_t \psi^{\varepsilon} + \ds\frac{\varepsilon^2}{2} \Delta \psi^{\varepsilon} = V \psi^{\varepsilon} + \lambda \, \eps^{\alpha} \, \psi^{\varepsilon} \, \log \, \lvert \psi^{\varepsilon} \lvert^2 \quad ; \quad \psi^{\varepsilon}_{|{t=0}} = \psi_{0}^{\varepsilon}
\end{equation}
\noindent with $\alpha \geqslant 1$ and where the initial data $\psi_{0}^{\varepsilon}$ is a coherent state of the form
\begin{equation*}
\psi_{0}^{\varepsilon}(x) = \varepsilon^{-\frac{d}{4}} \, u_0 \left(\dfrac{x-q_0}{\sqrt{\varepsilon}}\right) e^{\frac{i}{\varepsilon} p_0 \cdot (x-q_0)}.
\end{equation*}
Here, $(q_0,p_0)$ is a vector of $\R^{2d}$ and the function $u_0$ satisfies suitable regularity assumptions. In this setting, the authors of~\cite{carles2025log} identify a critical exponent $\alpha_{\mathrm{crit}}=1$, beyond which the influence of nonlinear effects on the propagation becomes negligible. In the case where~$\alpha>\alpha_{\mathrm{crit}}$, for such initial data, the solution behaves like the linear one at leading order, that is
\begin{equation}
\label{eq:approxlineaire}
    \forall t\geqslant 0, \quad \lim\limits_{\varepsilon \to 0} \normLp{\psi^{\varepsilon}(t) - \varphi^{\varepsilon}(t)}{2}{\R^d} = 0,
\end{equation}
where $\psi^{\varepsilon}(t)$ is the exact solution to \eqref{eq:logequationalpha} and $\varphi^{\varepsilon}$ is the solution to the linear equation
\begin{equation*}
i \varepsilon \, \partial_t \varphi^{\varepsilon} + \ds\frac{\varepsilon^2}{2} \Delta \varphi^{\varepsilon} = V \varphi^{\varepsilon} \quad ; \quad \varphi^{\varepsilon}_{|{t=0}} = \psi_{0}^{\varepsilon}.
\end{equation*}

\noindent This is what is called the \textit{subcritical case}. Moreover, in both \textit{subcritical case} and \textit{critical case} (that is the case considered in this article and which corresponds to $\alpha = \alpha_{\text{crit}} = 1$), in~\cite[Proposition 1.2]{carles2025log}, the authors also prove that the wave function $\psi^{\varepsilon}$ can be approximated in $L^2(\R^d)$, at leading order, by a wave packet. It holds that
\[ \forall t \geqslant 0, \quad \lim\limits_{\varepsilon \to 0} \normLp{\psi^{\varepsilon}(t) - \mathrm{WP}_{z(t)}^{\varepsilon} \, u(t)}{2}{\R^d} = 0,
\]
where $\mathrm{WP}_{z(t)}^{\varepsilon} \, u(t)$ denotes a wave packet, defined for all $x \in \R^d$ by
\[\mathrm{WP}_{z(t)}^{\varepsilon} \, u(t,x) = \varepsilon^{-d/4}
u\left(t,\dfrac{x-q_{_{\mathrm{WP}}}(t)}{\sqrt{\varepsilon}}\right)
e^{\frac{i}{\varepsilon}p_{_{\mathrm{WP}}}(t)\cdot(x-q_{_{\mathrm{WP}}}(t)) + \frac{i}{\varepsilon}S(t)},\]
with $z(t) = (q_{_{\mathrm{WP}}}(t),p_{_{\mathrm{WP}}}(t))$ are the usual \textit{classical trajectories} and $S(t)$ is the \textit{classical action}, satisfying
\[
\dot{q}_{_{\mathrm{WP}}}= p_{_{\mathrm{WP}}}, \qquad
\dot{p}_{_{\mathrm{WP}}}=-\nabla V(q_{_{\mathrm{WP}}}) \quad \text{and} \quad
S(t) = \int_0^t
\left(
\frac{|p_{_{\mathrm{WP}}}(s)|^2}{2}-V(q_{_{\mathrm{WP}}}(s))
\right)\,ds,
\]
with \( (q_{_{\mathrm{WP}}}(0),p_{_{\mathrm{WP}}}(0),S(0))=(q_0,p_0,0).\) Moreover, the \textit{profile} $u$ is the solution to the following linear Schrödinger equation with a time-dependent harmonic potential
$$i\d_t u + \dfrac{1}{2}\Delta u = \dfrac{1}{2} \left\langle y, \textrm{Hess\;}V(q_{_{\mathrm{WP}}})y\right\rangle u =0 \quad ; \quad u_{|{t=0}} = u_0.$$

\noindent The approximation can be proved up to the \textit{Ehrenfest time}~$T^\varepsilon= c \, \log (\varepsilon^{-1})$, although such a long time of validity is not typically achieved for other nonlinearities, due to technical obstructions (see~\cite{Adiab, carlesfermanian2011, hari2013coherent} for example). On the other hand, in the \textit{critical case} $\alpha=1$, nonlinear effects appear through the profile, which is now the solution to the nonlinear equation
\begin{equation*}
i\partial_t u + \dfrac{1}{2}\Delta u = \dfrac{1}{2}\left\langle y,\operatorname{Hess}V(q_{_{\mathrm{WP}}})y\right\rangle u + \lambda u\log |u|^2 \quad ; \quad u_{|{t=0}} = u_0.
\end{equation*}

\noindent Our aim is to recover the results of~\cite{carles2025log} for more general Gaussian initial states, using a \textit{variational method} based on the time-dependent Dirac--Frenkel/McLachlan variational principle, hereafter referred to as the DFMVP. This approach was described in~\cite{lasserlubich2020} for the linear scalar Schrödinger equation and in~\cite{lasserChimistesMag} for the magnetic case. For a more detailed discussion of the DFMVP, we refer the reader to the original works~\cite{dirac1930theory, frenkel1934, mclachlan1964} and to~\cite{KramerSaraceno,broeckhove1988equivalence,Raab} for a proof of the equivalence between the Dirac--Frenkel and McLachlan principles in specific settings. Pioneer rigorous mathematical analyses can be found in~\cite{LubichBook} and related results in \cite{FaouLubich}. Related methods based on these variational principles have been widely used in the chemistry literature, starting with Heller's work~\cite{Heller76, Heller81} and the subsequent developments in~\cite{CoalsonKarplus}. These methods include the multi-configuration time-dependent Hartree (MCTDH) approach and its variants (see for example~\cite{MCTDH, MCTDH99, MCTDHreport, vMCG}).

\smallskip

\noindent The Gaussian approximation that we build through this method is different from the one obtained in~\cite{carles2025log} as it does not concentrate around the same classical quantities. However, we will see that the quantities are close in the semiclassical limit and that the error terms are of the same order as in the aforementioned paper. Moreover, it is worth noting that our variational approximation enjoys energy conservation by construction. 

\smallskip 

\noindent Before introducing the method, we mention that the authors of~\cite{carles2025log} deal with the subcritical case~$\alpha > 1$ and, as we have seen in~\eqref{eq:approxlineaire}, they approached the exact solution to~\eqref{eq:logequationalpha} by the linear solution. In~\cite{lasserlubich2020}, this linear solution is approximated using a variational method. This explains why we restrict our analysis to the critical case~$\alpha = 1$. Moreover, as mentioned in \cite{carles2025log}, we focus on times $t\geqslant 0$ since the equation admits a time-reversal symmetry.

\subsection{Presentation of the variational method} \label{subsec:varapp} The purpose of this subsection is to introduce the \textit{variational approximation}. To that end, we consider a smooth submanifold \( \mathcal{M} \) of a Hilbert space \( \mathcal{H} \). For such a point~\( f \in \mathcal{M} \), we consider the tangent space to \( \mathcal{M} \) at \( f \), denoted by \( \mathcal{T}_f \mathcal{M} \) and defined as following
\[
\mathcal{T}_f \mathcal{M} := \left\{ \left. \frac{\diff \gamma}{\diff t}(0) \,\right|\, \gamma : (-\delta, \delta) \to \mathcal{M} \text{ is a differentiable path with } \delta > 0 \text{ and } \gamma(0) = f \right\},
\]
\noindent In other words, \( \mathcal{T}_f \mathcal{M} \) consists of the derivatives of all differentiable paths on~\(\mathcal{M}\) passing through the point~\(f\). The tangent space \( \mathcal{T}_f \mathcal{M} \) is a real-linear closed subspace of \( \mathcal{H} \). In this study, we will consider a submanifold $\mathcal{M}$ such that for all point $f \in \mathcal{M}$, the tangent space $\mathcal{T}_f \mathcal{M}$ is a complex vector space, meaning that if $v \in \mathcal{T}_f \mathcal{M}$, $iv \in \mathcal{T}_f \mathcal{M}$.

\smallskip

\noindent For all $t\in \R$, we seek a good approximation \( \psiapp(t) \in \mathcal{M} \) of the solution \( \psi^{\varepsilon}(t) \) to \eqref{eq:logequation}. To this end, we minimize a suitable functional. Within the framework of the \textit{variational approximation}, this optimization problem is reduced to computing the orthogonal projection of a relevant quantity onto the tangent space of \(\mathcal{M}\) at \(\psiapp(t)\).

\smallskip 

\noindent For simplicity, we omit the dependence on the time variable~\(t\) in the following. We want the function \( \psiapp \) to be a good approximate solution for \eqref{eq:logequation}, which leads us to require
\[i \varepsilon \partial_t \psiapp \approx \mathrm{H}^{\varepsilon} \psiapp + \lambda \,\eps \, \psiapp \, \log \, \lvert \psiapp \lvert^2 \qquad \text{with} \qquad \mathrm{H}^{\varepsilon} := - \dfrac{\varepsilon^2}{2} \Delta + V(x).\]

\noindent By construction, \( \partial_t \, \psiapp \in \mathcal{T}_{\psiapp} \mathcal{M} \) and it follows that \( i \varepsilon \partial_t \psiapp \in \mathcal{T}_{\psiapp} \mathcal{M} \). However, the linear and source terms~\( \mathrm{H}^{\varepsilon} \psiapp + \lambda \, \eps \, \psiapp \, \log \, \lvert \psiapp \lvert^2 \) do not necessarily belong to the tangent space. Therefore, it is natural to project it orthogonally onto the tangent space to obtain the most optimal approximation. Denoting the orthogonal projection onto the tangent space at \( \psiapp \) by ~\( \mathrm{P}_{\psiapp}\), we require the function \(\psiapp\) to satisfy 
\begin{equation}
\label{eq:varapp}
i \varepsilon\partial_t \psiapp = \mathrm{P}_{\psiapp} \left[ \mathrm{H}^{\varepsilon} \psiapp + \lambda \,\eps \, \psiapp \, \log \, \lvert \psiapp \lvert^2 \right] \quad ; \quad \psiapp (t=0) = \psi_{0}^{\varepsilon},
\end{equation}

\noindent with the same initial condition $\psi_{0}^{\varepsilon}$ as in \eqref{eq:logequation}. The following figure illustrates the process.

\begin{figure}[H]
\centering
\begin{tikzpicture}[xshift=1cm]
\def\L{6}
\def\H{2}
\def\s{1.5}
\coordinate (A) at (0,0);
\coordinate (B) at (\L,0);
\coordinate (C) at (\L+\s,\H);
\coordinate (D) at (\s,\H);
\draw[thick] (A)--(B)--(C)--(D)--cycle;
\node at ($(C)+(0.75,0)$) {$\mathcal{T}_{\psiapp} \mathcal{M}$};

\coordinate (P) at (2,0.8);
\fill (P) circle (2pt);
\node[above left] at (P) {$\psiapp$};

\coordinate (V) at (6,3);

\draw[->, thick, orange] (P) -- (V)
node[above right] {$\mathrm{H}^{\varepsilon} \psiapp + \lambda \, \eps \, \psiapp \, \log \lvert \psiapp \lvert^2$};

\coordinate (Proj) at (6, 0.8);
\draw[dashed, purple!85!black] (V) -- (Proj);

\fill[purple!85!black] (Proj) circle (2pt);
\node[right, purple!85!black] at ($(Proj)+(0.75,0)$) 
{$\mathrm{P}_{\psiapp}\!\left[\mathrm{H}^{\varepsilon} \psiapp + \lambda \, \eps \, \psiapp \, \log \lvert \psiapp \lvert^2\right]$};
\end{tikzpicture}
\caption{Orthogonal projection onto the tangent space $\mathcal{T}_{\psiapp}\mathcal{M}$.}
\label{fig:projection_tangent}
\end{figure}

\noindent An equivalent formulation of the equation satisfied by $\psiapp$ in~\eqref{eq:varapp} is given by the following identity: the function~$\psiapp$ is a solution to~\eqref{eq:varapp} if and only if, for all $v \in \mathcal{T}_{\psiapp} \mathcal{M}$, we have
\begin{equation}
    \label{eq:varappbis}
    \big\langle v , i \varepsilon \partial_t \psiapp - \mathrm{H}^{\varepsilon} \psiapp - \nl{\psiapp} \big\rangle = 0.
\end{equation}
Thanks to this reformulation, we can prove that the variational approximation enjoys both norm and energy conservation, as does the exact solution (see Proposition~\ref{prop:existencelog} in Appendix~\ref{appendix:existence}). This is a distinctive feature of this method, stated in the following proposition.

\smallskip

\begin{prop}[Norm and energy conservation by the variational method] \label{prop:normenergyconservation} The variational approximation conserves both the norm and the energy, that is, for all $t \geqslant 0$
\[\normLp{\psiapp(t)}{2}{\R^d} = \normLp{\psiapp(0)}{2}{\R^d} \quad \text{and} \quad E_{\psiapp}^\eps(t) := \big\langle \psiapp(t), \, \mathrm{H}^{\varepsilon}\psiapp(t) + \Big(\log\lvert \psiapp(t)\rvert^2-1\Big)\psiapp(t) \big\rangle = E_{\psiapp}^\eps(0).\]
\end{prop}

\noindent The proof of this proposition is presented in Appendix~\ref{appendix:existence}, Section \ref{sec:conservations}.

\subsection{Variational approximation by complex Gaussians} We consider the manifold of complex-valued Gaussian functions, denoted by $\mathcal{G}$ and defined by
\begin{equation*}
    \mathcal{G} = \Big\lbrace g \in \Lp{2}{\R^d}{\C} \;\big|\; g(x) = \exp \Bigg( \dfrac{i}{\varepsilon} \Big( \dfrac{1}{2} \transpose{(x-q)} C (x-q) \Big) + \transpose{p}(x-q) + \zeta \Bigg), (q,p) \in \R^{2d}, \, \, \zeta \in \C \, \, \text{and} \,\,C \in \siegel{d}{\C} \Big\rbrace
\end{equation*}

\noindent where \(\siegel{d}{\C}\) denotes the Siegel space
\[
\siegel{d}{\C}
= \left\{
C \in \mathbb{C}^{d \times d}
\;\middle|\;
C = C^{\mathsf{T}}, \;\;
\operatorname{Im}(C) \text{ is positive definite}
\right\}.
\] Several useful results about this space will be needed and are stated in Appendix \ref{appendix:siegel}. In the remainder of this work, a Gaussian function \( g \in \mathcal{G}\), parameterized by~\((q,p,\zeta,C) \in \mathbb{R}^{2d} \times \mathbb{C} \times \siegel{d}{\C}\), will be denoted by \(g := G(q,p,\zeta,C),\) where the map \(G\) is defined by
\begin{equation}
\label{eq:fonctionnelG}
    G : (q,p,\zeta,C) \longmapsto \left[ x \mapsto \exp\left( \frac{i}{\varepsilon}\left( \frac{1}{2} (x-q)^{\top} C (x-q) + p^{\top}(x-q) + \zeta \right) \right) \in L^2(\R^d) \right].
\end{equation}

\noindent Following the previous section, we wish to construct an approximate solution $\psiapp \in \mathcal{G}$ satisfying~\eqref{eq:varapp} when the initial condition $\psi_{0}^{\varepsilon}$ at time $t=0$ belongs to $\mathcal{G}$. For this purpose, we need to determine the parameters defining $\psiapp$ so that this function solves \eqref{eq:varapp}. This allows us to provide an explicit characterization of the
approximate solution~$\psiapp$ in terms of its parameters. Before stating the associated results, we introduce the following notation.

\begin{notation}
\label{nota:crochet}
We consider $g\in \mathcal{G}$ such that $\theta := 
\normLp{g}{2}{\R^d}\in (0,1].$ We define for any smooth function $f$
\begin{equation*}
    \left \langle f \right\rangle_g  := \left(\normLp{g}{2}{\R^d}\right)^{-2} \left\langle g \; , \; fg \right\rangle_{L^2(\R^d)} 
     = \dfrac{1}{\theta^2} \ds\int_{\R^d} f(x) |g(x)|^2 \, \diff x.
\end{equation*}
\end{notation}

\noindent Using this notation, we can now state the main theorems obtained to construct the variational approximation~$\psiapp$.

\subsection{Single Gaussian aproximation} The following theorem  characterizes the variational approximation~$\psiapp$ in terms of its parameters.

\begin{theorem}[Parameters of the variational approximation] \label{theo:parameters} Let $0 < \varepsilon \ll 1$ and $\lambda \in \mathbb{R}$. We consider a vector~$(q_0,p_0) \in \R^{2d}$, $\zeta_0 \in \C$ and $C_0 \in \siegel{d}{\C}$. We assume that the potential~$V$ satisfies Assumption~\ref{hyp:potential} and that the initial data $\psi_{0}^{\varepsilon} = G(q_0,p_0,\zeta_0,C_0) \in \mathcal{G}$ has a unit norm. Then, there exists a unique solution~\(\psiapp = G(q,p,\zeta,C) \in \mathcal{C}(\R, \mathcal{G})\) to~\eqref{eq:varapp} and the parameters of the Gaussian $\psiapp$ are the global solutions to the following system
\vspace{-0.2cm}
\begin{equation}
\label{eq:eqpara}
    \left\lbrace 
    \begin{array}{cll}
        \dot{q} & = & p \\
        \dot{p} & = & - \langle \nabla V \rangle_{\psiapp}\\
        \dot{C} & = & - C^2 - \langle \mathrm{Hess} \, V \rangle_{\psiapp}  + 2 \lambda \mathrm{Im}(C)\\
		\dot{\zeta} & = & \dfrac{\lvert p \lvert^2}{2} - \langle V \rangle_{\psiapp}  + \dfrac{i \varepsilon}{2} \mathrm{tr}(C) + \dfrac{\varepsilon}{4} \mathrm{tr}( (\mathrm{Im}(C))^{-1} \langle \mathrm{Hess} \, V \rangle_{\psiapp} ) + 2 \lambda \mathrm{Im}(\zeta)
    \end{array}
    \right.
\end{equation}
\noindent with the initial conditions at time $t=0$ given by \((q(0),p(0)) = (q_0,p_0)\), \(C(0) = C_0\), \(\zeta(0) = \zeta_0.\)
\end{theorem}

\noindent Now, the next theorem provides, for all $t \geqslant 0$, an error bound in the $\mathrm{L}^2$--norm between the approximate Gaussian~\(\psiapp(t) = G(q(t),p(t),\zeta(t),C(t))\), constructed in Theorem~\ref{theo:parameters}, and the solution~$\psi^\eps(t)$ of~\eqref{eq:logequation}.

\begin{theorem}[Error bound in the single Gaussian case] \label{th:errorbound} Let $0 < \varepsilon \ll 1$ and $\lambda \in \R$. We consider a vector~$(q_0,p_0) \in \R^{2d}$, $\zeta_0 \in \C$ and $C_0 \in \siegel{d}{\C}$. We assume that the potential $V$ satisfies Assumption~\ref{hyp:potential}. For $\psi^{\varepsilon}_0 = G(q_0,p_0,\zeta_0,C_0) \in \mathcal{G}$ with a unit norm, we consider~$\psi^{\varepsilon}$ the solution to the Cauchy problem \eqref{eq:logequation} and $\psiapp$ the variational approximation defined in Theorem~\ref{theo:parameters}. Then, there exists $C > 0$, independent of~$\varepsilon$, such that for all $t \geqslant 0 $, we have
\begin{equation*}
    \normLp{\psi^{\varepsilon}(t) - \psiapp(t)}{2}{\R^d} \leqslant C \sqrt{\varepsilon}\textrm{e}^{Ct}.
\end{equation*} 

\noindent In particular, there exists $c>0$, independent of $\varepsilon$, such that 
\begin{equation}
\label{eq:erreurEhrenfest}
    \sup_{|t| \leqslant c \log (\varepsilon^{-1})} \normLp{\psi^{\varepsilon}(t) - \psiapp(t)}{2}{\R^d} \underset{\varepsilon \to 0}{\longrightarrow} 0.
\end{equation}
\end{theorem}

\begin{remark} We can compare the result obtained in this theorem with some previous results.

\smallskip

\noindent The linear case, corresponding to $\lambda = 0$, was studied in~\cite{lasserlubich2020}. When $\lambda = 0$, the first part of Theorem~\ref{th:errorbound} established exactly the same estimate as the one obtained in~\cite{lasserlubich2020}.

\smallskip

\noindent As explained before, the approximation of solutions to \eqref{eq:logequation} using wave packets was studied in~\cite{carles2025log}. In Theorem~\ref{th:errorbound}, we obtain the same order in $\varepsilon$ in the error bound and the same time of validity (an \textit{Ehrenfest time}) for our approximation as those obtained in~\cite{carles2025log}, when the authors consider an initial Gaussian data (see the second part of~\cite[Theorem~$1.3$]{carles2025log}). We note that, in~\cite{carles2025log}, the potential $V$ is assumed to decouple the variables in order to consider a simplified setting and keep the notation as light as possible. In our work, the variational method allows us to avoid this assumption without significantly complicating the computations. A more detailed comparison between the approximation in~\cite{carles2025log} and the one from our paper can be found in Section~\ref{sec:comparaison}.
\end{remark}

\subsection{About superposition} A natural question would be studying the approximation of a solution to~\eqref{eq:logequation} evolved from a sum of initial Gaussian states. The authors of~\cite{carles2025log} prove a superposition principle for such data under the assumption that $V$ is a decoupling potential and that the initial states are centered at different points in phase space. Such a result is far from trivial, as the nonlinear interaction terms generated by the superposition generally do not remain negligible, making their control a central analytical difficulty. Under the previous assumptions, the authors are able to track the evolution of the individual states and construct an approximate solution as a sum of two Gaussian states. Similar phenomena have been studied in the context of other nonlinear Schrödinger equations, see~\cite{carlesfermanian2011} for example. Related results have also been obtained for Hartree equations in~\cite{CaoCarles,Carles12interaction} and for systems in~\cite{Adiab,hari2013coherent}. 

\smallskip

\noindent We are therefore interested in deriving an analogous superposition principle using our method and comparing it with the one obtained in~\cite[Theorem~$1.4$]{carles2025log}. For simplicity, we restrict our analysis to the case of two Gaussian states. The following theorem establishes the desired result and can be generalized to $N$ initial Gaussian states with pairwise distinct centers~$(q_j,p_j)_{j=1}^N$ in phase space.

\smallskip

\begin{theorem}[Error bound for superposed Gaussian initial data] \label{th:superposition} Let $0 < \varepsilon \ll 1$ and $\lambda \in \R$. We consider two initial data~$(q_1^0,p_1^0,\zeta_1^0,C_1^0)$ and $(q_2^0,p_2^0,\zeta_2^0,C_2^0)$ that belong to $\mathbb{R}^{2d} \times \mathbb{C} \times \siegel{d}{\C}$ and we denote $\psi^\varepsilon$ the solution to \eqref{eq:logequation} with~$\psi^\varepsilon_0$ given by \(\psi^\varepsilon_0 = \psi^\eps_{1,0} + \psi^\eps_{2,0}\), where \[\psi^\eps_{1,0} = G(q_1^0,p_1^0,\zeta_1^0,C_1^0), \quad \psi^\eps_{2,0} = G(q_2^0,p_2^0,\zeta_2^0,C_2^0) \in \mathcal{G}\] and satisfy \[\normLp{\psi^\eps_{1,0}}{2}{\R^d} = N_{0,1}, \quad \normLp{\psi^\eps_{2,0}}{2}{\R^d} = N_{0,2} \quad \text{with} \quad N_{0,1}^2+N_{0,2}^2=1.\] We assume that the potential \(V\) satisfies Assumption~\ref{hyp:potential} and that the centers are distincts, that is \[(q_1^0,p_1^0)\neq (q_2^0,p_2^0).\] If we consider~$\psi^\eps_{1,\mathrm{app}}$ (resp.~$\psi^\eps_{2,\mathrm{app}}$) the variational approximation associated with $\psi^\eps_{1,0}$ (resp. $\psi^\eps_{2,0}$) constructed in Theorem~\ref{theo:parameters}, then the following assertions hold.
\begin{enumerate}[leftmargin=*, labelindent=0pt]
\item For any $T > 0$ and $0 < \gamma < \frac{1}{2}$, there exists $C > 0$, independent of $\varepsilon$ such that 
\begin{equation*}
    \sup_{t \in [0,T]} \normLp{\psi^{\varepsilon}(t) - \psi^\eps_{1,\mathrm{app}}(t) - \psi^\eps_{2,\mathrm{app}}(t)}{2}{\R^d} \leqslant C \varepsilon^{\gamma}.
\end{equation*}

\item If $d=1$ and $$E:=\dfrac{(p_1^0)^2}{2}+V(q_1^0)- \dfrac{(p_2^0)^2}{2}-V(q_2^0)\neq 0,$$ then for any $0 < \gamma < \frac{1}{2}$, there exist $c,C > 0$ and $\varepsilon_0 > 0$ such that for all $\varepsilon \in (0,\varepsilon_0)$
\begin{equation*}
    \normLp{\psi^{\varepsilon}(t) - \psi^\eps_{1,\mathrm{app}}(t) - \psi^\eps_{2,\mathrm{app}}(t)}{2}{\R^d} \leqslant C \varepsilon^{\gamma^{-}} e^{Ct}, \quad \text{for } 0 \leqslant t \leqslant c \log (\varepsilon^{-1}).
\end{equation*}
\noindent In particular, there exists~$\kappa > 0$, independent of $\varepsilon$, such that 
\begin{equation*}
    \sup_{|t| \leqslant \kappa \log (\varepsilon^{-1})} \normLp{\psi^{\varepsilon}(t) - \psi^\eps_{1,\mathrm{app}}(t) - \psi^\eps_{2,\mathrm{app}}(t)}{2}{\R^d} \underset{\varepsilon \to 0}{\longrightarrow} 0.
\end{equation*}
\end{enumerate}
\end{theorem}

\noindent In the previous theorem, by~$\varepsilon^{\gamma^{-}}$, we mean $\varepsilon^{\gamma - \beta} $ with $\beta > 0$ sufficiently small.

\subsection{Ideas of the proof and organization of the article} The goal of Section~\ref{sec:parameters} is to construct the variational approximation. To this end, we first derive a condition on the parameters by studying the tangent space and the projection onto it (see Section~\ref{sec:gaussianapptangent}). This condition leads to 
the system of differential equations \eqref{eq:eqpara} that must be satisfied by the parameters (see Subsection~\ref{subsec:eqofpara}). The main remaining task is to establish the global existence of solutions to this system in
a suitable space, ensuring that the associated Gaussian belongs to the manifold G. This relies on the results
of Section ~\ref{subsec:hagedorn}, where we introduce an alternative parametrization of the Gaussian functions, known as the Hagedorn parametrization. The new parameters satisfy a different system of equations, given in~\eqref{eq:eqparasymp}, and
enjoy suitable properties that allow us to obtain required estimates more easily and thereby establishing
their global existence. We then return to the original parametrization in Section~\ref{subsec:globalpara} to conclude.

\smallskip

\noindent In Section~\ref{sec:single}, we prove the validity of the approximation as stated in Theorem~\ref{th:errorbound}.

\smallskip

\noindent Section~\ref{sec:comparaison} provides a detailed comparison between the approximation constructed in this paper and the wave packet approximation developed in~\cite{carles2025log}. In Proposition~\ref{prop:comparaisonpara}, we compare the parameters arising in the two approaches.

\smallskip

\noindent Finally, Section~\ref{sec:superposition} is devoted to the proof of the superposition result.

\medskip

\noindent \textbf{Acknowledgments.} The authors would like to express their sincere gratitude to Clotilde Fermanian Kammerer for her valuable advice, fruitful discussions and constructive suggestions throughout the preparation of this paper. The author would also like to thank Caroline Lasser for her lecture at the workshop \textit{Mathematical Challenges in Quantum Mechanics}, held in L'Aquila in February $2025$. This work is supported by the Région Pays de la Loire via the Connect Talent Project HiFrAn $2022$ $07750$, as well as by the France $2030$ program, Centre Henri Lebesgue ANR-$11$-LABX-$0020$-$01$ and the grants ANR-$25$-CE$40$-$7296$ (La Gabare) and ANR-$23$-CE$40$-$0005$ (Nonlinear Quantum Graphs). 

\section{Variational Gaussian approximation} \label{sec:parameters} In this section, we provide all the necessary details to construct the variational approximation~$\psiapp$ using a complex Gaussian. In what follows, $g$ denotes a Gaussian function belonging to the manifold $\mathcal{G}$. All the theoretical results established for $g$ will then be applied to our approximation $\psiapp$.

\subsection{Tangent space and projection} \label{sec:gaussianapptangent} To analyze the equation satisfied by the variational approximation, we have to compute the orthogonal projection onto a tangent space. We first recall the definition of the tangent space~$\mathcal{T}_{g}\mathcal{G}$ associated with a Gaussian function \( g \in \mathcal{G} \), stated in~\cite[Lemma 3.1]{lasserlubich2020}. For $g \in \mathcal{G}$, the tangent space is given by
\begin{equation*}
    \mathcal{T}_{g}\mathcal{G} = \Big\lbrace \varphi_2 \, g, \, \,  \varphi_2 \, \, \text{a complex $d$-variate polynomial of degree at most} \, \, 2 \Big\rbrace.
\end{equation*}

\noindent Consequently, we deduce that for all $g \in \mathcal{G}$, $\mathcal{T}_{g}\mathcal{G}$ is a complex-linear subspace of $\Lp{2}{\R^d}{\C}$. Moreover, for all differential operators $\mathrm{A}$ of order at most ~$2$ with constant coefficients, we have $\mathrm{A} g \in \mathcal{T}_{g}\mathcal{G}$.

\smallskip

\noindent We now wish to compute the projection
\( \mathrm{P}_g \big[ \mathrm{H}^{\varepsilon} g + \nl{g} \big]\), for $g \in \mathcal{G}$. Since the Laplacian is a second-order operator, it follows that
\begin{equation*}
    \mathrm{P}_g \Big( - \dfrac{\varepsilon^2}{2} \Delta g \Big) = - \dfrac{\varepsilon^2}{2} \Delta g.
\end{equation*}

\noindent It remains to compute the projection of the potential term $V g$ and of the nonlinear term $\nl{g}$. The first computation is carried out in~\cite[Proposition 3.14]{lasserlubich2020} and we recall the result here.

\begin{prop}[Potential term projected onto the tangent space] \label{prop:PuVu} 
For $g \in \mathcal{G}$ and for a potential $V : \R^d \longrightarrow \R$ satisfying Assumption~\ref{hyp:potential}, using Notation~\ref{nota:crochet}, we have 
\[ \mathrm{P}_g\big( V g \big) = \Big( a_0 + \transpose{a_1}(x -q) + \dfrac{1}{2} \transpose{(x-q)} a_2 (x-q) \Big) g \]
\noindent where
\[ a_0 = \langle V \rangle_g - \dfrac{\varepsilon}{4} \mathrm{tr}\big( (\mathrm{Im}(C))^{-1} \langle \mathrm{Hess} \,V \rangle_g \big), \quad a_1 = \langle \nabla \,V \rangle_g \quad \text{and} \quad a_2 = \langle \mathrm{Hess} \,V \rangle_g.\]
\end{prop}

\begin{proof}
We briefly summarize the main steps of the proof, following the proof of~\cite[Proposition~3.14]{lasserlubich2020}. We write
\[V = \left[\langle V \rangle_g + \langle \nabla V \rangle_g^{\mathsf{T}} (x - q) + \frac{1}{2} (x - q)^{\mathsf{T}} \langle \mathrm{Hess}\, V \rangle_g (x - q)\right] +  W.\]
\noindent Using the definition of the tangent space, we have \[ \left[ \langle V \rangle_g + \langle \nabla V \rangle_g^{\mathsf{T}} (x - q) + \frac{1}{2} (x - q)^{\mathsf{T}} \langle \mathrm{Hess}\, V \rangle_g (x - q) \right] g \in \mathcal{T}_g \mathcal{G}\] and it remains to compute \(\mathrm{P}_g W g \). It is easy to check that the function $W$ satisfies the required conditions in Lemma~\ref{lemma:3.13} in Appendix~\ref{app:proj}. This allows us to conclude the proof, see Appendix~\ref{app:proj} for further details.
\end{proof}

\noindent We next consider the projection of the nonlinear term onto the tangent space, which leads to the following proposition.

\begin{prop}[Nonlinear term projected onto the tangent space] \label{prop:Pulogu}
For~\( g \in \mathcal{G}\), the nonlinear term~\( \nl{g} \) belongs to the tangent space \(\mathcal{T}_g \mathcal{G}\). Therefore, we have \[\mathrm{P}_g \big( \nl{g} \big) = - \lambda \Big( \transpose{(x-q)} \mathrm{Im}(C) (x-q) + 2 \mathrm{Im}(\zeta) \Big) g,\]
with \(g = G(q,p,\zeta,C)\) and \((q,p,\zeta,C) \in \R^{2d} \times \C \times \siegel{d}{\C}.\)
\end{prop}

\begin{proof}
We consider a Gaussian function $g \in \mathcal{G}$. Then, there exist \( (q, p) \in \mathbb{R}^{2d} \), \( \zeta \in \mathbb{C} \) and \( C \in \siegel{d}{\C} \) such that
\[ g = \exp \left( \frac{i}{\varepsilon} \left( \frac{1}{2} \transpose{(x - q)} C (x - q \right) + \transpose{p}(x - q) + \zeta \right).\]
Thus, computations give \[\nl{g} = \lambda  \times \Big[ -  \transpose{(x-q)} \mathrm{Im}(C) (x-q) - 2 \mathrm{Im}(\zeta) \Big] g.\]
Noting that \( g \) is multiplied by a complex \( d \)-variate polynomial of degree 2 and using the definition of the tangent space, we conclude the desired result.
\end{proof}

\noindent Consequently, at this stage, we have obtained an explicit expression for the right-hand side of Equation~\eqref{eq:varapp} satisfied by the approximate solution. Therefore, the next subsection establishes the equations satisfied by the parameters of this Gaussian and study the local existence of solutions. In Subsection~\ref{subsec:globalpara}, the global existence is established.

\subsection{Analysis of the parameters of the Gaussian approximation} \label{subsec:eqofpara} The following proposition provides an exact characterization of the parameters of the unique Gaussian solution to \eqref{eq:varapp}. 
\begin{prop}
\label{prop:parasolapp} Let $0 < \varepsilon \ll 1$ and $\lambda \in \mathbb{R}$. We consider $(q_0,p_0) \in \R^{2d}$, $\zeta_0 \in \C$, $C_0 \in \siegel{d}{\C}$ and four smooth functions $(q,p,\zeta,C) : \R \mapsto \mathbb{R}^{2d} \times \mathbb{C} \times \C^{d \times d}$. We assume that the potential $V$ satisfies Assumption~\ref{hyp:potential}. The Gaussian $\psiapp = G(q,p,\zeta,C)$ satisfies \eqref{eq:varapp} with $\psi_0^{\varepsilon} := G(q_0,p_0,\zeta_0,C_0)$ as initial condition at time $t=0$ if and only if the functions~$(q,p,\zeta,C)$ are solution to the system of ordinary differential equations given in~\eqref{eq:eqpara}.
\end{prop}

\begin{proof}
We consider $(q,p,\zeta,C) : \R \mapsto \mathbb{R}^{2d} \times \mathbb{C} \times \C^{d \times d}$, four smooth functions. We denote by~$\psiapp$ the associated Gaussian function, defined for all $(t,x) \in \R \times \R^d$ by \[ \psiapp(t,x) := \exp \left( \frac{i}{\varepsilon} \left( \frac{1}{2} \transpose{(x - q(t))} C(t) (x - q(t)) \right) + \transpose{p(t)}(x - q(t)) + \zeta(t) \right).\] Now, let us compute the partial derivative with respect to time and the Laplacian of the Gaussian function~\((t,x) \mapsto \psiapp(t,x)\). For simplicity, we will omit the time variable in the notation throughout this proof. We have
    \begin{equation*}
       i \varepsilon \partial_t \psiapp = \Big(  \transpose{\dot{q}} C (x-q) - \dfrac{1}{2} \transpose{(x-q)} \dot{C} (x-q) - \transpose{\dot{p}}(x-q) + \transpose{p} \dot{q} - \dot{\zeta} \Big) \psiapp
    \end{equation*}
\noindent and 
    \begin{equation*}
      - \dfrac{\varepsilon^2}{2} \Delta \psiapp = \Big( \dfrac{1}{2} \transpose{(x-q)} C^2 (x-q) +  \transpose{p} C (x-q) + \dfrac{\lvert p \lvert^2}{2} - \dfrac{i \varepsilon}{2} \mathrm{tr}(C) \Big) \psiapp.
    \end{equation*}
We can apply Proposition \ref{prop:PuVu} and it follows that
\begin{equation*}
    \mathrm{P}_{\psiapp} \big( V \psiapp \big) = \Big( \langle V \rangle_{\psiapp} - \dfrac{\varepsilon}{4} \mathrm{tr}\big( (\mathrm{Im}(C))^{-1} \langle \mathrm{Hess} \,V \rangle_{\psiapp} \big) + \langle \transpose{\nabla \,V \rangle_{\psiapp}} (x -q) + \dfrac{1}{2} \transpose{(x-q)} \langle \mathrm{Hess} \,V \rangle_{\psiapp} (x-q) \Big) \psiapp.
\end{equation*}

\noindent Proposition \ref{prop:Pulogu} provides \[\mathrm{P}_{\psiapp} \big(\nl{\psiapp} \big) = - \lambda \Big( \transpose{(x-q)} \mathrm{Im}(C) (x-q) + 2 \mathrm{Im}(\zeta) \Big) \psiapp.\]

\noindent Then, we deduce that $\psiapp$ satisfies
\begin{multline*}
    i \varepsilon \partial_t \psiapp + \dfrac{\varepsilon^2}{2} \Delta \psiapp - \mathrm{P}_{\psiapp} \big( V \psiapp \big) - \mathrm{P}_{\psiapp} \big(\nl{\psiapp} \big) \\ = \transpose{(x-q)} \Big(- \dfrac{1}{2} \dot{C} - \dfrac{1}{2} C^2 -  \dfrac{1}{2} \langle \mathrm{Hess} \,V \rangle_{\psiapp} + \lambda \mathrm{Im}(C) \Big) (x-q) \, \psiapp + \Big(  \transpose{\dot{q}} C - \transpose{\dot{p}} - \transpose{p} C - \langle \transpose{\nabla \,V \rangle_{\psiapp}} \Big) (x -q) \, \psiapp \\ + \Big(  \transpose{p} \dot{q} - \dot{\zeta} - \dfrac{\lvert p \lvert^2}{2} + \dfrac{i \varepsilon}{2} \mathrm{tr}(C) - \langle V \rangle_{\psiapp} + \dfrac{\varepsilon}{4} \mathrm{tr}\big( (\mathrm{Im}(C))^{-1} \langle \mathrm{Hess} \,V \rangle_{\psiapp} \big) + 2 \lambda \mathrm{Im}(\zeta) \Big) \psiapp.
\end{multline*}
By setting the coefficient of the term involving $\Ima (C)$ in the linear part to zero, we obtain $\dot{q} = p$. The
remaining equations for~$p$, $C$ and~$\zeta$ then follow by vanishing the linear, quadratic and constant terms.
\end{proof}

\noindent To deduce the result stated in Theorem~\ref{theo:parameters}, the main difficulty is to prove that the Gaussian obtained with these parameters belongs to the manifold $\mathcal{G}$. In particular, we have to prove that for all $t \in \R$, the matrix~$C(t)$ belongs to the Siegel space. We first prove a local result in the following proposition. The
global existence will be proved in the next subsection.

\begin{prop} \label{prop:existenceparaloc} Let $0 < \varepsilon \ll 1$ and $\lambda \in \mathbb{R}$. We consider $(q_0,p_0) \in \R^{2d}$, $\zeta_0 \in \C$ and $C_0 \in \siegel{d}{\C}$. If the potential $V$ satisfies Assumption \ref{hyp:potential}, then, there exists $T > 0$ such that System~\eqref{eq:eqpara} admits a unique solution $(q,p,C,\zeta)$ on the interval~$[0,T)$, with values in $\R^{2d} \times \siegel{d}{\C} \times \C$.
\end{prop}

\begin{proof} To apply the Cauchy--Lipschitz theorem, we consider the following Banach space \[ \mathrm{E}=\mathbb{R}^{2d}\times\symetrique{d}{\mathbb{C}}\times\mathbb{C} \quad \text{with} \quad \symetrique{d}{\mathbb{C}}=\left\{M\in\mathrm{M}_d(\mathbb{C}) \; ; \; M=\transpose{M}\right\}.\]
We then define \(U=\mathbb{R}^{2d}\times\siegel{d}{\mathbb{C}}\times\mathbb{C}.\) According to Proposition~\ref{prop:Siegelopen}, $U$ is an open subset of $E$ and System~\eqref{eq:eqpara} can be rewritten in the form $\dot{X}=\mathrm{F}(t,X), $ where \( X=\transpose{(q,p,C,\zeta)} \) and the function \(\mathrm{F}\) is defined for all $\left( t, \transpose{(a,b,M,z)} \right) \in \R \times U$ as
\begin{equation*}
    \mathrm{F}\left( t, \begin{pmatrix} a\\ b\\ M\\ z \end{pmatrix} \right) = \begin{pmatrix} b \\ - \valeurmoyenne{\nabla V}{a}{\Ima (M)}{z} \\ - M^2 - \valeurmoyenne{\textrm{Hess} \, V}{a}{\Ima (M)}{z} + 2 \lambda \Ima (M)\\ \dfrac{|b|^2}{2} - \valeurmoyenne{V}{a}{\Ima (M)}{z} + \dfrac{i \varepsilon}{2} \mathrm{tr}(M) \\ + \dfrac{\varepsilon}{4} \mathrm{tr}\big( (\Ima (M))^{-1} \valeurmoyenne{\textrm{Hess} \, V}{a}{\Ima (M)}{z} \big) \\ + 2 \lambda \Ima (z) \end{pmatrix}.
\end{equation*}

\noindent Since the potential \(V\) satisfies Assumptions~\ref{hyp:potential} and \(M\in\siegel{d}{\C}\), the function \(\mathrm{F}\) is smooth on \(\R \times U\). In particular, it is continuous and locally Lipschitz with respect to the second variable. Therefore, the Cauchy--Lipschitz theorem yields the existence of a unique maximal solution \(X:[0,T)\longrightarrow U,\) for some~\(T>0\). This proves the proposition.
\end{proof}

\subsection{Hagedorn’s parametrization} \label{subsec:hagedorn} To deduce the global existence of the solution of \eqref{eq:eqpara}, we want to prove that the solution does not blow up in finite time. The main difficulty is to control the norm of the matrix $C$. To obtain this control, we introduce a different parametrization of complex Gaussians, in which a matrix $C$ of the Siegel space $\siegel{d}{\C}$ is factorized into two matrices with remarkable properties. 

\subsubsection{Introduction of new parameters} Thanks to Proposition~\ref{prop:caracSiegel} in Appendix~\ref{appendix:siegel}, we know that a matrix~$C \in \siegel{d}{\C}$ can be written using two new matrices $Q,P$ belonging to $\matinv{d}{\C}$ as follows \[C = PQ^{-1} \quad \text{with } \quad \Ima (C) = (Q Q^{\ast})^{-1}.\] Consequently, every Gaussian function \(g\in\mathcal{G}\) can be expressed, for every \(x\in\R^d\), as 
\begin{equation*} 
g(x) = \exp\left( \frac{i}{\varepsilon} \left( \frac{1}{2} (x-q)^{\top} PQ^{-1} (x-q) + p^{\top}(x-q) + \zeta \right) \right),
\end{equation*} 
where \((q,p)\in\R^{2d}\), \(\zeta\in\C\) and \((Q,P) \in \matinv{d}{\C}^2\) such that the symplectic relations are satisfied
\begin{equation} 
\tag{\ref{eq:caracsymp}} \left\{ \begin{array}{rcl} Q^{\top}P - P^{\top}Q &=& 0,\\ Q^{\ast}P - P^{\ast}Q &=& 2i\,\mathrm{I}_d. \end{array} \right. \end{equation}

\noindent Moreover, for all $x \in \R^d$, we have
\begin{align}
\label{eq:modulegaussiennehaguedorn}
    |g(x)| = \exp\left( - \frac{1}{2\varepsilon} (x-q)^{\top} (QQ^{\ast})^{-1} (x-q) - \frac{1}{\varepsilon} \Ima (\zeta) \right) .
\end{align}

\noindent This is what is called the \textit{Hagedorn parametrization} of a complex Gaussian \(g\in\mathcal{G}\). In analogy with \eqref{eq:fonctionnelG}, we introduce the notation 
\(g:=\widetilde{G}(q,p,Q,P,\zeta),\) where for all $x \in \R^d$ 
\[\widetilde{G}(q,p,Q,P,\zeta)(x) = \exp\left( \frac{i}{\varepsilon} \left( \frac12 (x-q)^{\top} PQ^{-1} (x-q) +p^{\top}(x-q) +\zeta \right) \right).\] We are led to study a new system of equations governing the evolution of the parameters~\((q,p,Q,P,\zeta)\) of the approximate solution. This is carried out in the following proposition.
\begin{prop}
\label{prop:existenceparameterssymp}
Let $0 < \varepsilon \ll 1$ and $\lambda \in \mathbb{R}$. We consider $(q_0,p_0) \in \R^{2d}$, $\zeta_0 \in \C$ and two matrices $Q_0,P_0$ that belong to $\C^{d \times d}$ and that satisfy the symplecticity condition given in \eqref{eq:caracsymp}. We assume that the potential \(V\) satisfies Assumption~\ref{hyp:potential}. We consider the following system of ordinary differential equations
\begin{equation}
\label{eq:eqparasymp}
    \left\{
    \begin{array}{rcl}
        \dot{q} &=& p,\\
        \dot{p} &=& -\langle \nabla V \rangle_u,\\
        \dot{Q} &=& P,\\
        \dot{P} &=& -\langle \mathrm{Hess}\,V \rangle_u Q
        +2\lambda (Q^\ast)^{-1},\\
        \dot{\zeta}
        &=&
        \dfrac{|p|^2}{2}
        -\langle V\rangle_u
        +\dfrac{i\varepsilon}{2}\,\mathrm{tr}(PQ^{-1})
        +\dfrac{\varepsilon}{4}\,
        \mathrm{tr}\!\left(
        QQ^\ast\langle \mathrm{Hess}\,V\rangle_u
        \right)
        +2\lambda\,\mathrm{Im}(\zeta).
    \end{array}
    \right.
\end{equation}
\noindent with \( u=\widetilde{G}(q,p,Q,P,\zeta)\) and where the initial conditions at \(t=0\) are given by
\[ (q(0),p(0))=(q_0,p_0),\qquad (Q(0),P(0))=(Q_0,P_0),\qquad \zeta(0)=\zeta_0.\]

\noindent Then, there exists $T>0$ such that~\eqref{eq:eqparasymp} admits a unique solution on $[0,T)$. Moreover, for all~$t \in [0,T)$, the matrices~$Q(t)$ and $P(t)$ satisfy the symplecticity condition given in~\eqref{eq:caracsymp}.
\end{prop}

\begin{proof}
As in the proof of Proposition~\ref{prop:existenceparaloc}, we apply the Cauchy--Lipschitz theorem. We consider the Banach space
\(\mathrm{E}=\R^{2d} \times \C^{d \times d} \times \C^{d \times d} \times \C\)
and the open subset \(U\) of \(\mathrm{E}\) defined by
\(U=\R^{2d} \times \matinv{d}{\C} \times \matinv{d}{\C} \times \C.\)
Since the potential~\(V\) satisfies Assumptions~\ref{hyp:potential} and the matrix \((AA^{*})^{-1}\) is positive definite for every matrix~\(A\in\matinv{d}{\C}\), the Cauchy--Lipschitz theorem guarantees the existence of a unique maximal solution
\[
X:[0,T)\longrightarrow U,
\]
for some \(T>0\). It remains to prove that if \(X\) is a solution to \eqref{eq:eqparasymp}, then the matrices~\(Q\) and \(P\) satisfy the symplecticity condition~\eqref{eq:caracsymp}. We compute
\[
\frac{d}{dt}\bigl(Q^{\ast}P - P^{\ast}Q\bigr)
= \dot{Q}^{\ast}P + Q^{\ast}\dot{P} - \dot{P}^{\ast}Q - P^{\ast}\dot{Q}.
\]
\noindent Using that \((Q,P)\) solves \eqref{eq:eqparasymp}, and that \(\mathrm{Hess}\,V\) is Hermitian, we obtain
\begin{align*}
\frac{d}{dt}\bigl(Q^{\ast}P - P^{\ast}Q\bigr)
&= P^{\ast}P
+ Q^{\ast}\Bigl[-\langle \mathrm{Hess}\,V\rangle_u Q + 2\lambda (Q^{\ast})^{-1}\Bigr] - \Bigl[-\langle \mathrm{Hess}\,V\rangle_u Q + 2\lambda (Q^{\ast})^{-1}\Bigr]^{\ast} Q
- P^{\ast}P \\
& = - Q^{\ast} \langle \mathrm{Hess}\,V\rangle_u Q
+ Q^{\ast}\langle \mathrm{Hess}\,V\rangle_u^{\ast} Q + 2\lambda Q^{\ast}(Q^{\ast})^{-1}
- 2\lambda Q^{-1}Q = 0.
\end{align*}

\noindent Hence \(Q^{\ast}P - P^{\ast}Q\) is constant in time and therefore
\[
Q^{\ast}P - P^{\ast}Q
=
Q_0^{\ast}P_0 - P_0^{\ast}Q_0
= 2i\,\mathrm{I}_d,
\]
so that the second symplecticity relation is satisfied. To prove the first relation, we first note that \((QQ^\ast)^{-1}\) is real-valued. Indeed, let us consider the matrix \(M := PQ^{-1}\). We compute
\begin{equation*}
\Ima(M) = \frac{M - M^\ast}{2i} = \frac{PQ^{-1} - (Q^\ast)^{-1}P^\ast}{2i} = (Q^\ast)^{-1} \, \frac{Q^\ast P - P^\ast Q}{2i} \, Q^{-1}.
\end{equation*}

\noindent Using the second symplecticity relation, we deduce that \( \Ima(M) = (QQ^\ast)^{-1}.\) Thus, the matrix~\((QQ^\ast)^{-1}\) is real-valued. We then compute 
\[
\frac{d}{dt}\bigl(\transpose{Q} P - \transpose{P}Q\bigr)
= \transpose{\dot{Q}} P + \transpose{Q} \dot{P} - \transpose{\dot{P}}Q - \transpose{P} \dot{Q}.
\]

\noindent Using that \((Q,P)\) solves \eqref{eq:eqparasymp} and the fact that the matrices \(\mathrm{Hess}\,V\) and \((QQ^\ast)^{-1}\) are real-valued, we obtain that~\(\transpose{Q}P - \transpose{P}Q\) is constant in time and therefore
\[ \transpose{Q}P - \transpose{P}Q = \transpose{Q}_0 P_0 - \transpose{P}_0 Q_0 = 0.\] 
\end{proof}

\noindent Let us now establish the global existence of the solution to \eqref{eq:eqparasymp}. Proposition \ref{prop:controlparameterssymp}, stated below, provides a bound on the norm of the solution and allows us to conclude about the global existence for the parameters~$(q,p,Q,P,\zeta)$. Before proving this proposition, we first establish several preliminary and useful lemmas on Gaussian functions. This is the goal of the following subsection.

\subsubsection{Some properties of Gaussian functions} First, we compute the \(L^r\)--norm of a Gaussian function for~$r \geqslant 1$. In the following, we will apply this result to our approximate solution $\psiapp$.

\begin{lemma}[$L^{r}$--norm of a Gaussian] \label{lemma:normLpgaussienne} Let $0 < \varepsilon \ll 1$ and $\lambda \in \mathbb{R}$. We introduce~$(q_0,p_0) \in \R^{2d}$, $\zeta_0 \in \C$ and two matrices~$Q_0,P_0 \in \C^{d \times d}$ that satisfy the symplecticity condition given in \eqref{eq:caracsymp}. We consider~$T > 0$ and the solution~$(q,p,Q,P,\zeta)$ to \eqref{eq:eqparasymp} obtained in Proposition \ref{prop:existenceparameterssymp}. For~\( g = \widetilde{G}(q,p,Q,P,\zeta)\), for $r \geqslant 1$ and for all $t \in [0,T)$
\begin{equation*} 
N_r(t) := \normLp{g(t)}{r}{\R^d} = \Bigg( \dfrac{2 \pi \varepsilon}{r} \Bigg)^{\frac{d}{2r}} \exp\left(-\frac{\Ima (\zeta(t))}{\varepsilon}\right) \Bigg(\det\Big[ Q(t)Q^{\ast}(t)\Big]\Bigg)^{\frac{1}{2r}}.
\end{equation*}
\end{lemma}

\begin{proof}
Let \(0 < \varepsilon \ll 1\) and \(\lambda \in \mathbb{R}\). We fix \(T>0\) and consider the solution \((q,p,Q,P,\zeta)\) of \eqref{eq:eqparasymp} given by Proposition~\ref{prop:existenceparameterssymp}. We define \(g\) as the Gaussian associated with these parameters: \(g = \widetilde{G}(q,p,Q,P,\zeta).\) For~$r \geqslant 1$ and for $t \in [0,T)$, we denote by $N_r(t)$ the $L^r$--norm of $g(t)$. Using~\eqref{eq:modulegaussiennehaguedorn}, we have
\begin{align*}
    N_r^r(t) & = \exp\left(-\frac{r}{\varepsilon} \Ima (\zeta(t)) \right) \times \int_{\mathbb{R}^d} \exp\left(-\frac{r}{2 \varepsilon} \Big((x-q(t))^{\top} \big(Q(t)Q^{\ast}(t)\big)^{-1} (x-q(t))\Big) \right)\,\diff x\\
    & = \Bigg( \dfrac{2\varepsilon}{r} \Bigg)^{\frac{d}{2}} \exp\left(-\frac{r}{\varepsilon} \Ima (\zeta(t)) \right) \times \int_{\mathbb{R}^d} \exp\left(- y^{\top} \big(Q(t)Q^{\ast}(t)\big)^{-1} y \right)\,\diff y.
\end{align*}

\noindent using the change of variables $y = \sqrt{\dfrac{r}{2\varepsilon}}\,(x-q)$. Since the matrix~ $(QQ^{\ast})^{-1}$ is symmetric definite positive, we can write
\begin{equation}
    \label{eq:diagoQ}
    (QQ^{\ast})^{-1} = \transpose{U} D U
\end{equation}

\noindent where $U$ is an orthogonal matrix and $D = \mathrm{diag}(\lambda_1, \cdots, \lambda_d)$ with $\lambda_1, \cdots, \lambda_d > 0$ the eigenvalues of the matrix $(QQ^{\ast})^{-1}$. Therefore, using the change of variables \(z = U y\), we obtain that for all $t \in [0,T)$
\begin{align*}
    N_r^r(t) & = \Bigg( \dfrac{2\varepsilon}{r} \Bigg)^{\frac{d}{2}} \exp\left(-\frac{r}{\varepsilon} \Ima (\zeta(t)) \right) \times \int_{\mathbb{R}^d} \exp\left(- \sum\limits_{i=1}^d \lambda_i(t) z_i^2 \right) \, \diff z\\
    & = \Bigg( \dfrac{2\varepsilon}{r} \Bigg)^{\frac{d}{2}} \exp\left(-\frac{r}{\varepsilon} \Ima (\zeta(t)) \right) \times \prod\limits_{i=1}^d \int_{\mathbb{R}} \exp\left(- \lambda_i(t) z_i^2 \right) \, \diff z_i.
\end{align*}

\noindent \noindent Using the fact that for all $a > 0$, 
\begin{equation}
\label{eq:integralegaussienne}
    \ds\int_{\R} \exp \left( - a x^2 \right) \diff x = \sqrt{\dfrac{\pi}{a}},
\end{equation}

\noindent we deduce
\begin{equation*}
     N_r^r(t) = \Bigg( \dfrac{2\varepsilon}{r} \Bigg)^{\frac{d}{2}} \exp\left(-\frac{r}{\varepsilon} \Ima (\zeta(t)) \right) \times \pi^{\frac{d}{2}} \times \sqrt{\prod\limits_{i=1}^d \dfrac{1}{\lambda_i(t)}}.
\end{equation*}

\noindent Since $\dfrac{1}{\lambda_1}, \ldots, \dfrac{1}{\lambda_d}$ are the eigenvalues of the matrix $QQ^{\ast}$, we obtain the result.
\end{proof}

\begin{remark} \label{rq:Nr}
    The \(L^2\)-norm is a particularly important case, as it is conserved (see Proposition~\ref{prop:normenergyconservation}). The proof, presented in Appendix~\ref{appendix:existence}, relies only on the definition of the variational approximation and does not require its explicit expression in terms of the parameters. Alternatively, this conservation can be established using the expression derived in Proposition~\ref{lemma:normLpgaussienne} together with \eqref{eq:eqpara} satisfied by the parameters.
\end{remark}

\noindent Using Proposition~\ref{lemma:normLpgaussienne}, we immediately obtain the following corollary, which expresses the \(L^r\)-norm of a Gaussian function solely in terms of \(\varepsilon\) and the matrix~\(Q\).

\begin{corollary} \label{lemma:normgaussienne} Under the assumptions of Lemma~\ref{lemma:normLpgaussienne}, if $N_2(0) = \normLp{g(0)}{2}{\R^d} = N_0 \in \R$, for all $t \in [0,T)$, we have
\begin{equation}
    \label{eq:relationzeta}
    \exp\left(-\frac{\Ima (\zeta(t))}{\varepsilon}\right)
    = N_0 \, (\pi \varepsilon)^{-\frac{d}{4}}
    \Big(\deter{Q(t)Q^{\ast}(t)}\Big)^{-\frac{1}{4}}.
\end{equation}
\noindent In particular, 
\begin{equation}
        \label{eq:Nr}
        N_r(t) = N_0 \, \Big( \dfrac{2}{r} \Big)^{\frac{d}{2r}} (\pi \varepsilon)^{\frac{d}{2r} - \frac{d}{4}} \Big(\deter{Q(t)Q^{\ast}(t)}\Big)^{\frac{1}{2r} - \frac{1}{4}}.
    \end{equation}
\end{corollary}

\noindent The next lemma provides estimates for some momenta of a Gaussian function.

\begin{lemma}\label{lem:gaussian_moments}
Let \(k\in\mathbb{N}\). Under the assumptions of Lemma~\ref{lemma:normLpgaussienne} and for $N_0 = \normLp{g(0)}{2}{\R^d} \leqslant 1$, there exists a constant~\(C_{k} >0\) (independent of \(\varepsilon\)) such that, for every \(t\in[0,T)\), we have
\[\Big\| \|x-q(t) \|^k g(t) \Big\|_{L^2(\R^d)}
\leqslant C_k \, \varepsilon^{\frac{k}{2}} \Big( \smalleigenvalue{(Q(t)Q^{\ast}(t))^{-1}} \Big)^{\frac{k}{2} - \frac{d}{4}} \Big(\deter{Q(t)Q^{\ast}(t)} \Big)^{-\frac{1}{4}},\]
where $\smalleigenvalue{(Q(t)Q^{\ast}(t))^{-1}}$ is the smallest eigenvalue of $(Q(t)Q^{\ast}(t))^{-1}$.
\end{lemma}

\begin{proof}
Let $0 < \varepsilon \ll 1$, $\lambda \in \mathbb{R}$ and $k \in \N$. We consider \(T > 0 \) and the solution~$(q,p,Q,P,\zeta)$ to \eqref{eq:eqparasymp} obtained in Proposition \ref{prop:existenceparameterssymp}. Using \eqref{eq:modulegaussiennehaguedorn} to compute $|g|$ and \eqref{eq:diagoQ} to diagonalize the matrix $(QQ^{\ast})^{-1}$ with~$U$ an orthogonal matrix and setting $y = \dfrac{1}{\sqrt{\varepsilon}} U^{-1}(x-q)$, for all~$t \in [0,T)$, we obtain 
\begin{equation}
    \label{eq:expressionmoments}
    \Big\| \|x-q(t) \|^k g(t) \Big\|^2_{L^2(\R^d)} = \varepsilon^{k +\frac{d}{2}} \exp\left(- \frac{2 \Ima (\zeta(t))}{\varepsilon} \right) \times \ds\int_{\R^d} \| z \|^{2k} \exp\big( - \sum\limits_{i=1}^{d} \lambda_i(t) y_i^2 \big) \diff y.
\end{equation}

\noindent Thanks to Relation~\eqref{eq:relationzeta}, we deduce
\begin{equation*}
    \Big\| \|x-q(t) \|^k g(t) \Big\|^2_{L^2(\R^d)} = N_0^2 \varepsilon^{k} \pi^{-\frac{d}{2}} \times \Big(\deter{Q(t)Q^{\ast}(t)} \Big)^{-\frac{1}{2}} \times \ds\int_{\R^d} \| y \|^{2k} \exp\Big( - \sum\limits_{i=1}^{d} \lambda_i(t) y_i^2 \Big) \diff y,
\end{equation*}
\noindent where $N_0 = \normLp{g(0)}{2}{\R^d}$. Using definitions of the eigenvalues $\lambda_1, \cdots, \lambda_d$ and of the quantity~$\smalleigenvalue{(QQ^{\ast})^{-1}}$, we know \[ \lambda_i \geqslant \smalleigenvalue{(QQ^{\ast})^{-1}}, \quad \forall i \in [\![1,d]\!].\] 

\noindent Then, setting $z = \sqrt{\smalleigenvalue{(QQ^{\ast})^{-1}}} \, y$, we obtain

\begin{multline*}
    \Big\| \|x-q(t) \|^k g(t) \Big\|^2_{L^2(\R^d)} \leqslant N_0^2 \varepsilon^{k} \pi^{-\frac{d}{2}} \Big(\deter{Q(t)Q^{\ast}(t)} \Big)^{-\frac{1}{2}} \times \Big(\smalleigenvalue{(Q(t)Q^{\ast}(t))^{-1}}\Big)^{k - \frac{d}{2}} \ds\int_{\R^d} \| z \|^{2k} e^{- \|z\|^2} \diff z.
\end{multline*}

\noindent Since the function \(z \mapsto \|z\|^{2k} e^{-\|z\|^2}\) is integrable on $\mathbb{R}^d$, the integral is a constant depending only on $\alpha$. Taking the square root and using the assumptions on $N_0$, we obtain the desired result. 
\end{proof}

\noindent We will prove later (see the proof of Proposition~\ref{prop:existenceparaglob}) that the parameters $(p,q,P,Q,\zeta)$ are in fact defined on $\mathbb{R}_+$. Hence, all the computations performed in this subsection for $t \in [0,T)$ remain valid for all $t \geqslant 0$ and we will use them hereafter.

\subsubsection{Global existence of the Hagedorn's parameters} We are now in a position to establish the global existence of the parameters $(q,p,Q,P,\zeta)$ which solve \eqref{eq:eqparasymp}.

\smallskip 

\noindent For convenience, we introduce the Frobenius norm $\normfrob{\cdot}$ and the spectral norm $\normspec{\cdot}$ on $\C^{d\times d}$, defined by 
\begin{equation*}
    \normfrob{M} = \sqrt{\mathrm{tr}(MM^\ast)}
    \quad \text{and} \quad
    \normspec{M} = \sqrt{\largeeigenvalue{MM^\ast}}, \qquad \text{for } M \in \C^{d\times d},
\end{equation*}
where $\largeeigenvalue{MM^\ast}$ denotes the largest eigenvalue of the matrix $MM^\ast$. Since all norms are equivalent on~$\C^{d\times d}$, we use these norms interchangeably throughout the rest of the paper and choose the one that is most convenient for each estimate. Nevertheless, the estimates remain valid for any norm on $\C^{d\times d}$, which we generically denote by $\normtriple{\cdot}$.

\smallskip

\noindent As a first step, we prove a preliminary lemma providing bounds for the average terms \[ \langle V \rangle_{g(t)}, \quad \langle \nabla V \rangle_{g(t)} \quad \text{and} \quad \langle \mathrm{Hess}\,V \rangle_{g(t)}\] appearing in \eqref{eq:eqparasymp}. These estimates rely in particular on Assumption~\ref{hyp:potential} and will play a crucial role in the proof of Proposition~\ref{prop:controlparameterssymp}, which establishes a control on the norm of the parameters $(q,p,Q,P,\zeta)$.

\begin{lemma}
\label{lem:controlmoyenne} Under the same assumptions as in Lemma \ref{lemma:normLpgaussienne} and for $N_0 = \normLp{g(0)}{2}{\R^d} \leqslant 1$, there exists~$C > 0$ such that for all $t \in [0,T)$, we have
\begin{equation*}
    | \langle V \rangle_{g(t)} | \leqslant C \Big( 1 + \|q(t)\|^2 + \normtriple{Q(t)}^2 \Big) , \quad \| \langle \nabla V \rangle_{g(t)}  \| \leqslant C \Big( 1 + \|q(t)\| + \normtriple{Q(t)} \Big), \quad \normtriple{\langle \mathrm{Hess}\,V \rangle_{g(t)} } \leqslant C.
\end{equation*}
\end{lemma}

\begin{proof}
We consider \(g = \widetilde{G}(q,p,Q,P,\zeta).\) Using Assumption \ref{hyp:potential}, there exists $C > 0$ such that for $x \in \R^d$, we have \[ |V(x)| \leqslant C(1 + \|x\|^2), \quad \|\nabla V(x) \| \leqslant C(1 + \|x\|) \quad \text{and} \quad \normtriple{\mathrm{Hess}\, V(x)} \leqslant C.\]

\noindent For the Hessian term, the desired control follows directly from Corollary \ref{lemma:normgaussienne} together with the fact that~\(N_0 \leqslant 1\). In fact, for all $t \in [0,T)$, we can write
\begin{equation*}
    \normtriple{\langle \mathrm{Hess}\,V \rangle_g(t)} = \Bigg|\Bigg|\Bigg| \ds\int_{\R^d} \mathrm{Hess}\, V(x) \, |g(t,x)|^2 \diff x \, \Bigg|\Bigg|\Bigg| \leqslant C N_2^2(t) = C N_2^2(0) \leqslant C_0,
\end{equation*}

\noindent with $C_0 > 0$ independent of $\varepsilon$. For the potential and gradient terms, for all $t \in [0,T)$, it holds that
\begin{align}
\label{eq:expressionmoyenne}
| \langle V \rangle_{g(t)} | & = \Bigg| \ds\int_{\R^d} V(x) \, |g(t,x)|^2 \diff x \Bigg|
\leqslant C N_2^2(t) + C \ds\int_{\R^d} \|x\|^2 \, |g(t,x)|^2 \diff x
\end{align}
\noindent and 
\begin{align}
\label{eq:expressionmoyennebis}
\| \langle \nabla V \rangle_{g(t)} \| & = \Bigg\| \ds\int_{\R^d} \nabla V(x) \, |g(t,x)|^2 dx \Bigg\|
\leqslant C N_2^2(t) + C \ds\int_{\R^d} \| x \| \, |g(t,x)|^2 \diff x.
\end{align}

\noindent Using the Cauchy--Schwarz inequality and the $L^2$--norm conservation of Corollary \ref{lemma:normgaussienne} in the gradient term, for all~$t \in [0,T)$, we deduce that

\begin{align*}
    \ds\int_{\R^d} \| x \| \, |g(t,x)|^2 \diff x & \leqslant \Bigg( \ds\int_{\R^d} |g(t,x)|^2 \diff x \Bigg)^{\frac{1}{2}} \Bigg( \ds\int_{\R^d} \| x \|^2 \, |g(t,x)|^2 \diff x \Bigg)^{\frac{1}{2}}\\
    & \leqslant N_2(0) \times \Bigg( \ds\int_{\R^d} \| x \|^2 \, |g(t,x)|^2 \diff x \Bigg)^{\frac{1}{2}}.
\end{align*}

\noindent Thus, to conclude, it remains to estimate the integral \[I(t) := \ds\int_{\R^d} \| x \|^2 \, |g(t,x)|^2 \diff x. \]

\noindent Using the fact that $\|a+b\|^2 \leqslant 2(\|a\|^2 + \|b\|^2)$, for all $t \in [0,T)$, we have
\begin{align*}
I(t) & = \ds\int_{\R^d} \| x - q + q \|^2 \, |g(t,x)|^2 \diff x\\
& \leqslant \ds\int_{\R^d} \| x - q(t) \|^2 \, |g(t,x)|^2 \diff x + \ds\int_{\R^d} \| q(t) \|^2 \, |g(t,x)|^2 \diff x.
\end{align*}

\noindent Now, let us estimate the first integral. For all $t \in [0,T)$, using the proof of Lemma~\ref{lem:gaussian_moments} and more precisely~\eqref{eq:expressionmoments} for~$\alpha =1$, it yields
\begin{align*}
    \ds\int_{\R^d} \| x - q(t) \|^2 \, |g(t,x)|^2 \diff x & = \varepsilon^{1 + \frac{d}{2}} \exp \left( - \dfrac{2 \Ima (\zeta(t))}{\varepsilon} \right) \ds\int_{\R^d}  \sum\limits_{i=1}^d z_i^2 \exp \left( - \sum\limits_{j=1}^d \lambda_j(t) z_j^2  \right) \diff z \\
    & = \varepsilon^{1 + \frac{d}{2}} \times \exp \left( - \dfrac{2 \Ima (\zeta(t))}{\varepsilon} \right) \times  \sum\limits_{i=1}^d \ds\int_{\R} z_i^2 e^{- \lambda_i(t) z_i^2} \diff z_i \times \prod\limits_{j \neq i} \ds\int_{\R} \exp \left( - \lambda_j(t) z_j^2 \right) \diff z_j.
\end{align*}

\noindent Using the explicit computation given in~\eqref{eq:integralegaussienne} and also the fact that for all $a > 0$, \[\ds\int_{\R} x^2 \exp \left( - a x^2 \right) \diff x = \dfrac{1}{2a} \sqrt{\dfrac{\pi}{a}},\]

\noindent we deduce that
\begin{equation*}
    \ds\int_{\R^d} \| x - q(t) \|^2 \, |g(t,x)|^2 \diff x = \dfrac{\varepsilon}{2} \, \mathrm{tr}\big[ Q(t)Q^{\ast}(t) \big] (\pi \varepsilon)^{\frac{d}{2}} \exp \Big(-\frac{2\Ima (\zeta(t))}{\varepsilon} \Big) \sqrt{\mathrm{det} [Q(t)Q^{\ast}(t)]}.
\end{equation*}

\noindent Thanks to Lemma~\ref{lemma:normLpgaussienne} with $r = 2$, we obtain that for all $t \in [0,T)$
\begin{equation*}
    \ds\int_{\R^d} \| x - q(t) \|^2 \, |g(t,x)|^2 \diff x = \dfrac{\varepsilon}{2} \, \mathrm{tr}\Big[ Q(t)Q^{\ast}(t) \Big] \times N_2^2(t).
\end{equation*}

\noindent Thus, using the conservation of the $L^2$--norm, we conclude that
\begin{equation*}
    I(t) =  N_2^2(0) \Big( \dfrac{\varepsilon}{2} \normfrob{Q(t)}^2 + \|q(t)\|^2 \Big).
\end{equation*}

\noindent Then, thanks to \eqref{eq:expressionmoyenne} and \eqref{eq:expressionmoyennebis}, the conservation of the norm and $\varepsilon \ll 1$, we conclude that
\begin{equation*}
| \langle V \rangle_{g(t)} | \leqslant C N_2^2(0) \Bigg[ 1 + \dfrac{\varepsilon}{2} \normfrob{Q(t)}^2 + \|q(t)\|^2 \Bigg] \leqslant C_1 ( 1 + \normfrob{Q(t)}^2 + \|q(t)\|^2)
\end{equation*}

\noindent and
\begin{equation*}
\| \langle \nabla V \rangle_{g(t)} \| \leqslant C N_2^2(0) \Bigg[ 1 +  \sqrt{ \dfrac{\varepsilon}{2} \normfrob{Q(t)}^2 + \|q(t)\|^2} \Bigg] \leqslant C_2(1 + \|q(t)\| + \normfrob{Q(t)}),
\end{equation*}

\noindent with $C_1, C_2 > 0$, independent of $\varepsilon$. The desired estimate then follows from the equivalence of norms.
\end{proof}

\noindent We are now in a position to establish the control of the Hagedorn's parameters of our Gaussians and the global existence of solutions to~\eqref{eq:eqparasymp}.

\begin{prop} \label{prop:controlparameterssymp} Let $0 < \varepsilon \ll 1$ and $\lambda \in \mathbb{R}$. We consider $(q_0,p_0) \in \R^{2d}$, $\zeta_0 \in \C$ and $Q_0,P_0 \in \C^{d \times d}$ that satisfy the symplecticity condition given in \eqref{eq:caracsymp}. There exist two constants $C_1, C_2 > 0$ such that the solution $(q,p,Q,P,\zeta)$ of \eqref{eq:eqparasymp} obtained in Proposition \ref{prop:existenceparameterssymp} satifies 
\begin{equation}
    \label{eq:controlpara}
    \|q(t) \| + \|p(t) \| + \normtriple{Q(t)} + \normtriple{Q^{-1}(t)}+ \normtriple{P(t)} \leqslant C_1 e^{C_2 t}, \quad \text{for all } t \in [0,T).
\end{equation}
\end{prop}

\begin{proof}
Let $0 < \varepsilon \ll 1$ and $\lambda \in \mathbb{R}$. We consider \(T > 0 \) and the solution $(q,p,Q,P,\zeta)$ to~\eqref{eq:eqparasymp} obtained in Proposition~\ref{prop:existenceparameterssymp}. For \(g = \widetilde{G}(q,p,Q,P,\zeta).\) For \(t \in [0,T)\), we define \[ E_1(t):=\|q(t)\|^2+\|p(t)\|^2, \quad E_2(t):=\normfrob{Q(t)}^2+\normfrob{P(t)}^2 \quad \text{and} \quad E(t):=E_1(t)+E_2(t). \] 

\noindent Our goal is to control the quantity \(E\). We begin by computing the derivatives of \(E_1\) and \(E_2\). Using \eqref{eq:eqparasymp}, for~\(t\in[0,T)\), we obtain that
\begin{equation*} 
|\dot{E_1}(t)| = \Big| 2\,\Rea\Big( \dot{q}(t)\cdot q(t) + \dot{p}(t)\cdot p(t) \Big) \Big| \leqslant 2\|p(t)\|\,\|q(t)\| + 2\|\langle \nabla V\rangle_{g(t)}\|\,\|p(t)\|
\end{equation*} and 
\begin{align*} 
|\dot{E_2}(t)| &= \Big| \mathrm{tr}\!\Big( \dot{Q^{\ast}}(t)Q(t) + Q^\ast(t)\dot{Q}(t) + \dot{P^{\ast}}(t)P(t) + P^\ast(t)\dot{P}(t) \Big) \Big| = \Big| 2\,\Rea\Big( \mathrm{tr}\!\big(Q^\ast(t)\dot{Q}(t)\big) + \mathrm{tr}\!\big(P^\ast(t)\dot{P}(t)\big) \Big) \Big| \\ &\leqslant C\Big( \normfrob{Q(t)}\,\normfrob{P(t)} + \normfrob{P(t)} \,\normfrob{\langle\mathrm{Hess}\,V\rangle_g} \,\normfrob{Q(t)} + \normfrob{P(t)} \,\normfrob{(Q^\ast(t))^{-1}} \Big) 
\end{align*} 
\noindent for some constant \(C>0\). By Lemma~\ref{lem:controlmoyenne}, we deduce that \[|\dot{E_1}(t)| \leqslant 2 \big( \|p(t)\|\,\|q(t)\| + C\bigl(1+\|q(t)\|+\normfrob{Q(t)}\bigr)\|p(t)\|\big)\]
\noindent and
\[|\dot{E_2}(t)| \leqslant C\Big( \normfrob{Q(t)}\,\normfrob{P(t)} + \normfrob{P(t)}\,\normfrob{(Q^\ast(t))^{-1}} \Big).\]

\noindent Using the inequalities \(2ab\leqslant a^2+b^2\) and \(a\leqslant \dfrac{1+a^2}{2}\), we obtain 
\[ |\dot{E_1}(t)| \leqslant C\bigl( 1+\|q(t)\|^2+\|p(t)\|^2+\normfrob{Q(t)}^2 \bigr) \quad \text{and} \quad |\dot{E_2}(t)| \leqslant C\bigl( \normfrob{Q(t)}^2 + \normfrob{P(t)}^2 + \normfrob{(Q^\ast(t))^{-1}}^2 \bigr). \]

\noindent Therefore, for all $t \in [0,T)$, we get
\begin{align} \label{eq:controlenergie} |\dot{E}(t)| \leqslant |\dot{E_1}(t)|+|\dot{E_2}(t)| \leqslant C\bigl(1+E(t)\bigr) + C\,\normfrob{(Q^\ast(t))^{-1}}^2. 
\end{align}

\noindent To conclude the proof, it remains to estimate \(\normfrob{(Q^\ast(t))^{-1}}^2\). To this end, we will use of the symplecticity relations~\eqref{eq:caracsymp}. Indeed, for every \(u\in\R^d\), the second symplecticity relation yields \[ Pu \cdot Qu - Qu \cdot Pu = 2i\,\|u\|^2. \] \noindent Let \(v\in\R^d\setminus\{0\}\). Applying the above identity with \(u=Q^{-1}v\), we obtain \[ 2i\,\|Q^{-1}v\|^2 = PQ^{-1}v\cdot v - v\cdot PQ^{-1}v. \] \noindent Taking the modulus and using the Cauchy--Schwarz inequality, we deduce that \begin{align*} \|Q^{-1}v\|^2 &\leqslant \|PQ^{-1}v\|\,\|v\| \leqslant \normtriple{P}\, \|Q^{-1}v\|\, \|v\|. \end{align*} \noindent Since \(v\neq 0\) and by the definition of the induced matrix norm, we deduce
\begin{equation}
\label{eq:controlQ-1}
    \normtriple{Q^{-1}} \leqslant \normtriple{P}.
\end{equation}

\noindent Since \((Q^\ast)^{-1}=(Q^{-1})^\ast\) and \(\normtriple{A^\ast}=\normtriple{A}\) for all $A \in \C^{d\times d}$, we obtain
\begin{equation}
    \label{eq:controlQ-1star}
    \normtriple{(Q^\ast)^{-1}} \leqslant \normtriple{P}.
\end{equation}

\noindent Returning to~\eqref{eq:controlenergie}, we deduce that \begin{equation*} \dot{E}(t) \leqslant |\dot{E}(t)| \leqslant C\,E(t). \end{equation*} \noindent Grönwall's lemma yields the exponential estimate stated in the proposition. The conclusion follows from the equivalence of norms on \(\C^{d\times d}\), which completes the proof.
\end{proof}

\begin{corollary} \label{prop:controlparameterszeta} Under the same assumptions as in Proposition \ref{prop:existenceparameterssymp}, there exist $C_1, C_2 > 0$ such that the function $\zeta : [0,T) \mapsto \C$, solution to \eqref{eq:eqparasymp}, satisfy for all~$t \in [0,T)$
\begin{equation*}
    |\zeta(t) | \lesssim C_1 e^{C_2 t}.
\end{equation*}
\end{corollary}

\begin{proof}
According to \eqref{eq:eqparasymp}, for every \(t\in[0,T)\), we have
\[\dot{\zeta}(t) = \frac{|p(t)|^2}{2} -\langle V\rangle_{g(t)} +\frac{i\varepsilon}{2}\,\mathrm{tr}\, \big(P(t)Q^{-1}(t)\big) +\frac{\varepsilon}{4}\, \mathrm{tr}\, \left( Q(t)Q^\ast(t)\langle \mathrm{Hess}\,V\rangle_{g(t)} \right) +2\lambda\,\Ima (\zeta(t)). \] 

\noindent We estimate each term appearing in \(\dot{\zeta}\). By Proposition~\ref{prop:controlparameterssymp}, there exist two constants~\(C_1,C_2>0\) such that, for all~\(t\in[0,T)\) \[ \frac{|p(t)|^2}{2} \leqslant C_1e^{C_2t}, \quad \left| \mathrm{tr}\!\big(P(t)Q^{-1}(t)\big) \right| \leqslant C_1e^{C_2t} \quad \text{and} \quad \left| \mathrm{tr}\!\left( Q(t)Q^\ast(t) \langle\mathrm{Hess}\,V\rangle_{g(t)} \right) \right| \leqslant C_1 \normtriple{\langle\mathrm{Hess}\,V\rangle_{g(t)}} e^{C_2t}. \] 

\noindent Using Lemma~\ref{lem:controlmoyenne} to control the term $\langle V\rangle_{g(t)}$ and Proposition~\ref{prop:controlparameterssymp}, there exist~\(C,\widetilde C>0\) such that for~\(t\in[0,T)\), we have \[ |\dot{\zeta}(t)| \leqslant \widetilde C \bigl( |\Ima(\zeta(t))| + e^{Ct} \bigr). \] \noindent Since \(|\Ima(\zeta)| \leqslant |\zeta|\) and \(\frac{\mathrm d}{\mathrm dt}|\zeta| \leqslant |\dot{\zeta}|,\) we obtain \[ \frac{\mathrm d}{\mathrm dt}|\zeta(t)| \leqslant \widetilde C \bigl( |\zeta(t)| + e^{Ct} \bigr), \qquad t\in[0,T). \]

\noindent Therefore, Grönwall's lemma yields the estimate stated in the proposition.
\end{proof}

\noindent The following corollary provides several useful bounds in the remainder of the article.

\begin{corollary}\label{coro:controldetetlambdamin}
    Under the same assumptions as in Proposition \ref{prop:existenceparameterssymp}, there exists a constant $C > 0$ such that for all $t \geqslant 0$, we have
\[\smalleigenvalue{(Q(t)Q^{\ast}(t))^{-1}} \leqslant C e^{Ct}, \quad \deter{Q(t)Q^{\ast}(t)} \leqslant C e^{Ct} \quad \text{and} \quad \left( \deter{Q(t)Q^{\ast}(t)} \right)^{-1} \leqslant C e^{Ct},\]
where $\smalleigenvalue{(QQ^{\ast})^{-1}}$ is the smallest eigenvalue of $\Ima \; C = (QQ^{\ast})^{-1}$.
\end{corollary}

\begin{proof}
Using the definition of the spectral norm, for all $t \geqslant 0$, we have that
\[\smalleigenvalue{(Q(t)Q^{\ast}(t))^{-1}} = \dfrac{1}{\normspec{Q(t)}^2}.\]

\noindent By applying the fact that $1 \leqslant \normspec{Q} \cdot \normspec{Q^{-1}}$, the control given in Proposition \ref{prop:controlparameterssymp} on the norm of $Q^{-1}$ and the equivalence of norm, we obtain the desired bound for $\smalleigenvalue{(QQ^{\ast})^{-1}}$. Regarding to the term involving the determinant, by definition~$\lambda_1, \cdots, \lambda_d$ are the eigenvalues of the matrix $(QQ^{\ast})^{-1}$. Then, we have
\[\deter{QQ^{\ast}} = \prod_{i=1}^d \dfrac{1}{\lambda_i} \quad \text{and} \quad \left( \deter{QQ^{\ast}} \right)^{-1} = \prod_{i=1}^d \lambda_i.\]

\noindent But, for all $i \in [\![1,d]\!]$, we have
\begin{equation*}
    | \lambda_i | \leqslant \normtriple{(QQ^{\ast})^{-1}} \leqslant \normtriple{(Q^{\ast})^{-1}} \cdot \normtriple{Q^{-1}} \quad \text{and} \quad \dfrac{1}{\lambda_i}\leqslant \dfrac{1}{\smalleigenvalue{(QQ^{\ast})^{-1}}} = \normspec{Q}^2.
\end{equation*} 

\noindent Thus, combining this two inequalities~\eqref{eq:controlQ-1} and~\eqref{eq:controlQ-1star} with the equivalence of norms and Proposition~\ref{prop:controlparameterssymp}, we obtain the desired result.
\end{proof}

\subsection{Global existence of the Gaussian parameters} \label{subsec:globalpara} The goal of this subsection is to prove the global existence of the parameters $(q,p,C,\zeta)$ solving \eqref{eq:eqpara}. We recall that their local existence was established in Proposition~\ref{prop:existenceparaloc}. The following proposition establishes their global existence by combining the global existence of the Hagedorn's parameters, obtained thanks to Proposition~\ref{prop:controlparameterssymp}, with the relation between~\eqref{eq:eqpara} for the initial parameters and~\eqref{eq:eqparasymp} for the Hagedorn's parameters.

\begin{prop}
\label{prop:existenceparaglob} Let \(0 < \varepsilon \ll 1\) and \(\lambda \in \mathbb{R}\). We consider \((q_0,p_0)\in\mathbb{R}^{2d}\), \(\zeta_0\in\mathbb{C}\) and~\(C_0\in\siegel{d}{\mathbb{C}}\). If~\(V\) satisfies Assumption~\ref{hyp:potential}, \eqref{eq:eqpara} admits a unique global solution \((q,p,C,\zeta) : [0,+\infty) \longrightarrow \mathbb{R}^{2d}\times\siegel{d}{\mathbb{C}}\times\mathbb{C}\). Moreover, there exist $C_1, C_2 > 0$ such that for all $t \in [0,+\infty)$, we have
\[\|q(t)\|+\|p(t)\|+\|C(t)\|+|\zeta(t)|\leqslant C_1 e^{C_2 t}.\]
\end{prop}

\begin{proof}
Let \(0 < \varepsilon \ll 1\) and \(\lambda \in \mathbb{R}\). By Proposition \ref{prop:controlparameterssymp} and Corollary \ref{prop:controlparameterszeta}, there exist~\(C_1, C_2 >0\) such that
\[
\|q(t)\|+\|p(t)\|+\|Q(t)\|+\|P(t)\|+|\zeta(t)|
\leqslant C_1 e^{C_2 t},
\]
for all $t$ in the maximal interval of existence $[0,T)$. Therefore, the solution cannot blow up in finite time. Hence, the standard continuation theorem for ordinary differential equations implies that the maximal solution can be extended beyond any finite time. Consequently, the solution is global, that is, its maximal interval of existence is \([0,+\infty)\). We denote \((q,p,Q,P,\zeta) : [0,+\infty) \to \mathbb{R}^{2d} \times \matinv{d}{\mathbb{C}}^2 \times \mathbb{C}\) this global solution. The matrices \(Q\) and \(P\) satisfy the symplecticity relations given in \eqref{eq:caracsymp}. We set \(M = P Q^{-1}\). According to Proposition \ref{prop:caracSiegel} in Appendix \ref{appendix:siegel}, we know that~$M$ belongs to $\siegel{d}{\C}$ and satisifies the relation \(\Ima (M) = (QQ^\ast)^{-1}. \) Moreover, since $Q$ and $P$ are solutions to \eqref{eq:eqparasymp}, for all $t \in [0,T)$, we deduce
\begin{align*}
    \dot{M}(t) & = \dot{P}(t) Q^{-1}(t) - P Q^{-1}(t) \dot{Q}(t) Q^{-1}(t) \\
            & = \Big[ - \langle \textrm{Hess} \, V \rangle_{g(t)} Q(t) + 2\lambda (Q(t)Q^\ast(t))^{-1} Q(t) \Big] Q^{-1}(t) - P(t) Q^{-1}(t) P(t) Q^{-1}(t) \\
            & = - \langle \textrm{Hess} \, V \rangle_{g(t)} + 2\lambda (Q(t)Q^\ast(t))^{-1} - M(t)^2 \\
            & = - \langle \textrm{Hess} \, V \rangle_{g(t)} + 2\lambda \Ima (M(t)) - M(t)^2,
\end{align*}

\noindent where \(g\) denotes the Gaussian associated with these parameters~\( g = \widetilde{G}(q,p,Q,P,\zeta).\) Thus, $M = PQ^{-1}$ satisfies the fourth equation in \eqref{eq:eqpara}. Proposition \ref{prop:existenceparaloc} allows us to obtain the local existence of solution for the system \eqref{eq:eqpara}: there exists $T > 0$ such that $(q,p,C,\zeta)$ exist uniquely on the time interval $[0,T)$. Consequently, using the unicity of the solution, the matrix $C$ coincide with the matrix $M = PQ^{-1}$ where $Q$ and $P$ are solutions to \eqref{eq:eqparasymp} on $[0,T)$. According to Proposition \ref{prop:controlparameterssymp}, we deduce that the norm of $C$ does not explose on finite time. For all $t \in [0,T)$, we have
\begin{equation*}
    \normtriple{C(t)} \leqslant \normtriple{P(t)} \cdot \normtriple{Q^{-1}(t)} \leqslant C_1 e^{C_2t}.
\end{equation*}
\noindent Then, the norms of all the parameters $(q,p,C,\zeta)$ remain bounded on every finite time interval (the equation on $q, p$ and $\zeta$ are the same than before). Therefore, the solution of \eqref{eq:eqpara} is global. Moreover, using the fact that $Q$ and $P$ satisfy the symplecticty relation (see Proposition \ref{prop:existenceparameterssymp}), by Proposition~\ref{prop:caracSiegel}, it yields that
\[C(t) \in \siegel{d}{\C}, \quad \forall t \in [0,+\infty).\]

\noindent The control of the solution follows from Lemma~\ref{prop:controlparameterssymp}.
\end{proof}

\begin{remark} \label{rq:controlImC} According to the preceding proof, we have seen that there exist~$Q,P \in \matinv{d}{\C}$, solving~\eqref{eq:eqparasymp}, such that \(\Ima(C)=(QQ^{\ast})^{-1}.\) Using Proposition~\ref{prop:controlparameterssymp} and more specifically the two inequalities~\eqref{eq:controlQ-1} and~\eqref{eq:controlQ-1star} established in its proof, we also obtain an exponential bound on the norm
of $\Ima(C)$: there exist two constants $C_1,C_2>0$ such that for all $t\geqslant0$
\[\normtriple{\Ima(C)}\leqslant C_1 e^{C_2t}.\]
\end{remark}

\noindent We can now prove Theorem~\ref{theo:parameters}, which gives an explicit expression for the approximate solution $\psiapp$ constructed using the variational approximation.

\begin{proof}
Let $0 < \varepsilon \ll 1$ and $\lambda \in \mathbb{R}$. We consider $(q_0,p_0) \in \R^{2d}$, $\zeta_0 \in \C$ and $C_0 \in \siegel{d}{\C}$ and we assume that the potential $V$ satisfies Assumption \ref{hyp:potential}. We set $\psi^{\varepsilon}_0 = G(q_0,p_0,\zeta_0,C_0) \in \mathcal{G}$. To obtain the result stated in Theorem \ref{theo:parameters}, we have to use Proposition \ref{prop:parasolapp}, the fact that the solution of \eqref{eq:eqparasymp} is global and that the matrix $C(t)$ belongs to the Siegel space $\siegel{d}{\C}$ for all~$t \in [0,+\infty)$. Indeed, according to Proposition~\ref{prop:existenceparaglob}, System~\eqref{eq:eqpara} has a unique global solution denoted by \[(q,p,C,\zeta) : \R \longmapsto \mathbb{R}^{2d}\times\siegel{d}{\mathbb{C}}\times\mathbb{C}.\]

\noindent For all $t \in [0,+\infty)$, we define \(\psiapp(t) := g(q(t),p(t),C(t),\zeta(t)) \in \mathcal{G}.\) According to Proposition \ref{prop:parasolapp}, we know that \( \psiapp\) is the unique Gaussian that satisfies the following equation
\begin{equation}
\tag{\ref{eq:varapp}}
    i \varepsilon\partial_t \psiapp = \mathrm{P}_{\psiapp} \left[ \mathrm{H}^{\varepsilon} \psiapp + \lambda \, \psiapp \, \log \, \lvert \psiapp \lvert^2 \right]
\end{equation}
\noindent with $\psiapp(t=0) = \psi^{\varepsilon}_0$, which conclude the proof.
\end{proof}

\section{Single Gaussian approximation: proof of the main result} \label{sec:single} In this section, we prove the error estimate stated in Theorem \ref{th:errorbound}. In this case, the initial data is a single Gaussian state of $\mathcal{G}$ and we compare the exact solution of \eqref{eq:logequation} with the variational approximation, denoted by~$\psiapp$, constructed in the previous section (see Theorem \ref{theo:parameters}).

\vspace{-0.2cm}
\begin{proof}
The beginning of the proof is the same for both parts of the theorem. We consider $0 < \varepsilon \ll 1$, $\lambda \in \R$, a vector $(q_0,p_0) \in \R^{2d}$, $\zeta_0 \in \C$ and $C_0 \in \siegel{d}{\C}$. We assume that the potential~$V$ satisfies Assumption \ref{hyp:potential}.  We consider $\psi^{\varepsilon}$ the solution to the Cauchy problem \eqref{eq:logequation} with the initial data given by~$\psi^{\varepsilon}_0 = G(q_0,p_0,\zeta_0,C_0) \in \mathcal{G}$ and $\psiapp$ the variational approximation defined in Theorem \ref{theo:parameters}. We introduce the error quantity defined by $e^{\varepsilon} = \psi^{\varepsilon} - \psiapp$. We know that the functions $\psi^{\varepsilon}$ and $\psiapp$ satisfy the following equations
     \[ i \varepsilon \partial_t \psi^{\varepsilon} = \mathrm{H}^{\varepsilon} \psi^{\varepsilon} + \nl{\psi^{\varepsilon}} \quad \text{and} \quad i \varepsilon \partial_t \psiapp =\mathrm{P}_{\psiapp} \mathrm{H}^{\varepsilon} \psiapp + \nl{\psiapp} \] with the same initial data at time $0$. We deduce that the error $e^{\varepsilon}$ is solution of the Cauchy problem 
\begin{equation*}
\label{eq:cauchypberreur}
i \varepsilon \partial_t e^{\varepsilon} = \mathrm{H}^{\varepsilon} e^{\varepsilon} + \mathrm{P}_{\psiapp}^{\perp} \, \mathrm{H}^{\varepsilon} \psiapp + \lambda \, \eps \, \Big( (e^{\varepsilon}+\psiapp) \log \, \lvert e^{\varepsilon}+\psiapp \lvert^2 -  \, \psiapp \log \, \lvert \psiapp \lvert^2 \Big) \quad ; \quad e^{\varepsilon} (0) = 0,
\end{equation*}
where $\mathrm{P}_{\psiapp}^{\perp} = \mathrm{I_d} - \mathrm{P}_{\psiapp}$ is the projection onto the orthogonal complement of $\mathcal{T}_{\psiapp} \mathcal{G}$. We compute
\begin{align*}
    \dfrac{\diff}{\diff t} \normLp{e^{\varepsilon}}{2}{\R^d}^2 & = 2 \mathrm{Re}\Big( \langle \partial_t e^{\varepsilon} , e^{\varepsilon} \rangle \Big) \\ & = 2 \mathrm{Re}\Bigg( \Big\langle \dfrac{1}{i\varepsilon} \mathrm{H}^{\varepsilon} e^{\varepsilon} + \dfrac{1}{i \varepsilon} \mathrm{P}_{\psiapp}^{\perp} \, \mathrm{H}^{\varepsilon} \psiapp + \dfrac{\lambda}{i} \Big( (e^{\varepsilon}+\psiapp) \log \, \lvert e^{\varepsilon}+\psiapp \lvert^2 - \psiapp \log \, \lvert \psiapp \lvert^2 \Big) , e^{\varepsilon} \Big\rangle \Bigg).
\end{align*}
Since $\mathrm{H}^{\varepsilon}$ is self-adjoint, $\langle \mathrm{H}^{\varepsilon} e^{\varepsilon} , e^{\varepsilon} \rangle$ is real and therefore
\begin{equation*}
    \dfrac{\diff}{\diff t} \normLp{e^{\varepsilon}}{2}{\R^d}^2 = 2 \mathrm{Re}\Bigg( \big\langle \dfrac{1}{i \varepsilon} \mathrm{P}_{\psiapp}^{\perp} \, \mathrm{H}^{\varepsilon} \psiapp , e^{\varepsilon} \big\rangle \Bigg) + 2\lambda \, \mathrm{Im}\Bigg( \big\langle (e^{\varepsilon}+\psiapp) \log \, \lvert e^{\varepsilon}+\psiapp \lvert^2 - \, \psiapp \log \, \lvert \psiapp \lvert^2 , e^{\varepsilon} \big\rangle \Bigg).
\end{equation*}

\noindent Using the fact that for all $z \in \C$, $\mathrm{Re}(z) \leqslant \lvert z \lvert$ and the Cauchy-Schwarz inequality , we obtain 
\[ 2  \mathrm{Re}\Bigg( \Big\langle \dfrac{1}{i \varepsilon} \mathrm{P}_{\psiapp}^{\perp} \, \mathrm{H}^{\varepsilon} \psiapp , e^{\varepsilon} \Big\rangle \Bigg) \leqslant 2 \Bigg\| \dfrac{1}{i \varepsilon} \mathrm{P}_{\psiapp}^{\perp} \, \mathrm{H}^{\varepsilon} \psiapp \Bigg\|_{L^2(\R^d)} \times \normLp{e^{\varepsilon}}{2}{\R^d}. \]
\noindent We note that
\[ \Bigg\| \dfrac{1}{i \varepsilon} \mathrm{P}_{\psiapp}^{\perp} \, \mathrm{H}^{\varepsilon} \psiapp \Bigg\|_{L^2(\R^d)} = \, \, \mathrm{dist} \Big( \dfrac{1}{i \varepsilon} \mathrm{H}^{\varepsilon} \psiapp , \mathcal{T}_{\psiapp} \mathcal{M} \Big) \]
\noindent and then 
\[ 2  \mathrm{Re}\Bigg( \Big\langle \dfrac{1}{i \varepsilon} \mathrm{P}_{\psiapp}^{\perp} \, \mathrm{H}^{\varepsilon} \psiapp , e^{\varepsilon} \Big\rangle \Bigg) \leqslant 2 \, \mathrm{dist} \Big( \dfrac{1}{i \varepsilon} \mathrm{H}^{\varepsilon} \psiapp , \mathcal{T}_{\psiapp} \mathcal{M} \Big) \, \normLp{e^{\varepsilon}}{2}{\R^d}. \]
For the second term that appears in $\dfrac{\diff}{\diff t} \normLp{e^{\varepsilon}}{2}{\R^d}^2$, we use Equation \eqref{lem:majlog} with $z_1 = \psiapp$ and $z_2 = \psi^{\varepsilon}$. We obtain
\begin{equation*} 
2\lambda \, \mathrm{Im}\Bigg( \big\langle (e^{\varepsilon}+\psiapp) \log \, \lvert e^{\varepsilon}+\psiapp \lvert^2 - \, \psiapp \log \, \lvert \psiapp \lvert^2 , e^{\varepsilon} \big\rangle \Bigg) \leqslant 4 | \lambda | \normLp{e^{\varepsilon}}{2}{\R^d}^2.
\end{equation*}

\noindent Consequently, we get
\[ \dfrac{\diff}{\diff t} \normLp{e^{\varepsilon}}{2}{\R^d}^2 \leqslant 2 \, \mathrm{dist} \Big( \dfrac{1}{i \varepsilon} \mathrm{H}^{\varepsilon} \psiapp , \mathcal{T}_{\psiapp} \mathcal{M} \Big) \, \normLp{e^{\varepsilon}}{2}{\R^d} + 4 | \lambda |   \normLp{e^{\varepsilon}}{2}{\R^d}^2. \]

\noindent Using $\dfrac{\diff}{\diff t} \normLp{e^{\varepsilon}}{2}{\R^d}^2 = 2 \normLp{e^{\varepsilon}}{2}{\R^d} \dfrac{\diff}{\diff t} \normLp{e^{\varepsilon}}{2}{\R^d}$, we deduce that
\[ \dfrac{\diff}{\diff t} \normLp{e^{\varepsilon}}{2}{\R^d} \leqslant \mathrm{dist} \Big( \dfrac{1}{i \varepsilon} \mathrm{H}^{\varepsilon} {\psiapp} , \mathcal{T}_{\psiapp} \mathcal{M} \Big) + 2 | \lambda |  \, \normLp{e^{\varepsilon}}{2}{\R^d}.\]

\noindent Integrating between $0$ and $t$ for $t \geqslant 0$, we have
\begin{equation}
\label{eq:majoerreurpreuve}
    \normLp{e^{\varepsilon}}{2}{\R^d} \leqslant \ds\int_{0}^{t} \mathrm{dist} \Big( \dfrac{1}{i \varepsilon} \mathrm{H}^{\varepsilon} {\psiapp}(s) , \mathcal{T}_{{\psiapp}(s)} \mathcal{M} \Big) \, \diff s + \ds\int_{0}^{t} 2 | \lambda | \, \normLp{e^{\varepsilon}}{2}{\R^d} \, \diff s.
\end{equation}

\noindent Using the notations of \cite{lasserlubich2020}, we write the second order Taylor expansion of $V$ at $q$ as \[V=U_q+W_q,\] where $U_q$ is a polynomial of degree 2 and $W_q$ is the corresponding remainder. We notice that the terms~$\Delta \psiapp$ and~$U_q \psiapp$ belong to $\mathcal{T}_{\psiapp} \mathcal{M}$. Hence,
\begin{align*}
    \mathrm{dist} \Big( \dfrac{1}{i \varepsilon} \mathrm{H}^{\varepsilon} {\psiapp}(s) , \mathcal{T}_{{\psiapp}(s)} \mathcal{M} \Big) & = \mathrm{dist} \Big( \dfrac{1}{i \varepsilon} W_q \;  {\psiapp}(s) , \mathcal{T}_{{\psiapp}(s)} \mathcal{M} \Big) \\ & \leqslant \dfrac{1}{\eps} \left\| W_q {\psiapp}(s) \right\|_{L^2(\R^d)}.
\end{align*}

\noindent Since $W_q$ involves a polynomial of degree 3, using Lemma \ref{lem:gaussian_moments} with $\alpha = 3$, we end up with
\[\left\| W_q {\psiapp}(s) \right\|_{L^2(\R^d)} \leqslant  C \varepsilon^{\frac{3}{2}} \times \left( \deter{Q(s)Q^{\ast}(s)} \right)^{-\frac{1}{4}} \left(\smalleigenvalue{(Q(s)Q^{\ast}(s))^{-1}}\right)^{\frac{3}{2}-\frac{d}{4}}\]

\noindent where $\smalleigenvalue{(QQ^{\ast})^{-1}}$ is the smallest eigenvalue of $\Ima \; C = (QQ^{\ast})^{-1}$. Consequently, we deduce that
\begin{multline*}
    \int_0^t \mathrm{dist} \Big( \dfrac{1}{i \varepsilon} \mathrm{H}^{\varepsilon} {\psiapp}(s) , \mathcal{T}_{{\psiapp}(s)} \mathcal{M} \Big)\leqslant \dfrac{C \varepsilon^{\frac{3}{2}} }{\eps} \int_0^t \left( \deter{Q(s)Q^{\ast}(s)} \right)^{-\frac{1}{4}} \times \left(\smalleigenvalue{(Q(s)Q^{\ast}(s))^{-1}}\right)^{\frac{3}{2}-\frac{d}{4}} \, \diff s.
\end{multline*}

\noindent According to Corollary~\ref{coro:controldetetlambdamin}, we get 
\[ \smalleigenvalue{(Q(s)Q^{\ast}(s))^{-1}} \leqslant C e^{Cs} \quad \text{and} \quad \left( \deter{Q(s)Q^{\ast}(s)} \right)^{-\frac{1}{4}} \leqslant C e^{Cs}.\]

\noindent Therefore, we have the following control
\begin{equation}
\label{eq:majointdistance}
\int_0^t \mathrm{dist} \Big( \dfrac{1}{i \varepsilon} \mathrm{H}^{\varepsilon} {\psiapp}(s) , \mathcal{T}_{{\psiapp}(s)} \mathcal{M} \Big) \leqslant C \sqrt{\varepsilon} \int_0^t e^{Cs}\, \diff s = C \sqrt{\varepsilon} (e^{Ct} - 1) \notag \leqslant C \sqrt{\varepsilon} e^{Ct}
\end{equation}

\noindent where the constant can be changed at each line. We now distinguish between the two cases of the theorem. Then, substituting the estimate established in~\eqref{eq:majointdistance} into~\eqref{eq:majoerreurpreuve} for $t \geqslant 0$, it infers that
\[ \normLp{e^{\varepsilon}}{2}{\R^d} \leqslant C \sqrt{\varepsilon} e^{Ct} + \ds\int_{0}^{t} 2 | \lambda | \, \normLp{e^{\varepsilon}}{2}{\R^d} \, \diff s.\]
\noindent Since $t \mapsto C \sqrt{\varepsilon} e^{Ct}$ is a non-decreasing function, the Grönwall lemma yields
\[ \normLp{e^{\varepsilon}}{2}{\R^d} \leqslant C \sqrt{\varepsilon}\,\exp\!\big(C t\big) \exp\left(2|\lambda|t\right) \leqslant \widetilde{C}\sqrt{\eps}\exp(\widetilde{C}t)\]
\noindent with $\widetilde{C} > 0$. To obtain \eqref{eq:erreurEhrenfest} stated in Theorem~\ref{th:errorbound}, it suffices to choose the constant~$c$ appropriately, depending on~$\widetilde{C}$. This completes the proof of the theorem.
\end{proof}

\section{Comparison with wave packets approximation} \label{sec:comparaison} This section is devoted to the comparison between the wave packet approximation constructed in~\cite{carles2025log} and the one obtained in this article with the variational approach.

\subsection{The wave packet approximation} We first recall how the wave packet approximation is constructed in~\cite{carles2025log} adapting the notations to our setting.

\begin{theorem}[From Theorem $1.3$ in \cite{carles2025log}] \label{th:WPapp} Let $0 < \varepsilon \ll 1$ and $\lambda \in \R$. We assume that the potential $V$ satisfies Assumption~\ref{hyp:potential} and decouples the variables, in the sense that
\vspace{-0.2cm}
\begin{equation}
\label{eq:Vdecouples}
    V(x) = \sum\limits_{i=1}^d V_i(x_i), 
\end{equation}

\noindent for all $x = (x_1,\ldots,x_d) \in \R^d$, where each $V_i : \R \longrightarrow \R$ is a smooth function. We suppose that $\psi_0^{\varepsilon}$ is a Gaussian function of unit norm, defined for all $x \in \R^d$ by
\[\psi_0^{\varepsilon}(x) = \varepsilon^{-\frac{d}{4}} b_0 \exp \Bigg( \dfrac{i}{\varepsilon} \transpose{p_0}(x - q_0) - \dfrac{1}{2\varepsilon} \transpose{(x-q_0)} D_0 (x-q_0) \Big) \Bigg),\]
where $(q_0,p_0) \in \R^{2d}$, $b_0 \in \C$ and $D_0 = \mathrm{diag}(d_{0,1}, ..., d_{0,d}) \in \C^{d\times d}$ satisfies \(\Rea (d_{0,i}) > 0 \) for all \(i \in [\![ 1,d]\!].\) Then, there exists~$C > 0$, independent of $\varepsilon$, such that
\[\normLp{\psi^{\varepsilon}(t) - \psi_{_{\mathrm{WP}}}^{\varepsilon}(t)}{2}{\R^d} \leqslant C \sqrt{\varepsilon} e^{Ct}, \quad \text{for all } t \geqslant 0,\]

\noindent where $\psi_{_{\mathrm{WP}}}^{\varepsilon}$ is a Gaussian wave packet approximation, defined for all $t \geqslant 0$ and all $x \in \R^d$ by
\begin{multline*}
    \psi_{_{\mathrm{WP}}}^{\varepsilon}(t,x) = \varepsilon^{-\frac{d}{4}} b(t) \exp \Bigg( \dfrac{i}{\varepsilon} \Big( \transpose{p_{_{\mathrm{WP}}}(t)}(x - q_{_{\mathrm{WP}}}(t)) + S(t) \Big) - \dfrac{1}{2\varepsilon} \transpose{(x-q_{_{\mathrm{WP}}}(t))} D(t) (x-q_{_{\mathrm{WP}}}(t)) \Bigg)
\end{multline*}

\noindent with $(q_{_{\mathrm{WP}}}(t),p_{_{\mathrm{WP}}}(t)) \in \R^{2d}$, $b(t) \in \C$, $S(t)\in \R$ and the diagonal matrix $D(t) = \mathrm{diag}(d_{1}(t), ..., d_{d}(t)) \in \C^{d\times d}$ satisfies \( \Rea (d_{i}(t)) > 0,\) for all~\(i \in [\![ 1,d]\!].\) Moreover, the parameters $(q_{_{\mathrm{WP}}},p_{_{\mathrm{WP}}},S,b,D)$ are solutions to the following system
\begin{equation}
\label{eq:eqparaWP}
    \left\lbrace 
    \begin{array}{cll}
        \dot{q}_{_{\mathrm{WP}}} & = & p_{_{\mathrm{WP}}} \\
        \dot{p}_{_{\mathrm{WP}}} & = & - \nabla V (q_{_{\mathrm{WP}}})\\
        \dot{S} & = & \dfrac{|p_{_{\mathrm{WP}}}|^2}{2} - V(q_{_{\mathrm{WP}}})\\
        i \dot{b} & = & \dfrac{1}{2} \sum\limits_{i=1}^d d_i b + \lambda \, b \, \log |b|^2\\
        i \dot{d_i} & = & d_i^2 - V_i''(q_{_{\mathrm{WP}}}) + 2 \Rea (d_i), \quad \forall i \in [\![ 1,d ]\!]
    \end{array}
    \right.
\end{equation}
\noindent where initial conditions at time $t=0$ are given by \((q_{_{\mathrm{WP}}}(0),p_{_{\mathrm{WP}}}(0)) = (q_0,p_0), \, S(0) = 0, \, b(0) = b_0\),and \(D(0) = D_0.\)
\end{theorem}

\noindent We observe that the error estimate between the exact solution $\psi^{\varepsilon}$ of~\eqref{eq:logequation} and the wave packet approximation $\psiWP$ is similar to the one obtained in our theorem (Theorem~\ref{th:errorbound}). Moreover, the exponential bounds for the parameters of wave packets evolving under an at most quadratic potential (such as those established for our parameters in Proposition~\ref{prop:existenceparaglob}) are already known. For the center $(q_{_{\mathrm{WP}}},p_{_{\mathrm{WP}}})$, a proof can be found in~\cite[Chapter 10]{carleslivre2021semi}, while the corresponding estimates for the parameters $(b,D)$ are established in~\cite{carles2025log} (see Lemma~$4.1$ and the proof of Lemma~$4.2$). We summarize these results in the following proposition.

\begin{prop} \label{prop:controlparaWP} Let $0 < \varepsilon \ll 1$ and $\lambda \in \R$. We assume that $V$ satisfies Assumption~\ref{hyp:potential} and decouples the variables as in~\eqref{eq:Vdecouples}. We consider $(q_0,p_0) \in \R^{2d}$, $b_0 \in \C$ and a matrix~$D_0 = \mathrm{diag}(d_{0,1}, ..., d_{0,d}) \in \C^{d\times d}$ with
\(\Rea (d_{0,i}) > 0\) for all \(i \in [\![ 1,d ]\!].\) Then, System~\eqref{eq:eqparaWP} has a unique global smooth solution $(q_{_{\mathrm{WP}}},p_{_{\mathrm{WP}}},S,b,D)$ with values in $\R^{2d + 1} \times \C \times \C^{d \times d}$, where $D = \mathrm{diag}(d_{1}, ..., d_{d})$ satisfies
\(\Rea (d_{i}) > 0\) for all \(i \in [\![ 1,d ]\!].\) Moreover, these parameters grow at most exponentially: there exists $C > 0$ such that, for all $t \geqslant 0$, we have
\vspace{-0.4cm}
\begin{equation*}
    \|q_{_{\mathrm{WP}}}(t) \| + \| p_{_{\mathrm{WP}}}(t) \| + |b(t)| + \sum\limits_{i=1}^d |d_i(t)| \leqslant C e^{Ct}.
\end{equation*}
\end{prop}

\noindent The main difference between the wave packet approximation and the variational approximation lies in the way the parameters of the approximations are defined. The system satisfied by the parameters $(q,p,C,\zeta)$ of $\psiapp$, given in~\eqref{eq:eqpara}, differs from the one satisfied by the parameters $(q_{_{\mathrm{WP}}},p_{_{\mathrm{WP}}},S,b,D)$ of $\psiWP$, given in~\eqref{eq:eqparaWP}. The goal of the next subsection is to compare these different sets of parameters.

\subsection{Comparison of the parameters} First, we state a useful lemma that provides the first-order asymptotics of some averaged quantities appearing in the equations satisfied by the variational parameters~$(q,p,C,\zeta)$ and those of the wave packet approximation~$(q_{_{\mathrm{WP}}},p_{_{\mathrm{WP}}},S,b,D)$.

\begin{lemma} \label{lem:asympaverage} Let $0 < \varepsilon \ll 1$ and $\lambda \in \R$. We assume that the potential $V$ satisfies Assumptions~\ref{hyp:potential} and decouples the variable as in~\eqref{eq:Vdecouples}. We consider $\psiapp$ the variational approximation constructed in Theorem~\ref{theo:parameters}. Then, there exists $C >0$ independent of $\varepsilon$ such that for all $t \geqslant 0$ and for $i\in[\![1,d]\!]$, the following asymptotics hold.
\begin{equation*} 
\langle V_i'\rangle_{\psi_{\mathrm{app}}(t)} = V_i'(q_i(t)) + R^{\varepsilon}_{i,1}(t), \quad \langle V_i''\rangle_{\psi_{\mathrm{app}}(t)} = V_i''(q_i(t)) + R^{\varepsilon}_{i,2}(t) \quad \text{and} \quad \langle V\rangle_{\psi_{\mathrm{app}}(t)} = V(q(t)) + R^{\varepsilon}_0(t),
\end{equation*}
where the remainder terms $R^{\varepsilon}_0$, $R^{\varepsilon}_{i,1}$ and $R^{\varepsilon}_{i,2}$ satisfy
\begin{equation*}
    | R^{\varepsilon}_{i,1}(t) | + |R_{i,\varepsilon}^{(2)}(t) | + |R_{\varepsilon}^{(0)}(t) | \leqslant C \varepsilon e^{Ct}.
\end{equation*}    
\end{lemma}

\begin{proof} 
We consider $\psiapp = g(q,p,C,\zeta) \in \mathcal{G}$ the variational approximation constructed in Theorem~\ref{theo:parameters}. Here, we only detail the proof for for the term $\langle V\rangle_{\psi_{\mathrm{app}}}$, while the other estimates are carried out using the same arguments. We start by writing the second order Taylor expansion of $V$ at $q$. For $x \in \R^d$ and $t \geqslant 0$
\[V(x) = V(q(t)) + (x-q(t)) \cdot \nabla V(q(t)) + r_{2}(t)\]
\noindent where $r_{2}$ is the corresponding remainder. Since $V$ is at most quadratic as stated in Assumption~\ref{hyp:potential} and $r_{2}$ involves a polynomial of degree 2, there exists a constant $C > 0$, independent of $\varepsilon$ such that for all~$t \geqslant 0$
\[ |r_{2}(t) | \leqslant C \| x - q(t) \|^2.\]

\noindent Then, for~$t \geqslant 0$, we can write 
\begin{align*}
    \langle V\rangle_{\psi_{\mathrm{app}}(t)} & = \ds\int_{\R^d} V(x) |\psiapp (t,x)|^2 \diff x\\
    & = V(q(t)) \normLp{\psiapp (t)}{2}{\R^d}^2 +  \ds\int_{\R^d} (x-q(t)) \cdot \nabla V(q(t)) |\psiapp (t,x)|^2 \diff x + \ds\int_{\R^d} R_{2,q}(t) |\psiapp (t,x)|^2 \diff x.
\end{align*}

\noindent Using the conservation of the $L^2$--norm stated in Corollary~\ref{lemma:normgaussienne}
and the fact that the first moment of a centered Gaussian vanishes, we deduce
\begin{equation*}
    \langle V \rangle_{\psi_{\mathrm{app}}(t)} = V(q(t)) + R_0^{\varepsilon}(t), \quad \text{where} \quad |R_0^{\varepsilon}(t)| \leqslant C \ds\int_{\R^d} \| x - q(t) \|^2 |\psiapp (t,x)|^2 \diff x.
\end{equation*}

\noindent Using Lemma~\ref{lem:gaussian_moments} with $\alpha = 1$ and the control stated in Proposition~\ref{prop:existenceparaglob}, we conclude using that
\begin{equation*}
    |R_0^{\varepsilon}(t)| \leqslant C \varepsilon e^{Ct}.
\end{equation*}
\noindent This completes the proof.
\end{proof}

\noindent Since Theorem~\ref{th:WPapp} is stated for a potential that decouples the variables and for a diagonal matrix $D$, we first rewrite System~\eqref{eq:eqpara} in this setting. If we assume that \(C_0 = \operatorname{diag}(c_{0,1},\ldots,c_{0,d})\) with \( \Ima (c_{0,i}) > 0\) for all \(i \in [\![ 1,d]\!],\)
the matrix~$C$, solution to \eqref{eq:eqpara} with $C(0) = C_0$, is also diagonal and for all $t\geqslant 0$, we can write
\[ C(t)=\operatorname{diag}(c_1(t),\ldots,c_d(t)),\]
with for all $i\in\{1,\ldots,d\}$, \[\dot c_i = -c_i^2 -\left\langle V_i''\right\rangle_{\psi_{\rm app}} +2\lambda\operatorname{Im}(c_i), \qquad c_i(0)=c_{0,i}.\]

\noindent We are thus able to extend the comparison result obtained in~\cite{lasserlubich2020}, in Subsection~$3.6$.

\begin{prop} \label{prop:comparaisonpara} Let $0 < \varepsilon \ll 1$ and $\lambda \in \R$. We consider $(q_0,p_0) \in \R^{2d}$, $\zeta_0 \in \C$ and a diagonal matrix~$C_0 = \mathrm{diag}(c_{0,1}, ..., c_{0,d}) \in \C^{d\times d}$ where for all \(i \in [\![ 1,d]\!]\), \( \Ima (c_{0,i}) > 0.\) We assume that the potential $V$ satisfies Assumption \ref{hyp:potential} and decouples the variables as in Theorem~\ref{th:WPapp}. We consider $(q,p,C,\zeta)$ solving~\eqref{eq:eqpara} and $(q_{_{\mathrm{WP}}},p_{_{\mathrm{WP}}},S,b,D)$ solving~\eqref{eq:eqparaWP}, with initial conditions at $t=0$ given by
\begin{gather*}
(q(0),p(0))=(q_0,p_0), \qquad \zeta(0)=\zeta_0, \qquad C(0)=C_0,\\
\text{and} \quad (q_{\mathrm{WP}}(0),p_{\mathrm{WP}}(0))=(q_0,p_0), \qquad S(0)=0, \qquad b(0)=\varepsilon^{\frac{d}{4}}\exp\left(\dfrac{i}{\varepsilon}\zeta_0\right),
\qquad D(0)=iC_0.
\end{gather*}

\noindent If $\psi_{0}^{\varepsilon} = G(q_0,p_0,\zeta_0,C_0)$ has a unit norm, for all $t\geqslant0$, there exists $C>0$, independent of $\varepsilon$, such that
\begin{gather*}
    \|q(t)-q_{_{\mathrm{WP}}}(t)\| + \|p(t)-p_{_{\mathrm{WP}}}(t)\| \leqslant C\varepsilon e^{Ct}, \qquad \sum_{i=1}^d |ic_i(t)-d_i(t)| \leqslant C \varepsilon e^{e^{Ct}}\\
\text{and} \quad \left|\zeta(t)-S(t)+i\varepsilon\log\left(\varepsilon^{-d/4}b(t)\right)\right| \leqslant C\varepsilon e^{Ct}\left(1+|\log\varepsilon|\right).
\end{gather*}
\end{prop}

\begin{proof} Let $0 < \varepsilon \ll 1$ and $\lambda \in \R$. We consider the parameters $(q,p,C,\zeta)$ solution to~\eqref{eq:eqpara} and $(q_{_{\mathrm{WP}}},p_{_{\mathrm{WP}}},S,b,D)$ solution~\eqref{eq:eqparaWP}, with initial conditions at $t=0$ given in the statement of Proposition~\ref{prop:comparaisonpara}. We denote by~\(\psiapp\) the Gaussian associated with the parameters $(q,p,C,\zeta)$. So, with the notation introduced in~\eqref{eq:fonctionnelG}, we write \(\psiapp = G(q,p,C,\zeta).\) For the remainder of the proof, we introduce the transformed variables
$\widetilde{C}$ and $\widetilde{\zeta}$, defined by the following relations
\[\widetilde{C} = iD \quad \text{and} \quad \exp \Big( \dfrac{i}{\varepsilon} \widetilde{\zeta} \Big) = \varepsilon^{-d/4} \exp \Big( \dfrac{i}{\varepsilon} S \Big) b.\]
\noindent By the assumptions on the initial conditions, we have \(\widetilde{C}(0)=C(0) = C_0\) and \(\widetilde{\zeta}(0)= \zeta (0) = \zeta_0.\) We split the proof into several steps.

\smallskip

\noindent \textbf{Step 1: Estimate for $q-q_{_{\mathrm{WP}}}$ and
$p-p_{_{\mathrm{WP}}}$.} From \eqref{eq:eqpara} and \eqref{eq:eqparaWP}, we have
\[
\frac{d}{dt} \delta_q = p-p_{_{\mathrm{WP}}} \quad \text{and} \quad \frac{d}{dt} \delta_p = -  \langle V \rangle_{\psi_{\mathrm{app}}} + V(q_{_{\mathrm{WP}}}).
\]

\noindent where $\delta_q = q-q_{_{\mathrm{WP}}}$ and $\delta_p = p-p_{_{\mathrm{WP}}}$. Using the estimate from Step 1, we obtain
\[ \frac{d}{dt} \delta_p = - \nabla V(q) + \nabla V(q_{_{\mathrm{WP}}}) + R_1^{\varepsilon}.\]

\noindent Since $V$ is at most quadratic, a Taylor expansion yields the existence of a constant $C > 0$ such that
\begin{equation*}
    \| \nabla V(q) - \nabla V(q_{_{\mathrm{WP}}}) \| \leqslant C \| q - q_{_{\mathrm{WP}}}\|.
\end{equation*}

\noindent It follows that
\begin{align*}
    \frac{d}{dt} \Big( \|\delta_p\| + \|\delta_q\| \Big) & \leqslant \Big\|  \frac{d}{dt} \delta_p \Big\| + \Big\|  \frac{d}{dt} \delta_q \Big\| \\
    & \leqslant C \Big( \|\delta_q\| + \|R_1^{\varepsilon}\| \Big) + \|\delta_p\| \leqslant \widetilde{C} \Big( \|\delta_q\| + \|\delta_p\| + \varepsilon \, e^{Ct} \Big).
\end{align*}

\noindent Since $\delta_q(0) = \delta_p(0) = 0$, for all $t \geqslant 0$, this implies that 
\begin{equation*}
    \|\delta_p(t)\| + \|\delta_q(t)\| \leqslant C \varepsilon e^{Ct} + C \ds\int_0^t \big( \|\delta_p(s)\| + \|\delta_q(s)\| \big) \diff s.
\end{equation*}

\noindent The desired estimates for $q-q_{_{\mathrm{WP}}}$ and $p-p_{_{\mathrm{WP}}}$ is obtained using Grönwall's lemma.

\smallskip

\noindent \textbf{Step 2: Estimate for $C-\widetilde C$.} By definition, $\widetilde c_i=id_i$. Then, using Systems~\eqref{eq:eqpara} and~\eqref{eq:eqparaWP}, for all $i \in [\![1,d]\!]$, we have
\[\dot c_i = - c_i^2 - \langle V_i''\rangle_{\psi_{\mathrm{app}}} + 2\lambda\operatorname{Im}(c_i) \quad \text{and} \quad \dot{\widetilde c}_i = -\widetilde c_i^2 - V_i''(q_{_{\mathrm{WP}},i})+
2\lambda\operatorname{Im}(\widetilde c_i).\]

\noindent Consequently, for $i \in [\![1,d]\!]$, denoting $\delta_{c_i} = c_i - \widetilde{c}_i$, we obtain
\[\frac{d}{dt} \delta_{c_i} = -\left[ c_i^2-\widetilde c_i^2\right] -
\left[ \langle V_i''\rangle_{\psi_{\mathrm{app}}} -
V_i''(q_{_{\mathrm{WP}},i}) \right] +
2\lambda \operatorname{Im}(c_i - \widetilde c_i).\]
Using Lemma~\ref{lem:asympaverage}, we know that
\[\langle V_i''\rangle_{\psi_{\mathrm{app}}} - V_i''(q_{_{\mathrm{WP}},i}) = V_i''(q_i) - V_i''(q_{_{\mathrm{WP}},i}) + R^{\varepsilon}_{i,2}.\]

\noindent As in Step 1, using the assumption on the potential $V$, a Taylor expansion and the bound on the remainder~$R^{\varepsilon}_{i,2}$ established in Lemma~\ref{lem:asympaverage}, we obtain a constant $C > 0$ such that
\[
\left| \langle V_i''\rangle_{\psi_{\mathrm{app}}} -
V_i''(q_{_{\mathrm{WP}},i}) \right| \leqslant C
\left( \|q_i-q_{_{\mathrm{WP}},i}\| + \varepsilon e^{Ct} \right).
\]
Moreover, since $c_i$ and $\widetilde c_i$ satisfy exponential bounds for all $i \in [\![1,d]\!]$ (see Propositions~\ref{prop:existenceparaglob} and~\ref{prop:controlparaWP}), we have
\[
|c_i^2-\widetilde c_i^2| = |c_i-\widetilde c_i| |c_i+\widetilde c_i| \leqslant C e^{Ct} |\delta_{c_i}|.
\]

\noindent Therefore, using Step 1, we deduce that for all $t \geqslant 0$
\[|\delta_{c_i}(t) | \leqslant C \varepsilon e^{Ct} + C \ds\int_0^t e^{Cs} |\delta_{c_i}(s)| \diff s.\]

\noindent Applying Grönwall's lemma, we obtain that, for all $i \in [\![1,d]\!]$ and $t \geqslant 0$
\[
| \delta_{c_i}(t) | \leqslant C \varepsilon e^{e^{Ct}}.
\]

\noindent Summing over $i$ yields the desired result for $C - \widetilde{C}$.

\smallskip

\noindent \textbf{Step 3: Estimate for $\zeta-\widetilde{\zeta}$.} By definition, we have
\(\widetilde{\zeta} =S-i\varepsilon
\log\left(\varepsilon^{-d/4}b\right).\) By differentiating with respect to time, we obtain
\begin{equation*}
    \dot{\widetilde{\zeta}} = \dot S - i\varepsilon\frac{\dot b}{b}.
\end{equation*}

\noindent Using the equation satisfied by $S$ and $b$ given in \eqref{eq:eqparaWP}, we deduce
\[\dot{\widetilde{\zeta}} = \dfrac{|p_{_{\mathrm{WP}}}|^2}{2} - V(q_{_{\mathrm{WP}}}) -\frac{\varepsilon}{2}\sum_{i=1}^d d_i - \varepsilon\lambda\log|b|^2.\]
\noindent Since
\(\widetilde c_i=id_i\) and \(\log|b|^2 = \frac{d}{2}\log\varepsilon - \frac{2}{\varepsilon} \operatorname{Im}(\widetilde{\zeta})\)
we get
\[\dot{\widetilde{\zeta}} = \dfrac{|p_{_{\mathrm{WP}}}|^2}{2} - V(q_{_{\mathrm{WP}}}) + \frac{i\varepsilon}{2}
\operatorname{tr}(\widetilde C) - \frac{\lambda d \varepsilon}{2}\log\varepsilon + 2 \lambda \operatorname{Im}(\widetilde{\zeta}).\]

\noindent
Consequently, setting $\delta_{\zeta}=\zeta-\widetilde{\zeta}$, we obtain
\begin{equation*}
    \frac{d}{dt}\delta_{\zeta}
    = \frac{|p|^2 - |p_{_{\mathrm{WP}}}|^2}{2}
    - \langle V\rangle_{\psi_{\mathrm{app}}}
    + V(q_{_{\mathrm{WP}}})
    + \frac{i\varepsilon}{2} \operatorname{tr}(C - \widetilde C) + 2\lambda \operatorname{Im}(\zeta - \widetilde{\zeta}) +\frac{\varepsilon}{4}
    \operatorname{tr}\left(
    (\operatorname{Im}(C))^{-1}
    \langle\operatorname{Hess}V\rangle_{\psi_{\mathrm{app}}}
    \right)
    +\frac{\lambda d\varepsilon}{2}\log\varepsilon.
\end{equation*}

\noindent First, using the exponential control of Propositions~\ref{prop:existenceparaglob} and~\ref{prop:controlparaWP}, we have 
\begin{equation}
\label{eq:controledifferencp}
    \left| \frac{|p|^2}{2} -\frac{|p_{_{\mathrm{WP}}}|^2}{2} \right| = \left| (p+p_{_{\mathrm{WP}}})\cdot(p-p_{_{\mathrm{WP}}}) \right| \leqslant \frac{|p|+|p_{_{\mathrm{WP}}}|}{2} \|\delta_p\| \leqslant C e^{Ct} \|\delta_p \|.
\end{equation}

\noindent Next, we split the potential contribution as
\[
\begin{aligned}
    V(q_{_{\mathrm{WP}}})
    -\langle V\rangle_{\psi_{\mathrm{app}}}
    &=
    V(q_{_{\mathrm{WP}}})-V(q)
    +V(q)-\langle V\rangle_{\psi_{\mathrm{app}}}.
\end{aligned}
\]
By the mean value theorem, Assumption~\ref{hyp:potential} and using Lemma~\ref{lem:asympaverage}, we obtain \[|V(q_{_{\mathrm{WP}}})-V(q)| \leqslant C \|\delta_q \| (1 + \|\delta_q \|) \quad \text{and} \quad |V(q)-\langle V\rangle_{\psi_{\mathrm{app}}}| \leqslant C\varepsilon e^{Ct}.\] Consequently,
\begin{equation} 
\label{eq:controledifferenceV}
|V(q_{_{\mathrm{WP}}})-\langle V\rangle_{\psi_{\mathrm{app}}}| \leqslant C \Big( \varepsilon e^{Ct} + \|\delta_q \| (1 + \|\delta_q \|) \Big).
\end{equation}

\noindent
Moreover, since $\operatorname{Hess}V$ is bounded by
Assumption~\ref{hyp:potential}, and using the control of the matrix~$(\operatorname{Im}(C))^{-1}$ established in
Remark~\ref{rq:controlImC}, we obtain
\[\left|\frac{\varepsilon}{4} \operatorname{tr}\left( (\operatorname{Im}(C))^{-1}\langle\operatorname{Hess}V\rangle_{\psi_{\mathrm{app}}}\right)\right| \leqslant C\varepsilon e^{Ct}.\]

\noindent
Finally, using the estimates
\[\left|2\lambda\left(\operatorname{Im}(\zeta)-\operatorname{Im}(\widetilde{\zeta})\right)\right| \leqslant C|\delta_{\xi}| \quad \text{and} \quad \left| \frac{i\varepsilon}{2}\left(\operatorname{tr}(C)-\operatorname{tr}(\widetilde C)\right)\right|\leqslant C\varepsilon\sum_{i=1}^d|\delta_{c_i}|\]
and the bounds obtained in Steps~1 and~2, we infer that
\begin{equation*}
\left|\frac{d}{dt}\delta_{\zeta}\right| \leqslant C |\delta_{\zeta} | + C \varepsilon e^{Ct} + C\varepsilon|\log\varepsilon|.
\end{equation*}

\noindent Since $\delta_{\zeta}(0) = 0$, for all $t \geqslant 0$, this implies that
\[|\delta_{\zeta}(t)| \leqslant C \ds\int_0^t |\delta_{\zeta}(s)| \diff s + C \varepsilon (e^{Ct} + t |\log \varepsilon |).\]
 
\noindent By Grönwall's lemma, we therefore obtain
\[
|\delta_{\zeta}(t)|
\leqslant
C\varepsilon e^{Ct}
\left(e^{Ct} + t|\log\varepsilon|\right).
\]
\noindent After possibly changing the constant $C$, we conclude that
\[|\delta_{\zeta}(t)|\leqslant C\varepsilon e^{Ct}\left(1+|\log\varepsilon|\right).\]
\noindent This concludes the proof.
\end{proof}

\section{Nonlinear superposition} \label{sec:superposition} In this section, the proof of the nonlinear superposition result (Theorem \ref{th:superposition}) is established. We follow the strategy of the proof presented in~\cite{carles2025log}. The initial data is taken as a sum of two Gaussian functions of~$\mathcal{G}$ and we approximate the solution by the superposition of two Gaussian functions, evolving according to the variational dynamics associated with each corresponding initial condition. We proceed as in the proof of Theorem \ref{th:errorbound}, handling additional terms coming from the interaction between the two Gaussian states $\psi^\eps_{1,\mathrm{app}}$ and $\psi^\eps_{2,\mathrm{app}}$. In fact, if we write $e^\eps = \psi^\eps - \psi^\eps_{1,\mathrm{app}} - \psi^\eps_{2,\mathrm{app}}$, since
\begin{align*}
i \varepsilon \partial_t \psi^{\varepsilon} - \mathrm{H}^{\varepsilon} \psi^{\varepsilon} &= \nl{\psi^{\varepsilon}} \quad ; \quad \psi^\eps_0=\psi^\eps_{1,0}+\psi^\eps_{2,0} \\
 i\varepsilon \partial_t \psi^\eps_{j,\mathrm{app}} - \mathrm{P}_{\psi^\eps_{j,\mathrm{app}}} \mathrm{H}^{\varepsilon} \psi^\eps_{j,\mathrm{app}} &= \nl{\psi^\eps_{j,\mathrm{app}}}\quad ; \quad \psi^\eps_{j,\mathrm{app}}(0)=\psi^\eps_{j,0}, \quad j=1,2
\end{align*}

\noindent we have 
\begin{equation*}
i \varepsilon \partial_t e^{\varepsilon} - \mathrm{H}^{\varepsilon} e^{\varepsilon} = \mathrm{P}_{\psi^\eps_{1,\mathrm{app}}}^{\perp} \, \mathrm{H}^{\varepsilon} \psi^\eps_{1,\mathrm{app}} + \mathrm{P}_{\psi^\eps_{2,\mathrm{app}}}^{\perp} \, \mathrm{H}^{\varepsilon} \psi^\eps_{2,\mathrm{app}} + \eps N_S^\eps + \eps N_I^\eps \quad ; \quad e^{\varepsilon} (0) = 0,
\end{equation*}
where the nonlinear term has two components
\begin{align*}
N_S^\eps & = \lambda \Big( \psi^{\varepsilon} \log \, \lvert \psi^{\varepsilon} \rvert^2 -  \, (\psi^\eps_{1,\mathrm{app}} +\psi^\eps_{2,\mathrm{app}}) \log \, \lvert \psi^\eps_{1,\mathrm{app}} +\psi^\eps_{2,\mathrm{app}} \rvert^2 \Big)\\
N_I^\eps & =   \lambda \Big( (\psi^\eps_{1,\mathrm{app}} +\psi^\eps_{2,\mathrm{app}} ) \log \, \lvert \psi^\eps_{1,\mathrm{app}} +\psi^\eps_{2,\mathrm{app}} \rvert^2 -  \, \psi^\eps_{1,\mathrm{app}} \log \, \lvert \psi^\eps_{1,\mathrm{app}} \rvert^2 -  \, \psi^\eps_{2,\mathrm{app}} \log \, \lvert \psi^\eps_{2,\mathrm{app}} \rvert^2\Big).
\end{align*}

\noindent We refer to \(N_I^\varepsilon\) as the \textit{interaction source term} and to \(N_S^\varepsilon\) as the \textit{separated source term}. Proceeding as in the proof of Theorem \ref{th:errorbound}, we obtain
\begin{equation*}
    \dfrac{\diff}{\diff t} \normLp{e^{\varepsilon}}{2}{\R^d} \leqslant \sum\limits_{i=1}^2 \mathrm{dist} \Big( \dfrac{1}{i \varepsilon} \mathrm{H}^{\varepsilon} {\psi^\eps_{i,\mathrm{app}}} , \mathcal{T}_{\psi^\eps_{i,\mathrm{app}}} \mathcal{M} \Big) + \dfrac{1}{\normLp{e^{\varepsilon}}{2}{\R^d}} \Ima \Big( \langle N_S^\varepsilon, e^{\varepsilon} \rangle \Big) + \normLp{N_I^\varepsilon}{2}{\R^d}.
\end{equation*}

\noindent Hence, for all $t \geqslant 0$, we deduce
\begin{equation}
\label{eq:superposition_error}
\normLp{e^{\varepsilon}(t)}{2}{\R^d} \leqslant \int_0^t  \sum\limits_{i=1}^2 \mathrm{dist} \Big( \dfrac{1}{i \varepsilon} \mathrm{H}^{\varepsilon} {\psi^\eps_{i,\mathrm{app}}(s)} , \mathcal{T}_{\psi^\eps_{i,\mathrm{app}}(s)} \mathcal{M} \Big) \, \diff s + \int_0^t \dfrac{\Ima \big( \langle N_S^\varepsilon(s), e^{\varepsilon}(s) \rangle \big)}{\normLp{e^{\varepsilon}(s)}{2}{\R^d}} \,  \, \diff s + \int_0^t \normLp{N_I^\varepsilon(s)}{2}{\R^d} \, \diff s.
\end{equation}

\noindent We now proceed to control separately all the terms that appear on the right-hand side of this inequality.

\subsection{Control of the linear source term} For~$i\in\lbrace 1,2\rbrace$ and $s \in [0,t]$, we have $\psi^\eps_{i,\mathrm{app}}\in\mathcal{G}$, then the term \[ \int_0^t \sum\limits_{i=1}^2 \mathrm{dist} \Big( \dfrac{1}{i\varepsilon}\mathrm{H}^{\varepsilon}\psi^\eps_{i,\mathrm{app}}(s),\mathcal{T}_{\psi^\eps_{i,\mathrm{app}}(s)}\mathcal{M}\Big) \, \diff s \] has already been controlled in the proof of Theorem~\ref{th:errorbound}. Using the control given in~\eqref{eq:majointdistance} and the preceding step, for $i\in\lbrace 1,2\rbrace$, it follows that there exists a constant $C >0$ such that
\begin{equation}
\label{eq:majlinearST}
\int_0^t \sum\limits_{i=1}^2 \mathrm{dist} \Big( \dfrac{1}{i\varepsilon}\mathrm{H}^{\varepsilon}\psi^\eps_{i,\mathrm{app}}(s),\mathcal{T}_{\psi^\eps_{i,\mathrm{app}}(s)}\mathcal{M}\Big) \, \diff s \leqslant C \sqrt{\varepsilon} e^{Ct}.
\end{equation}

\subsection{Control of the separated source term} Let us now consider the second term of the right-hand side in~\eqref{eq:superposition_error}. Replacing $N_S^\varepsilon$ with its expression, we have 
\begin{equation*}
    \Ima \big( \langle N_S^\varepsilon, e^{\varepsilon} \rangle \big) = \int_{\R^d} \Ima \Big( \lambda \Big( \psi^{\varepsilon}(x) \log \, \lvert \psi^{\varepsilon}(x) \rvert^2 - \, [(\psi^\eps_{1,\mathrm{app}}+\psi^\eps_{2,\mathrm{app}})(x)] \log \, \lvert (\psi^\eps_{1,\mathrm{app}}+\psi^\eps_{2,\mathrm{app}})(x) \rvert^2 \Big) \overline{e^{\varepsilon}(x)} \, \diff x \Big).
\end{equation*}

\noindent Applying~\eqref{lem:majlog} with \(z_1 = \lambda \left( \psi^\eps_{1,\mathrm{app}} + \psi^\eps_{2,\mathrm{app}} \right)\) and \(z_2 = \lambda \, \psi^{\varepsilon} = \lambda \left( e^\eps + \psi^\eps_{1,\mathrm{app}} + \psi^\eps_{2,\mathrm{app}} \right)\), for all $s \in [0,t]$, we obtain
\begin{equation*}
    \Ima \Big( \langle N_S^\varepsilon(s) , e^{\varepsilon}(s) \rangle \Big) \leqslant 2 | \lambda | \, \normLp{e^{\varepsilon}(s)}{2}{\R^d}^2.
\end{equation*}

\noindent Consequently, we deduce 
\begin{equation}
    \label{eq:majseparatedST}
    \int_0^t \dfrac{\Ima \Big( \langle N_S^\varepsilon(s), e^{\varepsilon}(s) \rangle \Big)}{\normLp{e^{\varepsilon}(s)}{2}{\R^d}} \, \diff s \leqslant 2 |\lambda| \int_0^t \normLp{e^{\varepsilon}(s)}{2}{\R^d} \, \diff s.
\end{equation}

\noindent Substituting estimates~\eqref{eq:majlinearST} and~\eqref{eq:majseparatedST} into~\eqref{eq:superposition_error}, we obtain that there exists $C > 0$ such that for all $t \geqslant 0$
\begin{equation}
\label{eq:superposition_error2}
\normLp{e^{\varepsilon}(t)}{2}{\R^d} \leqslant C \sqrt{\varepsilon} e^{Ct} + 2 | \lambda |  \int_0^t \,  \normLp{e^{\varepsilon}(s)}{2}{\R^d} + \int_0^t \normLp{N_I^\varepsilon(s)}{2}{\R^d}.
\end{equation}

\noindent Thus, it remains to control the term \(\ds\int_0^t \normLp{N_I^\varepsilon(s)}{2}{\R^d} \, \diff s.\) As in \cite{carles2025log,carlesfermanian2011}, we split the analysis of this interaction source term into two cases, depending on the closeness of the trajectories $q_1$ and $q_2$. Therefore, we need to analyze the set where our trajectories $(q_1,p_1)$ and $(q_2,p_2)$, solutions to~\eqref{eq:eqpara}, are close to each other. We first provide a control of the length of the set on which the trajectories remain close to each other, both for finite times (Proposition~\ref{prop:Ieps_court}) and for large times (Proposition~\ref{prop:Ieps_long}).

\subsection{A key proposition on the set of close trajectories} For $T>0$ and $\gamma\in \big(0,\frac{1}{2}\big)$, we introduce
\begin{equation*}
    I^\varepsilon(T)=\left\lbrace t\in [0,T], \; |q_1(t)-q_2(t)|\leqslant \varepsilon^\gamma \right\rbrace.
\end{equation*}

\noindent The next proposition provides an estimate for the length of the set $I^{\varepsilon}(T)$ for a fixed~$T > 0$ independent of~$\varepsilon$. The proof relies on arguments similar to those used in~\cite[Lemma~6.2]{carlesfermanian2011} for the trajectories arising from the wave packet approximation, solutions to~\eqref{eq:eqparaWP} and denoted by $(q_{_{\mathrm{WP}}},p_{_{\mathrm{WP}}})$ in this paper.

\begin{prop}
\label{prop:Ieps_court}
Let $\gamma \in \big(0,\frac{1}{2}\big)$ and $T > 0$. For $k=1,2$, we consider $(q_k,p_k)$ two solutions to \eqref{prop:existenceparaglob} with $(q_1^0,p_1^0)\neq(q_2^0,p_2^0)$ as initial data. There exists $C > 0$, there exists $\varepsilon_0 >0$ such that for all $\varepsilon \in (0,\varepsilon_0)$ $$\left|I^\varepsilon(T)\right|\leqslant C \varepsilon^\gamma.$$
\end{prop}

\begin{proof}
We fix \(\gamma \in \big(0,\frac{1}{2}\big)\) and \(T>0\). By Proposition~\ref{prop:existenceparaglob}, for \(k=1,2\), let \((q_k,p_k)\) denote the global solutions to~\eqref{eq:eqpara} with initial data
\((q_k(0),p_k(0))=(q_k^0,p_k^0),\) where \((q_1^0,p_1^0)\neq(q_2^0,p_2^0)\). We also consider, for \(k=1,2\), the wave packet parameters \((q_{_{\mathrm{WP}}}^k,p_{_{\mathrm{WP}}}^k)\) introduced in Proposition~\ref{prop:controlparaWP} and associated with the same initial data. For \(c>0\), we define \(I_{_{\mathrm{WP}},c}^\varepsilon(T)\) as the set of times at which the two wave packet trajectories remain close \[I_{_{\mathrm{WP},c}}^\varepsilon(T) = \left\{t\in[0,T], \;\; \left|q_{_{\mathrm{WP}}}^1(t) -  q_{_{\mathrm{WP}}}^2(t)\right| \leqslant c\,\varepsilon^\gamma \right\}.\]
For \(c=1\), this set coincides with the one considered in~\cite[Lemma~6.1]{carlesfermanian2011}. For arbitrary \(c>0\), the proof proceeds in the same way, yielding a constant \(\widetilde{C}>0\) such that
\[ \left|I_{_{\mathrm{WP}},c}^\varepsilon(T)\right|
\leqslant \widetilde{C}\,\varepsilon^\gamma.\]
\noindent Then, by Proposition~\ref{prop:comparaisonpara} and since $\gamma < \frac{1}{2}$, we deduce that there exists $c>0$ such that \(I^\varepsilon(T) \subset I_{_{\mathrm{WP}},c}^\varepsilon(T),\) which concludes the proof.
\end{proof}

\noindent The task is now to control the length of the set $I^{\varepsilon}$ for large times. This is established in Proposition~\ref{prop:Ieps_long} below. We restrict our analysis to the case \(d=1\). 
Before proving this proposition, we establish a useful lemma concerning some energy terms.

\begin{lemma}
\label{lem:energy}
For $d = 1$, we consider $(q,p)$ the solution to~\eqref{eq:eqpara} and $(q_{_{\mathrm{WP}}},p_{_{\mathrm{WP}}})$ the solution to~\eqref{eq:eqparaWP} with the same initial data $(q_0,p_0)$ at time $t = 0$. For $t \geqslant 0$, we define
\begin{equation*}
    E(t) = \dfrac{|p(t)|^2}{2} + V(q(t)) \quad \text{and} \quad E_{_{\mathrm{WP}}}(t) = \dfrac{|p_{_{\mathrm{WP}}}(t)|^2}{2} + V(q_{_{\mathrm{WP}}}(t)).
\end{equation*}

\noindent The following assertions hold.

\smallskip

\begin{enumerate}
    \item For all $t \geqslant 0$, we have \[E_{_{\mathrm{WP}}}(t) = E_{_{\mathrm{WP}}}(0) = \dfrac{|p_0|^2}{2} + V(q_0).\]
    \item For all $t \geqslant 0$, we have \[E(t) = E_{_{\mathrm{WP}}}(0) + R^{\varepsilon}(t),\] where the remainder~$R^{\varepsilon}$ satisfies the following control: there exists a constant~$C > 0$ such that \[R^{\varepsilon}(t) \leqslant C \varepsilon e^{Ct}.\]
\end{enumerate}
\end{lemma}

\begin{proof}
Let $(q_0,p_0)\in\mathbb{R}^2$. We consider $(q,p)$ the solution to~\eqref{eq:eqpara} and $(q_{_{\mathrm{WP}}},p_{_{\mathrm{WP}}})$ the solution to~\eqref{eq:eqparaWP}, with the same initial condition \((q(0),p(0))=(q_{_{\mathrm{WP}}}(0),p_{_{\mathrm{WP}}}(0))=(q_0,p_0).\)

\smallskip

\begin{enumerate}[leftmargin=*, labelindent=0pt]
    \item Since $(q_{_{\mathrm{WP}}},p_{_{\mathrm{WP}}})$ solves the Hamiltonian system~\eqref{eq:eqparaWP}, the corresponding energy $E_{\mathrm{WP}}$ is conserved along the flow. This establishes the first assertion.
\smallskip
\item From the first point, we have that for $t \geqslant 0$, \(E(t) - E_{_{\mathrm{WP}}}(0) = E(t) - E_{_{\mathrm{WP}}}(t).\)
\noindent But, by definition
\begin{equation*}
    E(t) - E_{_{\mathrm{WP}}}(t) = \frac{|p(t)|^2}{2} - \frac{|p_{_{\mathrm{WP}}}(t)|^2}{2} + V(q_{_{\mathrm{WP}}}(t)) - V(q(t)).
\end{equation*}

\noindent Therefore, using the inequalities~\eqref{eq:controledifferencp} and~\eqref{eq:controledifferenceV} established in the proof of Proposition~\ref{prop:comparaisonpara}, we deduce that for all $t \geqslant 0$, we have
\begin{equation*}
    \Big| E(t) - E_{_{\mathrm{WP}}}(0) \Big| \leqslant C e^{Ct} \|p(t) - p_{_{\mathrm{WP}}}(t) \| + C \|q(t) - q_{_{\mathrm{WP}}}(t) \| (1 + \|q(t) - q_{_{\mathrm{WP}}}(t)\|),
\end{equation*}

\noindent with $C$ a positive constant. To conclude, it remains to apply the estimates on $\|p - p_{_{\mathrm{WP}}} \|$ and $\|q - q_{_{\mathrm{WP}}} \|$ established in Proposition~\ref{prop:comparaisonpara}.
\end{enumerate}
\end{proof}

\vspace{-0.8cm}
\begin{remark} The quantity $E_{_{\mathrm{WP}}}$ corresponds to the energy associated with the wave packets parameters considered in~\cite{carlesfermanian2011}. In our setting, we decide to consider the same energy but evaluated in our trajectories: this is what we denote by $E$. In the previous lemma, we first observe that the energy $E$ is close to $E_{_{\mathrm{WP}}}$, up to an error of $\varepsilon$. Moreover, using Remark~\ref{rq:energypara} of Appendix~\ref{appendix:existence} (where we conpute the energy of the variational approximation in terms of its parameters, denoted by $E^{\varepsilon}(q,p,Q,P,\zeta)$) and the asymptotic of the average of $V$ stated in Lemma~\ref{lem:asympaverage}, we also notice that $E$ is the dominant term that appears in the energy~$E^{\varepsilon}(q,p,Q,P,\zeta)$.
\end{remark}

\noindent We now prove Proposition~\ref{prop:Ieps_long}, following the strategy developed in the proof of Proposition~$6.3$ in~\cite{carlesfermanian2011}, where a similar comparison result is established for the wave packet trajectories~$(q_{_{\mathrm{WP}}},p_{_{\mathrm{WP}}})$. For completeness and to make the article self-contained, we reproduce the argument here in the present setting, adapting it to the quantities $(q,p)$ considered in this work.

\begin{prop}
\label{prop:Ieps_long}
Let $\gamma \in \big(0,\frac{1}{2}\big)$. For $d=1$ and $k=1,2$, we consider~$(q_k,p_k)$ solutions to \eqref{eq:eqpara} and $(q^k_{_{\mathrm{WP}}},p^k_{_{\mathrm{WP}}})$ solutions to~\eqref{eq:eqparaWP}, with the same initial data $(q_1^0,p_1^0)\neq(q_2^0,p_2^0)$ at time $t=0$. For $k=1,2$, we consider the following associated quantities \[ E_{_{\mathrm{WP}}}^k = \dfrac{(p_k^0)^2}{2}+V(q_k^0)\]

\noindent and we recall that for all $t \geqslant 0$, \[I^\varepsilon(t)=\left\lbrace s \in [0,t], \; |q_1(s)-q_2(s)|\leqslant \varepsilon^\gamma \right\rbrace.\] If $E := E_{_{\mathrm{WP}}}^1 - E_{_{\mathrm{WP}}}^2 \neq 0$, there exist $C, \widetilde{C}>0$, independent of $\varepsilon$, such that for all $t\in \left[0, C \log \varepsilon^{-1}\right]$
$$|I^\varepsilon(t)| \leqslant \widetilde{C} \dfrac{\varepsilon^\gamma}{E^2}  \textrm{e}^{\widetilde{C}t}.$$
\end{prop}

\begin{proof}
Let us consider $\gamma \in \big(0,\frac{1}{2}\big)$ and $t \geqslant 0$. We introduce~$J^{\varepsilon}(t)$ an interval of maximal length included in~$I^{\varepsilon}(t)$. We denote $N^{\varepsilon}(t)$ the number of such intervals. Therefore, we have
\begin{equation}
    \label{eq:tailleIeps}
    |I^\varepsilon(t)| \leqslant N^{\varepsilon}(t) \times |J^\varepsilon(t)|.
\end{equation}

\noindent \textbf{Step 1: Control of the length of $J^{\varepsilon}(t)$.} Let $(t_1,t_2) \in J^{\varepsilon}(t) \times J^{\varepsilon}(t)$. Our goal is to obtain an upper bound for the quantity $|t_1-t_2|$. Without loss of generality, we may assume that~$t_1 < t_2$. According to the mean value theorem, there exists $t^{\star} \in [t_1,t_2]$ such that 
\[q_1(t_1) - q_2(t_1) = q_1(t_2) - q_2(t_2) + (t_1 - t_2) \Big(p_1(t^{\star}) - p_2(t^{\star}) \Big).\]

\noindent Since $J^\varepsilon(t) \subset I^\varepsilon(t)$, by definition, it follows that
\[ \Big| \left( q_1(t_1) - q_2(t_1) \right) - \left( q_1(t_2) - q_2(t_2) \right) \Big| \leqslant 2 \varepsilon^{\gamma}.\]

\noindent On the other hand, by the reverse triangle inequality, we have
\begin{equation}
\label{eq:conjuguep1moinsp2}
    | p_1(t^{\star}) - p_2(t^{\star}) | \geqslant \Big| |p_1(t^{\star})| - | p_2(t^{\star})| \Big| = \dfrac{\Big| |p_1(t^{\star})|^2 - | p_2(t^{\star})|^2 \Big|}{|p_1(t^{\star})| + | p_2(t^{\star})|}.
\end{equation}

\noindent Thanks to Proposition~\ref{prop:existenceparaglob}, there exists $C > 0$ such that
\[ |p_1(t^{\star})| + | p_2(t^{\star})| \leqslant C e^{C t^{\star}} \leqslant C e^{C t} .\]

\noindent For $k = 1,2$, we set \(E^k = \dfrac{|p_k|^2}{2} + V(q_k).\) Using Lemma~\ref{lem:energy}, the mean value theorem and the reverse triangle inequality, we get
\begin{align*}
    \Big| |p_1(t^{\star})|^2 - | p_2(t^{\star})|^2 \Big| & = 2 \Bigg| E_1(t^{\star}) - E_2(t^{\star}) - \Big( V(q_1(t^{\star})) - V(q_2(t^{\star})) \Big) \Bigg|\\
    & = 2 \Bigg| E_{_{\mathrm{WP}}}^1 - E_{_{\mathrm{WP}}}^2 + R_1^{\varepsilon}(t^{\star}) - R_2^{\varepsilon}(t^{\star}) - \Big( V'(q^{\star}) (q_1(t^{\star}) - q_2(t^{\star})) \Big)\Bigg|\\
    & = 2 \Bigg| E - \Big(R_2^{\varepsilon}(t^{\star}) - R_1^{\varepsilon}(t^{\star}) \Big) - \Big( V'(q^{\star}) \times (q_1(t^{\star}) - q_2(t^{\star})) \Big)\Bigg|\\
    & \geqslant 2 |E| - 2 \Big| R_2^{\varepsilon}(t^{\star}) - R_1^{\varepsilon}(t^{\star}) \Big| - 2 \Big| V'(q^{\star}) \times (q_1(t^{\star}) - q_2(t^{\star})) \Big|
\end{align*}

\noindent with $q^{\star} \in \Big[ \min ( q_1(t^{\star}),q_2(t^{\star}) ), \max ( q_1(t^{\star}),q_2(t^{\star}) ) \Big]$ and for $i \in \lbrace 1, 2 \rbrace$
\[ \Big| R_2^{\varepsilon}(t^{\star}) - R_1^{\varepsilon}(t^{\star}) \Big| \leqslant C \varepsilon e^{Ct}. \]

\noindent Since the potential $V$ is at most quadratic and by the fact that $t^{\star} \in J^{\varepsilon}(t) \subset I^{\varepsilon}(t)$, there exists a constant $C > 0$ such that 
\[ \Big| V'(q^{\star}) \times (q_1(t^{\star}) - q_2(t^{\star})) \Big| \leqslant C \varepsilon^{\gamma} e^{Ct}.\]

\noindent Consequently, since $\varepsilon \leqslant \varepsilon^{\gamma}$, there exist constants $C,\widetilde{C}>0$ such that
\begin{equation*}
\Big| |p_1(t^{\star})|^2 - |p_2(t^{\star})|^2 \Big| \geqslant C\Big( |E| - \varepsilon^{\gamma} e^{Ct} \Big) \geqslant \widetilde{C}|E|,
\end{equation*}
provided that \(\varepsilon^{\gamma}e^{Ct}\leqslant \widetilde{C}|E|.\) Hence, there exist $c,\varepsilon_0>0$ such that this condition is satisfied for all $\varepsilon\in(0,\varepsilon_0)$ and for
all~\(t\in\left[0,c\log (\varepsilon^{-1}) \right].\) Therefore, there exists $C > 0$ such that for all~$\varepsilon\in(0,\varepsilon_0)$ and $t\in\left[0,c\log (\varepsilon^{-1}) \right]$, we have
\begin{equation*}
    | p_1(t^{\star}) - p_2(t^{\star}) | \geqslant \dfrac{C |E|}{e^{Ct}}.
\end{equation*}

\noindent Thus, the exist $C,\widetilde{C}>0$, independent of $\varepsilon$, such that for all \( t \in \left[0,C\log\left(\varepsilon^{-1}\right)\right]\) and $\varepsilon \in (0,\varepsilon_0)$, we conclude with
\begin{equation*}
    |t_1-t_2| = \frac{\Big|\big(q_1(t_1)-q_2(t_1)\big) - \big(q_1(t_2)-q_2(t_2)\big)\Big|}{|p_1(t^\star)-p_2(t^\star)|} \leqslant \frac{\widetilde{C} \varepsilon^\gamma e^{Ct}}{|E|}.
\end{equation*}

\noindent \textbf{Step 2: Control of $N^{\varepsilon}(t)$.} The number $N^{\varepsilon}(t)$ is comparable, up to an additive constant, to
$\widetilde{N}^{\varepsilon}(t)$, the number of maximal intervals on which \(|q_1(s)-q_2(s)|\geqslant \varepsilon^\gamma.\) Indeed, by continuity of $|q_1-q_2|$, the intervals on which
$|q_1-q_2|\leqslant\varepsilon^\gamma$ and those on which
$|q_1-q_2|\geqslant\varepsilon^\gamma$ alternate, with equality $|q_1-q_2|=\varepsilon^\gamma$ at each transition point. We therefore turn to the estimation of the number~$\widetilde{N}^{\varepsilon}(t)$. Let
\(\widetilde{J}^{\varepsilon}(t) = \big[ \widetilde{t}_1,\widetilde{t}_2 \big]\) be one of the maximal intervals counted by $\widetilde{N}^{\varepsilon}(t)$. Since we work on the interval $[0,t]$, of length $t$, we deduce that \[\widetilde{N}^{\varepsilon}(t) \times |\widetilde{J}^{\varepsilon}(t)| \leqslant t.\]

\noindent Then, it remains to obtain a lower bound for
$|\widetilde{J}^{\varepsilon}(t)|$, that is, for the quantity \(\left|\widetilde{t}_1-\widetilde{t}_2\right|.\) Using~\eqref{eq:eqpara} and the expansion for $V'$ given in Lemma~\ref{lem:asympaverage}, there exists $C > 0$ such that
\begin{align*}
    p_1(\widetilde{t}_1) - p_2(\widetilde{t}_1) & = p_1(\widetilde{t}_2) - p_2(\widetilde{t}_2) + \ds\int_{\widetilde{t}_1}^{\widetilde{t}_2} \big( \dot{p_1}(s) - \dot{p_2}(s) \big) \, \diff s \\
    & = p_1(\widetilde{t}_2) - p_2(\widetilde{t}_2) + \ds\int_{\widetilde{t}_1}^{\widetilde{t}_2} \Big( V'(q_1(s)) - V'(q_2(s)) + R^{\varepsilon}(s) \Big) \, \diff s
\end{align*}
where for all $s \in \big[ \widetilde{t}_1,\widetilde{t}_2 \big]$, we have
\[ | R^{\varepsilon}(s) | \leqslant C \varepsilon e^{Ct} \leqslant C e^{Ct}.\]

\noindent Then, using Assumption~\ref{hyp:potential} on $V$ and the exponential control of $q_1$ and $q_2$ established in Proposition~\ref{prop:existenceparaglob}, we deduce that there exists $C >0$ such that 
\begin{equation}
\label{eq:equationdifftempstilde}
    \Big| p_1(\widetilde{t}_1) - p_2(\widetilde{t}_1) - \big( p_1(\widetilde{t}_2) - p_2(\widetilde{t}_2) \big) \Big| \leqslant C e^{Ct} \left|\widetilde{t}_1-\widetilde{t}_2\right|.
\end{equation}

\noindent Now, let us look at the term \( \Big| p_1(\widetilde{t}_1) - p_2(\widetilde{t}_1) - \big( p_1(\widetilde{t}_2) - p_2(\widetilde{t}_2) \big) \Big|.\) By continuity and uniqueness of the solution to~\eqref{eq:eqpara}, the function $q_1-q_2$ has a constant sign on $\widetilde{J}^{\varepsilon}(t)$. Without loss of generality, we may assume
that $q_1-q_2$ is positive on $\widetilde{J}^{\varepsilon}(t)$. Since
\[ (q_1-q_2)(\widetilde{t}_1) = (q_1-q_2)(\widetilde{t}_2) = \varepsilon^{\gamma} \quad \text{and} \quad q_1 - q_2 \geqslant \varepsilon^{\gamma}\]
we necessarily have
\[ (\dot q_1-\dot q_2)(\widetilde{t}_1) = (p_1-p_2)(\widetilde{t}_1) >0 \quad \text{and} \quad
(\dot q_1-\dot q_2)(\widetilde{t}_2) = (p_1-p_2)(\widetilde{t}_2) < 0. \]

\noindent Therefore,
\begin{align*}
    \Big| p_1(\widetilde{t}_1) - p_2(\widetilde{t}_1) - \big( p_1(\widetilde{t}_2) - p_2(\widetilde{t}_2) \big) \Big| & \geqslant \Big| p_1(\widetilde{t}_1) - p_2(\widetilde{t}_1) \Big| - \Big| p_1(\widetilde{t}_2) - p_2(\widetilde{t}_2) \Big|\\
    & = \Big( p_1(\widetilde{t}_1) - p_2(\widetilde{t}_1) \Big) + \Big( p_1(\widetilde{t}_2) - p_2(\widetilde{t}_2) \Big).
\end{align*}

\noindent Using~\eqref{eq:conjuguep1moinsp2} and the same reasoning as in Step~1, we deduce that
\begin{equation*}
    \Big| p_1(\widetilde{t}_1) - p_2(\widetilde{t}_1) - \big( p_1(\widetilde{t}_2) - p_2(\widetilde{t}_2) \big) \Big| \geqslant C e^{-Ct} |E|.
\end{equation*}

\noindent for $\varepsilon\in(0,\varepsilon_0)$ and $t\in\left[0,c\log (\varepsilon^{-1}) \right]$ with $\varepsilon_0, c > 0$. Substituting into~\eqref{eq:equationdifftempstilde}, we obtain 
\begin{equation*}
    \left|\widetilde{t}_1-\widetilde{t}_2\right| \geqslant C e^{-Ct} |E|.
\end{equation*}

\noindent Thus, there exists $\widetilde{C}>0$ such that for all $\varepsilon\in(0,\varepsilon_0)$ and $t\in\left[0,c\log (\varepsilon^{-1}) \right]$, we have
\begin{align*}
    \widetilde{N}^{\varepsilon}(t)
    &\leqslant
    \frac{t}{|\widetilde{J}^{\varepsilon}(t)|}\\
    &\leqslant
    \frac{C t e^{Ct}}{|E|} \leqslant
    \frac{\widetilde{C}e^{\widetilde{C}t}}{|E|}.
\end{align*}

\noindent To conclude, it suffices to substitute the estimate for $J^{\varepsilon}$ obtained in Step~1 and the estimate for $\widetilde{N}^{\varepsilon}$ obtained in Step~2 into inequality~\eqref{eq:tailleIeps}, which completes the proof of the proposition.
\end{proof}

\vspace{-0.3cm}
\subsection{Control of the interaction term}

\noindent For $t\geqslant 0$, we are now in a position to control~\(\ds\int_0^t \normLp{N_I^\varepsilon(s)}{2}{\R^d}\,\diff s.\) For all $t\geqslant 0$, we first write
\[ \int_0^t \normLp{N_I^\varepsilon(s)}{2}{\R^d} \, \diff s = \int_{I^\varepsilon(t)} \normLp{N_I^\varepsilon(s)}{2}{\R^d} \, \diff s + \int_{[0,t] \setminus I^\varepsilon(t)} \normLp{N_I^\varepsilon(s)}{2}{\R^d} \, \diff s\]
\noindent and we will estimate the two integrals on the right-hand side separately. The goal of the following subsection is to study the first term, thanks to Propositions~\ref{prop:Ieps_court} and~\ref{prop:Ieps_long}.

\subsubsection{Region of close trajectories} In order to control the term $\normLp{N_I^\varepsilon}{2}{\R^d}$ on $I^{\varepsilon}(t)$, for $t \geqslant 0$, we first recall the following estimate, stated in~\cite{carles2025log}: for all~$\eta >0$, there exists a constant $C := C(\eta)$ such that
\begin{equation*}
|f|^2 \left(\log|f|^2 \right)^2 \leqslant C(\eta) \Big( |f|^{2-\eta} + |f|^{2+\eta} \Big).
\end{equation*}

\noindent This inequality will be used to estimate the \(L^2\)-norm of the interaction term \(N_I^\varepsilon\) in terms of the \(L^r\)-norms of Gaussian functions, with \(r=2-\eta\) and \(r=2+\eta\). These norms can be computed explicitly using Lemma~\ref{lemma:normLpgaussienne} with $\eta \leqslant 1$, which allows us to avoid using the Gagliardo--Nirenberg inequality, as is done in~\cite{carles2025log}. The Gagliardo--Nirenberg inequality imposes an upper bound on the parameter \(\eta\) that depends on the dimension \(d\). No such restriction is needed in our setting.

\smallskip

\noindent For the remainder of the paper, $a \lesssim b$ means that there exists a constant $C > 0$, independent of $\varepsilon$ and that may depend on $\eta$, such that $a \leqslant C  b $. Moreover, all constants will be denoted by \(C\) or \(\widetilde{C}\), even when their values change from one step to another.

\medskip

\noindent We now control $|N_I^\varepsilon|^2$. We get
\begin{align*}
|N_I^\varepsilon|^2 & \lesssim \sum_{j=1,2}
|\psi^\eps_{j,\mathrm{app}}|^2 \left(\log |\psi^\eps_{j,\mathrm{app}}|^2\right)^2 + |\psi^\eps_{1,\mathrm{app}}+\psi^\eps_{2,\mathrm{app}}|^2 \left(\log |\psi^\eps_{1,\mathrm{app}}+\psi^\eps_{2,\mathrm{app}}|^2\right)^2 \\
& \lesssim \sum_{j=1,2} |\psi^\eps_{j,\mathrm{app}}|^{2-\eta} +|\psi^\eps_{j,\mathrm{app}}|^{2+\eta} + (|\psi^\eps_{1,\mathrm{app}}| + |\psi^\eps_{1,\mathrm{app}}|)^{2-\eta} + (|\psi^\eps_{1,\mathrm{app}}| + |\psi^\eps_{1,\mathrm{app}}|)^{2+\eta}\\
& \lesssim  \sum_{j=1,2}|\psi^\eps_{j,\mathrm{app}}|^{2-\eta} +|\psi^\eps_{j,\mathrm{app}}|^{2+\eta},
\end{align*}
\noindent using the fact that for all $x,y \geqslant 0$, for all $\alpha > 0$, we have \( (x + y)^{\alpha} \lesssim x^{\alpha} + y^{\alpha}.\) Then, for $t\geqslant 0$ and for $0 < \eta < 2$, we get
\begin{equation*}
\normLp{N_I^\varepsilon(t)}{2}{\R^d}^2\lesssim \sum_{j=1,2} \left(N_{2-\eta}(\psi^\eps_{1,\mathrm{app}}(t))\right)^{2-\eta}+\left(N_{2+\eta}(\psi^\eps_{2,\mathrm{app}}(t))\right)^{2+\eta}
\end{equation*}

\noindent where \(N_r\) denotes the \(L^r\)-norm, as introduced in Lemma~\ref{lemma:normLpgaussienne}. Therefore, using Lemma~\ref{lemma:normLpgaussienne} (and more precisely Equation~\eqref{eq:Nr} of Remark~\ref{rq:Nr}) with \(r=2-\eta\) and \(r=2+\eta\) (chosen \(\eta\leqslant 1\) to ensure that \(r\geqslant 1\)), together with the control of the matrix $Q$ given in~\eqref{eq:controlpara}, we obtain
\begin{equation*}
\normLp{N_I^\varepsilon(t)}{2}{\R^d}^2 \lesssim \sum_{j=1,2} \eps^{\frac{d\eta}{4}} \left(\deter{Q(t)Q^\ast(t)}\right)^{\frac{\eta}{4}}+\eps^{-\frac{d\eta}{4}} \left(\deter{Q(t)Q^\ast(t)}\right)^{-\frac{\eta}{4}}
\end{equation*}

\noindent computing explicitly the two powers \(\frac{d}{2r} - \frac{d}{4}\) and \(\frac{1}{2r} - \frac{1}{4}\) that appear in the norm computed in~\eqref{eq:Nr}. Using Corollary~\ref{coro:controldetetlambdamin}, we deduce that there exists $C > 0$ such that 
\[\left(\deter{Q(t)Q^\ast(t)}\right)^{-\frac{\eta}{4}} \lesssim e^{Ct} \quad \text{and} \quad \left(\deter{Q(t)Q^\ast(t)}\right)^{\frac{\eta}{4}}\lesssim e^{Ct}.\]

\noindent Consequently, for all $t \geqslant 0$, we get
\begin{equation*}
\normLp{N_I^\varepsilon(t)}{2}{\R^d}^2 \lesssim \eps^{\frac{d\eta}{4}} e^{Ct} + \eps^{-\frac{d\eta}{4}} e^{Ct} \lesssim \eps^{-\frac{d\eta}{4}} e^{\widetilde{C}t}
\end{equation*}
\noindent with $C, \widetilde{C} > 0$. Therefore, we deduce
\begin{equation*}
\int_{I^\eps(t)} \normLp{N_I^\varepsilon(s)}{2}{\R^d}\; \diff s \lesssim \int_{I^\eps(t)} \eps^{-\frac{d\eta}{8}} e^{Cs} \; \diff s \lesssim |I^\eps(t)| \times \eps^{-\frac{d\eta}{8}} e^{CT}.
\end{equation*}

\noindent For $t \in [0,T]$ with $T>0$, $I^\eps(t) \subset I^\eps(T)$, and thus using Proposition~\ref{prop:Ieps_court}, we obtain
\begin{equation}
    \label{eq:estimationNI1_court}
    \int_{I^\eps(t)} \normLp{N_I^\varepsilon(s)}{2}{\R^d}\; \diff s \lesssim \eps^{\gamma - \frac{d\eta}{8}}.
\end{equation}

\noindent Moreover, by Proposition~\ref{prop:Ieps_long}, if $d=1$ and $E \neq 0$, there exists $c_c>0$ such that for~$t\leqslant c_I\log(\eps^{-1})$, we have
\begin{equation}
\label{eq:estimationNI1_long}
\int_{I^\eps(t)} \normLp{N_I^\varepsilon(s)}{2}{\R^d} \; \diff s \lesssim \eps^{\gamma-\frac{d\eta}{8}} e^{Ct}.
\end{equation}

\subsubsection{Region of separated trajectories} In this subsection, we analyze the source term $N_I^\eps$ involving well-separated Gaussian states, that is, states whose centers follow trajectories that remain sufficiently far apart. In other words, we aim to estimate \(\ds\int_{[0,t] \setminus I^\varepsilon(t)} \normLp{N_I^\varepsilon(s)}{2}{\R^d}\,\diff s.\) As in~\cite[Lemma~5.2]{carles2025log}, we need to state a result from \cite{Ferriere20}. We adapt~\cite[Lemma~3.2]{Ferriere20} to our setting and to a sum of two Gaussian functions. This result can, however, be generalized to a sum of \(N\) Gaussian states.

\begin{lemma}[From Lemma~3.2 in~\cite{Ferriere20}]
\label{lem:guillaume}
There exists $C >0$ such that the following holds. For $k \in \lbrace 1,2 \rbrace$, we consider $q_k\in \R^d$, $\omega_k \in \R$, $\theta_k$ a real-valued measurable function and $\Lambda_k \in \C^{d \times d}$ such that $i\Lambda_k\in \siegel{d}{\C}$. For~$x \in \R^d$, we introduce $$\varphi(x) =  \varphi_1(x)+\varphi_2(x), \quad \textrm{with} \quad \varphi_k(x) = \textrm{e}^{i\theta_k(x)+\omega_k-\transpose{(x-q_k)}\Lambda_k (x-q_k)}. $$
Setting $\delta_\omega = |\omega_1-\omega_2|$, $\ds \lambda_+= \max_{k=1,2} (\sigma(\Rea \; \Lambda_k))$ and $\ds \lambda_- = \min_{k=1,2} (\sigma(\Rea \; \Lambda_k))$ and assuming that $$\rho := |q_1-q_2|^{-1}< \rho_0 := \min \left(\dfrac{\sqrt{\lambda_+}}{\max\left(\sqrt{\delta_\omega +1},\sqrt{\log 2}\right)}, \sqrt{\dfrac{\lambda_-}{d+2}} \right),$$
\noindent we have $$\Big\|\varphi \log |\varphi|^2 - \varphi_1 \log |\varphi_1|^2 -\varphi_2 \log |\varphi_2|^2 \Big\|_{L^2(\R^d)} \leqslant C\; \dfrac{\lambda_+}{\rho^{1+\frac{d}{2}}\sqrt{\lambda_-}}{e}^{-\frac{\lambda_-}{4\rho^2}+\underset{k=1,2}{\max}\, \omega_k}.$$
\end{lemma}

\noindent To avoid any possible confusion, we slightly modify the statement above to correct a minor typo that does not affect the result: after carefully reading the proof of~\cite[Lemma~3.2]{Ferriere20}, we noticed that the symbols $\sigma$ and $\Rea$ had been interchanged. The same notation is used in~\cite[Theorem~5.5]{PhDFerriere}. This typo affects the statement of~\cite[Lemma~3.2]{Ferriere20} but not the result itself.

\smallskip

\noindent We will now apply Lemma~\ref{lem:guillaume} to control $\normLp{N_I^\eps(s)}{2}{\R^d}$ for $s \in [0,t] \setminus I^{\varepsilon}(t)$ with $t \geqslant 0$. First, we identify
\begin{equation*}
\theta_k(x) = \dfrac{1}{\varepsilon}\left(\transpose{p_k}\cdot(x-q_k)+ \Rea(\zeta_k) \right), \quad \Lambda_k = -\dfrac{i}{2\varepsilon}C_k \quad \text{and} \quad \omega_k = -\dfrac{1}{\varepsilon}\Ima(\zeta_k),
\end{equation*}

\noindent To apply Lemma~\ref{lem:guillaume} in our setting, we want the condition on $\rho$ be satisfied on the set~$[0,t] \setminus I^{\varepsilon}(t)$. To determine when it holds, we will first establish a lower bound for $\rho_0$. We first notice that since the matrices $C_1$ and $C_2$ belong to the Siegel space $\siegel{d}{\C}$, using Proposition \ref{prop:caracSiegel} of Appendix~\ref{appendix:siegel}, there exist two matrices $Q_1$ and $Q_2$ such that \[ \Ima\left(C_k\right)=\left(Q_kQ_k^\ast\right)^{-1}, \quad \text{for } k \in \lbrace 1, 2 \rbrace.\]

\noindent Therefore, we have $\lambda_-=\dfrac{1}{2\varepsilon}\min_{k=1,2} \;\sigma \left(\left(Q_kQ_k^\ast\right)^{-1}\right)$, $\lambda_+=\dfrac{1}{2\varepsilon}\max_{k=1,2} \; \sigma \left(\left(Q_kQ_k^\ast\right)^{-1}\right)$ and $\lambda_-\leqslant \lambda_+$. To obtain bounds on \(\lambda_-\) and \(\lambda_+\), we will use the following estimates on the matrices~\(Q_k\), \(Q_k^{-1}\) and \((Q_k^{\ast})^{-1}\), which follow from Proposition~\ref{prop:controlparameterssymp} and its proof, together with the fact that \(1\leqslant \normtriple{Q_k(t)}\cdot\normtriple{Q_k^{-1}(t)}\): there exists $C > 0$ such that
\begin{equation*}
e^{-Ct} \lesssim \normtriple{Q_k(t)} + \normtriple{Q^{-1}_k(t)} \lesssim e^{Ct} \quad \text{and} \quad \normtriple{(Q_k^{\ast})^{-1}(t)} \lesssim e^{Ct}.
\end{equation*}

\noindent We keep the same notation as before. For \(k=1,2\), we introduce the eigenvalues of the matrix~\((Q_kQ_k^{\ast})^{-1}\), denoted and ordered as follows
\(0 < \lambda_k^1 \leqslant \cdots \leqslant \lambda_k^d.\) For \(k=1,2\) and $j \in [\![1,d]\!]$, we have
\begin{equation*}
| \lambda_k^j | =  \lambda_k^j \leqslant \normtriple{(Q_kQ_k^\ast)^{-1}} \leqslant  \normtriple{(Q_k^\ast)^{-1}} \,\normtriple{Q_k^{-1}} \lesssim e^{Ct}.
\end{equation*}

\noindent Hence, we deduce that \[ \lambda_- \leqslant \lambda_+ \lesssim {\eps}^{-1} e^{Ct}.\] 

\noindent The task is now to find a lower bound for those quantities, which is obtained by considering $\lambda_k^1$, the lowest eigenvalue of each $(Q_kQ_k^\ast)^{-1}$ for \(k=1,2\). By definition of the spectral norm, these eigenvalues are given by 
$$\lambda_k^1 =  \dfrac{1}{2\eps}\cdot \dfrac{1}{\normspec{Q_k(t)}^2}, $$

\noindent and using the equivalence of norms, for $k=1,2$, we obtain \[ e^{-Ct} \lesssim \lambda_k^1.\]

\noindent Thus, we have a lower bound for $\lambda_-$, that is $\eps^{-1} e^{-Ct} \lesssim \lambda_-$. Then, it is then easy to check that 
$$\dfrac{1}{\sqrt{\eps}} e^{-Ct} \lesssim \sqrt{\dfrac{\lambda_-}{d+2}}.$$

\noindent On the other hand, Relation~\eqref{eq:relationzeta} allows us to write  
$$\dfrac{1}{\eps} \Ima \, \zeta_k = \dfrac{d}{4}\log \left(\pi \eps\right)+ \dfrac{1}{4} \det \left(Q_k(t)Q_k^\ast(t)\right)- \log(N_{0,k}), \quad \text{for } k=1,2.$$

\noindent Proposition \ref{prop:controlparameterssymp} gives
\begin{equation*}
\delta_\omega \leqslant C+\exp(Ct) \lesssim e^{\widetilde{C}t}
\end{equation*}
and this allows us to write 
$$\dfrac{1}{\sqrt{\eps}} e^{-Ct} \lesssim \dfrac{\sqrt{\lambda_+}}{\max\left(\sqrt{\delta_\omega +1},\sqrt{\log 2}\right)}$$
and so, using the definition of $\rho_0$, we obtain
$$\dfrac{1}{\sqrt{\eps}} e^{-Ct} \lesssim \rho_0.$$

\noindent Using \( \rho = \frac{1}{|q_1 - q_2|} \), the definition of $I^\eps(t)$ and Proposition~\ref{prop:controlparameterssymp}, for all $ s \in [0,t] \setminus I^\eps(t)$, we get $$e^{-Cs} \lesssim \rho(s) \leqslant \eps^{-\gamma}.$$

\noindent The assumption of Lemma~\ref{lem:guillaume} is satisfied on $[0,t]\setminus I^\eps(t)$ as long as $\rho<\rho_0$, which holds when
\begin{equation}
\label{eq:conditionlemGuillaume}
    \eps^{-\gamma}\lesssim \dfrac{1}{\sqrt{\eps}} e^{-Cs}, \qquad \text{for } s \in [0,t]\setminus I^\eps(t).
\end{equation}

\noindent If $t \in [0,T]$ with $T > 0$, this condition is always satisfied since $0<\gamma<\frac{1}{2}$. Otherwise, it is fulfilled for \(0 \leqslant t \leqslant c_f \log(\eps^{-1}),\) with $c_f>0$ independent of $\eps$. Under those conditions, we can apply the lemma. We first estimate the bound obtained in the lemma, thanks to several estimates written above. For $s \in [0,t] \setminus I^{\varepsilon}(t)$, we have
\begin{equation*}
\dfrac{\lambda_+(s)}{\rho(s)^{1+\frac{d}{2}}\sqrt{\lambda_-(s)}}{e}^{-\frac{\lambda_-(s)}{4\rho(s)^2}+\underset{k=1,2}{\max} \, \omega_k(s)} \lesssim \dfrac{\eps^{-1} e^{Cs} \times e^{Cs(1+\frac{d}{2})}}{\eps^{-\frac{1}{2}}e^{-Cs}} \times \eps^{-\frac{d}{4}} e^{-C\eps^{2\gamma-1}e^{-Cs}} \times \varepsilon^{-\frac{d}{4}} e^{Cs} \lesssim \eps^{-\frac{1}{2}-\frac{d}{4}} e^{Cs} e^{-C\eps^{2\gamma-1}e^{-Cs}}.
\end{equation*}

\noindent If $t \in [0,T]$, we get
\begin{equation}
\label{eq:estimationNI2_court}
\int_{[0,t] \setminus I^\eps(t)} \normLp{N_I^\eps(s)}{2}{\R^d} \, \diff s \lesssim  \eps^{-\frac{1}{2}-\frac{d}{4}} e^{-C\eps^{2\gamma-1}}.
\end{equation}

\noindent Otherwise, for $t \geqslant 0$, provided that condition~\eqref{eq:conditionlemGuillaume} is satisfied, we obtain
\begin{align*}
\int_{[0,t] \setminus I^{\varepsilon}(t)} \normLp{N_I^\eps(s)}{2}{\R^d} \, \diff s & \lesssim \eps^{-\frac{1}{2}-\frac{d}{4}}  \int_{[0,t] \setminus I^\eps(t)} e^{Cs} e^{-C\eps^{2\gamma-1}e^{-Cs}} \; \diff s \\
& \lesssim \eps^{-\frac{1}{2}-\frac{d}{4}} e^{-C\eps^{2\gamma-1}e^{-Ct}} \int_{[0,t]\setminus I^\eps(t)} e^{Cs} \; \diff s.
\end{align*}

\noindent Using the fact that $[0,t] \setminus I^\eps(t) \subset [0,t]$, we deduce 
\begin{equation}
    \label{eq:estimationNI2_long}
    \int_{[0,t] \setminus I^{\varepsilon}(t)} \normLp{N_I^\eps(s)}{2}{\R^d} \, \diff s \lesssim \eps^{-\frac{1}{2}-\frac{d}{4}} e^{Ct} e^{-C\eps^{2\gamma-1}e^{-Ct}}.
\end{equation}

\noindent We summarize the estimates obtained for the interaction term in the two previous sections. If $t \in [0,T]$ with $T > 0$, using~\eqref{eq:estimationNI1_court} and \eqref{eq:estimationNI2_court}, we obtain the following final control
\begin{equation}
\label{eq:controlNIfinalcourt}
    \int_0^t \normLp{N_I^\varepsilon(s)}{2}{\R^d} \, \diff s = \int_{I^\varepsilon(t)} \normLp{N_I^\varepsilon(s)}{2}{\R^d} \, \diff s + \int_{[0,t] \setminus I^\varepsilon(t)} \normLp{N_I^\varepsilon(s)}{2}{\R^d} \, \diff s \lesssim \eps^{\gamma-\frac{d\eta}{8}} +  \eps^{-\frac{1}{2}-\frac{d}{4}} e{-C\eps^{2\gamma-1}}.
\end{equation}

\noindent Otherwise, if $d = 1$ and $E \neq 0$, for $t \geqslant 0$ that satisfies \( t \leqslant c \log (\varepsilon^{-1})\) with \(c = \min (c_c, c_f),\) we can apply ~\eqref{eq:estimationNI1_long} and \eqref{eq:estimationNI2_long} and we obtain 
\begin{equation}
\label{eq:controlNIfinallong}
    \int_0^t \normLp{N_I^\varepsilon(s)}{2}{\R^d} \, \diff s = \int_{I^\varepsilon(t)} \normLp{N_I^\varepsilon(s)}{2}{\R^d} \, \diff s + \int_{[0,t] \setminus I^\varepsilon(t)} \normLp{N_I^\varepsilon(s)}{2}{\R^d} \, \diff s \lesssim \eps^{\gamma-\frac{d\eta}{8}} e^{Ct} +  \eps^{-\frac{1}{2}-\frac{d}{4}} e^{Ct}  e^{-C\eps^{2\gamma-1} e^{-Ct}}.
\end{equation}

\subsection{End of the proof} We are now in a position to conclude the proof of Theorem~\ref{th:superposition}. We distinguish the two cases of the theorem.

\begin{enumerate}[leftmargin=*, labelindent=0pt]
    \item We consider $T>0$, $0 < \gamma < \dfrac{1}{2}$ and $0 < \eta \leqslant 1$. For $t\in[0,T]$, substituting the estimate established in~\eqref{eq:controlNIfinalcourt} into~\eqref{eq:superposition_error2}, we get
\begin{equation*}
    \normLp{e^{\varepsilon}(t)}{2}{\R^d} \lesssim  \sqrt{\eps} + \int_0^t \normLp{e^{\varepsilon}(s)}{2}{\R^d} \; \diff s + \eps^{\gamma-\frac{d\eta}{8}} +  \eps^{-\frac{1}{2}-\frac{d}{4}} \exp \left(-C\eps^{2\gamma-1} \right).
\end{equation*}

\noindent Applying Grönwall lemma, for all $t\in[0,T]$, we obtain
\begin{equation*}
\normLp{e^{\varepsilon}(t)}{2}{\R^d} \lesssim  \sqrt{\eps} + \eps^{\gamma-\frac{d\eta}{8}} +  \eps^{-\frac{1}{2}-\frac{d}{4}} \exp \left(-C\eps^{2\gamma-1} \right).
\end{equation*}

\noindent By comparison of the growth rates, since $0<\gamma<1/2$ and $\eta>0$ can be chosen arbitrarily small and after readjusting the constants, we conclude in the finite-time case.

    \item We consider \( 0 \leqslant t \leqslant c \log (\varepsilon^{-1})\), $0 < \gamma < \dfrac{1}{2}$ and $0 < \eta \leqslant 1$. Substituting the estimate~\eqref{eq:controlNIfinallong} into~\eqref{eq:superposition_error2}, we obtain
    \begin{equation*}
\normLp{e^{\varepsilon}(t)}{2}{\R^d} \lesssim \sqrt{\eps} e^{Ct} + \int_0^t \normLp{e^{\varepsilon}(s)}{2}{\R^d} \; \diff s + \eps^{\gamma-\frac{d\eta}{8}} e^{Ct} + \eps^{-\frac{1}{2}-\frac{d}{4}} e^{Ct} e^{-C\eps^{2\gamma-1} e^{-Ct}}.
    \end{equation*}

\noindent Grönwall's lemma allows us to write
\begin{equation*}
\normLp{e^{\varepsilon}(t)}{2}{\R^d} \lesssim \sqrt{\varepsilon} e^{Ct} + \eps^{\gamma-\frac{d\eta}{8}} e^{Ct} + \eps^{-\frac{1}{2}-\frac{d}{4}} e^{Ct} e^{-C\eps^{2\gamma-1} e^{-Ct}}.
\end{equation*}

\noindent By comparison of growth rates, the leading term is $\eps^{\gamma-\frac{d\eta}{8}}e^{Ct}$, which tends to zero for \(t \leqslant \kappa \log \eps^{-1}.\) Then, the second part of Theorem \ref{th:superposition} is true, choosing a constant smaller than $c$ and $\kappa$ for the Ehrenfest time.
\end{enumerate}

\appendix

\section{Existence of the exact solution and conserved quantities for the variational approximation} 
\label{appendix:existence}

\subsection{Existence of solutions for the nonlinear logarithmic equation} In this subsection, we justify the existence and uniqueness of solutions to the Cauchy problem \eqref{eq:logequation} for any fixed $\varepsilon$. To this end, we apply a general result established in~\cite{carles2025log} and in~\cite{CarlesFerriere21} (see also \cite{ACS20} for further details) where the authors study the more general logarithmic Schrödinger equation \[i \varepsilon \partial_t \psi^{\varepsilon} + \frac{\varepsilon^2}{2} \Delta \psi^{\varepsilon} = V \psi^{\varepsilon} + \lambda \varepsilon^{\alpha} \, \psi^{\varepsilon} \log \left| \psi^{\varepsilon} \right|^2, \] with a potential $V$ satisfying Assumption~\ref{hyp:potential}, $(\lambda,\alpha)\in\mathbb{R}^2$ and an initial data belonging to the space \[\Sigma = \left\{ f \in H^1(\mathbb{R}^d) \; ; \; x \mapsto |x|\,f(x) \in L^2(\mathbb{R}^d;\mathbb{C}) \right\}.\] 

\noindent We recall the following important result, stated in \cite[Proposition~$1.3$]{CarlesFerriere21}, which provides the existence of a unique global solution to the Cauchy problem~\eqref{eq:logequation} with initial data~$\psi_0^{\varepsilon}$ in the space $\Sigma$.

\begin{prop}
\label{prop:existencelog}
    Let $0 < \varepsilon \ll 1$ be fixed and $(\lambda,\alpha) \in \mathbb{R}^2$. If the potential~$V$ satisfies Assumption \ref{hyp:potential}, if the initial data~$\psi^{\varepsilon}_{0} \in \Sigma$, \eqref{eq:logequation} admits a unique global solution \( \psi^{\varepsilon} \in \mathcal{C}(\mathbb{R}, \Lp{2}{\R^d}{\C}).\) Moreover, the mass and the energy are independent of time: for all $t \in \R$, we have
    \begin{align*}
    & \normLp{\psi^{\varepsilon}(t)}{2}{\R^d} = \normLp{\psi^{\varepsilon}_{0}}{2}{\R^d}\\
    & E_{\psi^{\varepsilon}}^\eps(t) = \dfrac{1}{2} \normLp{\eps \nabla \psi^{\varepsilon}(t)}{2}{\R^d}^2 + \int_{\R^d} V(x) |\psi^\eps (t,x)|^2 \; \diff x + \lambda \eps^\alpha\int_{\R^d} |\psi^\eps (t,x)|^2 \; \left(\log \left(|\psi^\eps (t,x)|^2\right)-1\right)\diff x =E_{\psi^{\varepsilon}}^\eps(0).
    \end{align*}
\end{prop}

\noindent In this paper, we study the particular case, called the \textit{critical regime}, corresponding to $\alpha = 1$. Moreover, we consider initial conditions given by complex Gaussian functions, which belong to $\Sigma$. This allows us to apply the previous proposition and establish the global existence and uniqueness of the solution for fixed~$\varepsilon$ in our setting.

\subsection{Mass and energy conservations of the variational approximation} \label{sec:conservations} Here, we prove the mass and energy conservation of the variational approximation, that is, we make the proof of Proposition \ref{prop:normenergyconservation}.

\begin{proof}
The proof relies on the characterization of the solution to \eqref{eq:varapp} given in \eqref{eq:varappbis}, that is: $\psiapp \in \mathcal{M}$ is solution to~\eqref{eq:varapp} if and only if, for all $v \in \mathcal{T}_{\psiapp} \mathcal{M}$, we have
\begin{equation}
    \tag{\ref{eq:varappbis}}
    \big\langle v , i \varepsilon \partial_t \psiapp - \mathrm{H}^{\varepsilon} \psiapp - \nl{\psiapp} \big\rangle = 0.
\end{equation}

\smallskip

\noindent In this proof, for all $u$ that belongs to the Hilbert space $\mathcal{H}$, we use the following notation 
\[ \mathrm{NL}^{\varepsilon} u = \nl{u}.\]

\noindent  We consider $t \geqslant 0$. Since the function~$\psiapp(t)$ belongs to $\mathcal{T}_{\psiapp(t)}\mathcal{M}$, using \eqref{eq:varappbis}, we obtain
\begin{equation*}
    \dfrac{\diff}{\diff t} \normLp{\psiapp(t)}{2}{\R^d}^2
     = 2 \Rea \Big( \big\langle \psiapp(t), \partial_t \psiapp(t) \big\rangle \Big)  = 2 \Rea \Big( i \varepsilon \big\langle \psiapp(t), \,
    \mathrm{H}^{\varepsilon} \psiapp(t) + \mathrm{NL}^{\varepsilon} \psiapp(t) \big\rangle \Big).
\end{equation*}
\noindent
As $\mathrm{H}^{\varepsilon}$ is a self-adjoint operator and
\begin{equation*}
\big\langle \psiapp(t), \mathrm{NL}^{\varepsilon} \psiapp(t) \big\rangle
= \lambda \varepsilon \ds\int_{\R^d} |\psiapp(t,x)|^2
\, \log \, \lvert \psiapp(t,x)\rvert^2 \in \R,
\end{equation*}

\noindent we deduce that
\(\dfrac{\diff}{\diff t} \normLp{\psiapp(t)}{2}{\R^d}^2 = 0.\) Thus, the norm is conserved. Regarding the energy, we compute
\begin{align*}
    \dfrac{\diff}{\diff t} E_{\psiapp}^\eps(t)
    & =  \dfrac{\diff}{\diff t} \Big[ \big\langle \psiapp(t), \, \mathrm{H}^{\varepsilon}\psiapp(t) + \mathrm{NL}^{\varepsilon} \psiapp(t) - \psiapp(t) \big\rangle \Big]\\
    & = 2 \Rea \big( \big\langle \psiapp(t), \partial_t \psiapp(t) \big\rangle + \big\langle \partial_t \psiapp(t), \, \mathrm{H}^{\varepsilon}\psiapp(t) \big\rangle \big) + \big\langle \partial_t \psiapp(t), \, \mathrm{NL}^{\varepsilon} \psiapp(t) \big\rangle + \big\langle \psiapp(t), \, \partial_t \big[ \mathrm{NL}^{\varepsilon} \psiapp(t) \big] \big\rangle.
\end{align*}

\noindent Using the norm conservation established above, we have
\( 2 \Rea \Big( \big\langle \psiapp(t), \partial_t \psiapp(t) \big\rangle \Big) = 0.\) Moreover, since the function~$\partial_t \psiapp(t)$ belongs to $\mathcal{T}_{\psiapp(t)}\mathcal{M}$, using \eqref{eq:varappbis}, we obtain 
\[\Big\langle \partial_t \psiapp(t), \, \mathrm{H}^{\varepsilon}\psiapp(t) \Big\rangle = \Big\langle \partial_t \psiapp(t), \, i \varepsilon \partial_t \psiapp(t) - \mathrm{NL}^{\varepsilon} \psiapp(t) \Big\rangle\]
and then
\begin{equation*}
    2 \Rea \Big( \Big\langle \partial_t \psiapp(t), \, \mathrm{H}^{\varepsilon}\psiapp(t) \Big\rangle \Big) + \Big\langle \partial_t \psiapp(t), \, \mathrm{NL}^{\varepsilon} \psiapp(t) \Big\rangle = 2 \Rea \Big( \Big\langle \partial_t \psiapp(t) , i \varepsilon \partial_t \psiapp(t) \Big\rangle \Big) = 0.
\end{equation*}

\noindent Consequently, we obtain
\[\dfrac{\diff}{\diff t} E_{\psiapp}^\eps(t)=\Big\langle \psiapp(t), \,\partial_t \Big[ \mathrm{NL}^{\varepsilon} \psiapp(t) \Big]\Big\rangle.\]

\noindent A direct computation gives
\[\Big\langle \psiapp(t), \,\partial_t \Big[ \mathrm{NL}^{\varepsilon} \psiapp(t) \Big]\Big\rangle = \partial_t \Big[ \normLp{\psiapp(t)}{2}{\R^d}^2  \Big]\]
\noindent and using the conservation of the norm, we conclude that \(\dfrac{\diff}{\diff t} E_{\psiapp}^\eps(t)=0.\) Thus, the energy is also conserved.
\end{proof}

\begin{remark} \label{rq:energypara} Using the definition of the energy, we can obtain it for the variational approximation $\psiapp$, introduced in Theorem~\ref{theo:parameters}, in terms of its Hagedorn's parameters $(q,p,Q,P,\zeta)$. We denote $E^{\varepsilon}(q,p,Q,P,\zeta)$ this energy and the previous proposition implies that \(E^{\varepsilon}(q,p,Q,P,\zeta)\) is constant over time. Thanks to computations, we have
\begin{multline*}
    E^{\varepsilon}(q,p,Q,P,\zeta)=N_0^2 \Bigg[ \dfrac{|p|^2}{2}+\left \langle V \right\rangle_{\psiapp} + \dfrac{\eps}{4}\mathrm{tr}\left(P^\ast P\right) - \lambda \eps \left(\dfrac{1}{2} \log \left[\det (Q^\ast Q)\right] -\log\left[N_0^2 (\pi \varepsilon)^{-\frac{d}{2}} \right] + \dfrac{d}{2} + 1\right)\Bigg]
\end{multline*}
\noindent where \( N_0 := \normLp{\psiapp(t=0)}{2}{\R^d} = \normLp{\psi_0^{\varepsilon}}{2}{\R^d}.\) This motivates the choice of energy considered in our analysis of nonlinear superpositions, as discussed more precisely in Lemma~\ref{lem:energy}.
\end{remark}

\vspace{-0.6cm}
\section{Siegel space and symplectic structure} \label{appendix:siegel} This appendix is based on notes taken during a lecture by Caroline Lasser. We recall that the Siegel space is defined as follows
\[\siegel{d}{\C} = \left\{ C \in \mathbb{C}^{d \times d} \;\middle|\; C = C^{\mathsf{T}}, \;\; \operatorname{Im}(C) \text{ is positive definite} \right\}. \]

\begin{prop} \label{prop:Siegelopen} The Siegel space is an open subset of the space of symmetric matrices~$\symetrique{d}{\C}$.
\end{prop}

\begin{proof} We consider the continuous map \[ \begin{aligned} \Phi:\symetrique{d}{\C}&\longrightarrow \symetrique{d}{\R},\\ C&\longmapsto \Ima (C). \end{aligned} \] We introduce the set of positive definite symmetric matrices \[ \mathcal{S}^+_{d}(\R) = \left\{ M \in \symetrique{d}{\R} \;\middle|\; M \text{ is positive definite} \right\}. \] 

\noindent The set $\mathcal{S}^+_{d}(\R)$ is open, as it is the preimage of $]0,+\infty[$ under the continuous map that associates to a matrix $M \in \symetrique{d}{\R}$ its smallest eigenvalue. Finally, noticing that \(\siegel{d}{\mathbb C} = \Phi^{-1}\left(\mathcal{S}^+_{d}(\R)\right), \) we conclude that $\siegel{d}{\mathbb C}$ is open in $\symetrique{d}{\mathbb C}$.
\end{proof}

\begin{defi}[Symplectic matrix] \label{def:matsymplectic} Let $Y \in \R^{2d \times 2d}$. $Y$ is called symplectic if it satisfies the following relation
\begin{equation}
\label{eq:defmatsymplectique}
\transpose{Y}JY = J \quad \text{with} \quad J = \begin{pmatrix} 0 & - \mathrm{I}_d  \\ \mathrm{I}_d & 0 \end{pmatrix}.
\end{equation}
\end{defi}

\noindent Here is a characterization of symplectic matrices.

\begin{prop} \label{prop:caracsympmat} We consider $Q$ and $P$ be two matrices that belong to $\C^{d \times d}$. The matrix~$Y = \begin{pmatrix} \Rea (Q) & \Ima (Q) \\ \Rea (P) & \Ima (P) \end{pmatrix}$ is symplectic if and only if the matrices $(Q,P)$ satisfy
\begin{equation}
\label{eq:caracsymp}
    \left\lbrace 
    \begin{array}{cll}
       \transpose{Q}P - \transpose{P}Q & = & 0 \\
       Q^{\ast} P - P^{\ast} Q & = & 2i \mathrm{I}_d.
    \end{array}
    \right. 
\end{equation}
\end{prop}

\begin{proof}
We consider $Q$ and $P$ two matrices of $\C^{d \times d}$. To prove the equivalence, we have to use the following computations. On the one hand, according to Definition~\ref{def:matsymplectic}, we assume that the matrix $Y$ satisfies~\eqref{eq:defmatsymplectique}. Writing \(Q=\Rea(Q)+i\,\Ima(Q)\) and \(P=\Rea(P)+i\,\Ima(P)\), this implies \[ \begin{pmatrix} \transpose{\Rea(P)}\Rea(Q)-\transpose{\Rea(Q)}\Rea(P) & & &\transpose{\Rea(P)}\Ima(Q)-\transpose{\Rea(Q)}\Ima(P) \\[0.2cm] \transpose{\Ima(P)}\Rea(Q)-\transpose{\Ima(Q)}\Rea(P) & & &\transpose{\Ima(P)}\Ima(Q)-\transpose{\Ima(Q)}\Ima(P) \end{pmatrix} = \begin{pmatrix} 0 & -\mathrm{I}_d \\ \mathrm{I}_d & 0 \end{pmatrix}. \] On the other hand, we obtain \[ \left\{ \begin{array}{rcl} \transpose{Q}P-\transpose{P}Q &=& \Big[ \transpose{\Rea(Q)}\Rea(P) -\transpose{\Rea(P)}\Rea(Q) +\transpose{\Ima(P)}\Ima(Q) -\transpose{\Ima(Q)}\Ima(P) \Big] \\ && +i\, \Big[ \transpose{\Ima(Q)}\Rea(P) -\transpose{\Ima(P)}\Rea(Q) +\transpose{\Rea(Q)}\Ima(P) -\transpose{\Rea(P)}\Ima(Q) \Big], \\[0.2cm] Q^{\ast}P-P^{\ast}Q &=& \Big[ \transpose{\Rea(Q)}\Rea(P) -\transpose{\Rea(P)}\Rea(Q) +\transpose{\Ima(P)}\Ima(Q) -\transpose{\Ima(Q)}\Ima(P) \Big] \\ && +i\, \Big[ \transpose{\Ima(Q)}\Rea(P) -\transpose{\Ima(P)}\Rea(Q) +\transpose{\Rea(P)}\Ima(Q) -\transpose{\Rea(Q)}\Ima(P) \Big]. \end{array} \right. \] Combining these identities with the relations obtained from $\transpose{Y}JY=J$, we obtain the desired result.
\end{proof}

\noindent We can thus derive a characterization of the Siegel space.

\begin{prop}[Characterization of the Siegel space] \label{prop:caracSiegel} Let $Q, P$ and $Y$ be matrices as given in Proposition~\ref{prop:caracsympmat}. If the matrix $Y$ is symplectic, then $Q$ and $P$ are invertible and the matrix $C=PQ^{-1}$ belongs to the Siegel space $\siegel{d}{\C}$ with
\begin{equation}
\label{eq:imaC}
\Ima(C)=(QQ^{\ast})^{-1}.
\end{equation}

\noindent Conversely, for every matrix $C\in\siegel{d}{\C}$, there exist two matrices $Q,P\in \C^{d \times d}$ such that the matrix \(Y\) is symplectic, the identity \eqref{eq:imaC} holds and \( C=PQ^{-1}.\)
\end{prop}

\begin{proof}
First, we assume that the matrix~$Y$ is symplectic. According to Proposition~\ref{prop:caracsympmat}, the matrices $Q$ and $P$ satisfy~\eqref{eq:caracsymp}. Using the second equation, for any vector $v \in \C^d$, we have
\begin{equation*}
\langle Pv, Qv \rangle - \langle Qv, Pv \rangle = 2i \| v \|^2.
\end{equation*}
\noindent Consequently, neither $Q$ nor $P$ has a nontrivial kernel. Indeed, if $Qv=0$, then \( 2i \|v\|^2=0 \) and $v=0$. The same argument holds for $P$ and it follows that the matrices $Q$ and $P$ are invertible. Now, to prove that $C = PQ^{-1}$ belongs to $\siegel{d}{\C}$, we need to show that~$\transpose{C} = C$ and that $\Ima (C)$ is positive definite. For the first point, if we multiply the first equation of~\eqref{eq:caracsymp} on the left by the matrix~$\transpose{(Q^{-1})}$ and on the right by the matrix~$Q^{-1}$ and we obtain \[PQ^{-1} - \transpose{(Q^{-1})}\transpose{P} = 0,\] which implies that $C = \transpose{C}$. For the second point, using the expression of $C$ and the second equation of~\eqref{eq:caracsymp}, we compute \[ \Ima(C):=\frac{C-C^{\ast}}{2i} = (Q^{\ast})^{-1} \times \dfrac{Q^{\ast}P - P^{\ast}Q}{2i} \times Q^{-1} = (QQ^{\ast})^{-1}.\] Consequently, we have obtained~\eqref{eq:imaC}. Since the matrix $Q^{\ast}$ is invertible (because $Q$ is invertible), $QQ^{\ast}$ is positive definite. Consequently, $\Ima(C)$ is also positive definite and we have proved the first implication of the proposition.

\noindent To prove the second implication, let us consider $C \in \siegel{d}{\C}$. Since $\Ima(C)$ is positive definite, we can set $Q = (\Ima(C))^{-\frac{1}{2}}$. Defining $P = CQ$, we obtain $C = PQ^{-1}$. Using the relation $C = \transpose{C}$, we obtain the first equation of~\eqref{eq:caracsymp}. Moreover, combining the identities $\Ima(C) = (Q^{\ast})^{-1}Q^{-1}$, which holds by construction, and $\Ima(C) = \dfrac{C - C^{\ast}}{2i}$, which follows from the definition, we obtain the second equation in~\eqref{eq:caracsymp}. This shows that the matrix $Y$ is symplectic. By the first part of the proof, we then deduce that $C$ satisfies~\eqref{eq:imaC}, which concludes the proof.
\end{proof}

\begin{remark} \label{rq:ImaCnotunique} We note that the decomposition of $\Ima(C)$ obtained in the previous proposition is not unique, since $Q$ and $P$ can be multiplied on the right by a unitary matrix $U$. Indeed, setting $\widetilde{Q} = QU$ and $\widetilde{P} = PU$, we obtain that \(C = \widetilde{P} \widetilde{Q}^{-1}\) and \(\Ima (C) = (\widetilde{Q} \widetilde{Q}^{\ast})^{-1}.\)
\end{remark}

\vspace*{-0.4cm}
\section{Computation of projections onto the tangent space}
\label{app:proj} Since, in the superposition case, we consider Gaussian states whose $L^2$--norms are not necessarily equal to~$1$, we need to extend the computations of the projections onto the tangent space carried out in~\cite{lasserlubich2020}, where Gaussian states are assumed to have unit norm, to the more general setting in which their norms belong to~$(0,1]$. Let $g \in \mathcal{G}$. There exist a vector~\((q, p) \in \mathbb{R}^{2d} \), \( \zeta \in \mathbb{C} \) and~\( C \in \siegel{d}{\C} \) such that for~$x \in \R^d$
\[ g(x) = \exp \left( \frac{i}{\varepsilon} \left( \frac{1}{2} \transpose{(x - q)} C (x - q) \right) + \transpose{p}(x - q) + \zeta \right).\]

\noindent As in Notation~\ref{nota:crochet}, we set $\theta := \normLp{g}{2}{\R^d}$ and we assume that $\theta$ satisfies
\(0 < \theta \leqslant 1\)  (see Lemma~\ref{lemma:normLpgaussienne} with $r =2$ for an explicit formula).

\noindent For $x \in \mathbb{R}^d$, we set
\( y = \frac{1}{\sqrt{\varepsilon}} M^{-1}(x-q)\), where $M$ is a Cholesky factor in the decomposition of the positive definite matrix~$(\operatorname{Im}(C))^{-1}$ \[\operatorname{Im}(C) = (M \transpose{M})^{-1}.\] 

\noindent According to Remark~\ref{rq:ImaCnotunique}, $M$ and $Q$ can be related through multiplication by a unitary matrix. In particular \(\det Q = \pm \det M.\) For \( k \in \mathbb{N}^d \) with \( |k| := \sum_i k_i \leqslant 2 \), for $x \in \R^d$, we define the following multivariate scaled Hermite functions \( \varphi_k : \mathbb{R}^d \to \mathbb{C}\) by
\[
\varphi_k(x) =
\begin{cases}
\dfrac{g(x)}{\theta} & \text{if } |k| = 0, \\ & \\
\dfrac{\sqrt{2}}{\theta}\, y_m\, g(x) & \text{if } |k| = 1 \text{ and } k_m = 1, \\ & \\
\dfrac{1}{\theta \sqrt{2}} (-1 + 2y_m^2)\, g(x) & \text{if } |k| = 2 \text{ and } k_m = 2, \\ & \\
\dfrac{2}{\theta}\, y_m y_n\, g(x) & \text{if } |k| = 2 \text{ and } k_m = k_n = 1.
\end{cases}
\]

\noindent Using~\cite[Lemma 3.12]{lasserlubich2020}, we conclude that \(\{\varphi_k\}_{|k| \leqslant 2}\) is an orthonormal basis of \( \mathcal{T}_g \mathcal{M}\). Thus, in this setting, we can rewrite~\cite[Lemma~3.13]{lasserlubich2020} as follows.

\begin{lemma}
\label{lemma:3.13}
Considering a Gaussian $g \in \mathcal{G}$ and a smooth function $W : \R^d \rightarrow \R$ with at most polynomial growth, satisfying \(\langle \nabla W \rangle_g = 0 \) and \(\langle \mathrm{Hess}\, W \rangle_g = 0.\) Then, we have
\begin{equation*}
\mathrm{P}_g W g = \langle W \rangle_g g.
\end{equation*}
\end{lemma}

\begin{proof}
The proof follows the same arguments as in~\cite{lasserlubich2020}, while taking into account the normalization by $\theta$ in the functions \(\{\varphi_k\}_{|k| \leqslant 2}\) and using the modified definition of the bracket given in Notation~\ref{nota:crochet}.
\end{proof}

\clearpage
\pagestyle{empty}

\bibliography{bibliographie}

@book{livreRobertComb2012,
  title={Coherent States and Applications in Mathematical Physics},
  author={Combescure, Monique and Robert, Didier},
  isbn={9789400701960},
  lccn={2012931507},
  series={Theoretical and Mathematical Physics},
  url={https://books.google.fr/books?id=mJgh0XAemNIC},
  year={2012},
  publisher={Springer Netherlands}
}

@book{carleslivre2021semi,
  title={Semi-classical analysis for nonlinear Schr{\"o}dinger equations: WKB analysis, focal points, coherent states},
  author={Carles, R{\'e}mi},
  year={2021},
  publisher={World Scientific}
}

@book{frenkel1934,
  author    = {Frenkel, Janko Il'ich},
  title     = {Wave Mechanics: Advanced General Theory},
  publisher = {Clarendon Press},
  year      = {1934},
  address   = {Oxford},
  series    = {International Series of Monographs on Physics},
  pages     = {524},
  language  = {English}
}

@article{broeckhove1988equivalence,
  title={On the equivalence of time-dependent variational principles},
  author={Broeckhove, J and Lathouwers, L and Kesteloot, E and Van Leuven, P},
  journal={Chemical physics letters},
  volume={149},
  number={5-6},
  pages={547--550},
  year={1988},
  publisher={Elsevier}
}

@article{dirac1930theory,
  title   = {Note on exchange phenomena in the {T}homas atom},
  author  = {Dirac, Paul Adrien Maurice},
  journal = {Mathematical Proceedings of the Cambridge Philosophical Society},
  volume  = {26},
  number  = {3},
  pages   = {376--385},
  year    = {1930},
  publisher = {Cambridge University Press},
}

@ARTICLE{Hagedorn1980,
       author = {Hagedorn, George Allan},
        title = {Semiclassical quantum mechanics: I. {T}he {\ensuremath{\hbar}}{\textrightarrow}0 limit for coherent states},
      journal = {Commun. Math. Phys.},
         year = {1980},
       volume = {71},
       number = {1},
        pages = {77-93},
}

@article{hari2013coherent,
  title={Coherent states for systems of ${L}^2-$supercritical nonlinear {S}chr{\"o}dinger equations},
  author={Hari, Lysianne},
  journal={Commun. Partial Differ. Equ.},
  volume={38},
  number={3},
  pages={529--573},
  year={2013},
  publisher={Taylor \& Francis}
}

@article{carlesfermanian2011,
  title={Nonlinear coherent states and {E}hrenfest time for {S}chr{\"o}dinger equation},
  author={Carles, R{\'e}mi and Fermanian Kammerer, Clotilde},
  journal={Commun. Math. Phys.},
  volume={301},
  number={2},
  pages={443--472},
  year={2011},
  publisher={Springer}
}

@article{mclachlan1964,
  title={A variational solution of the time-dependent {S}chrodinger equation},
  author={McLachlan, Andrew D},
  journal={Molecular Physics},
  volume={8},
  number={1},
  pages={39--44},
  year={1964},
  publisher={Taylor \& Francis}
}

@article{Carles12interaction,
  title         = {Interaction of coherent states for {H}artree equations},
  author        = {Carles, R{\'e}mi},
  journal       = {Arch. Ration. Mech. Anal.},
  fjournal      = {Archive for Rational Mechanics and Analysis},
  volume        = {204},
  number        = {2},
  pages         = {559--598},
  year          = {2012},
}

@article{CaoCarles,
  title         = {Semiclassical wave packet dynamics for {H}artree equations},
  author        = {Cao, Pei and Carles, R{\'e}mi},
journal       = {Rev. Math. Phys},
  fjournal       = {Reviews in Mathematical Physics},
  volume        = {24},
  number        = {2},
  pages         = {1250004, 40~pp.},
  year          = {2012},
}

@article{APPP11,
  title         = {Semiclassical propagation of coherent states for the {H}artree equation},
  author        = {Athanassoulis, Agissilaos and Paul, Thierry and Pezzotti, Federica and Pulvirenti, Mario},
  fjournal       = {Ann. Henri Poincar{\'e}},
    journal       = {Annales Henri Poincar{\'e}},
  volume        = {12},
  number        = {8},
  pages         = {1613--1634},
  year          = {2011},
}

@article{Adiab,
  title         = {A nonlinear adiabatic theorem for coherent states},
  author        = {Carles, R{\'e}mi and Fermanian Kammerer, Clotilde},
  journal       = {Nonlinearity},
  volume        = {24},
  number        = {8},
  pages         = {2143--2164},
  year          = {2011},
}

@article {CarlesEMS,
    AUTHOR = {Carles, R\'emi},
     TITLE = {Logarithmic {S}chr\"odinger equation and isothermal fluids},
   JOURNAL = {EMS Surv. Math. Sci.},
  FJOURNAL = {EMS Surveys in Mathematical Sciences},
    VOLUME = {9},
      YEAR = {2022},
    NUMBER = {1},
     PAGES = {99--134},
}

@article{HayashiOzawa25,
    author = {Hayashi, M. and Ozawa, T.},
    title = {The {C}auchy problem for the logarithmic {S}chrödinger equation revisited},
    journal = {Ann. Henri Poincaré},
    volume = {26},
    pages = {1209--1238},
    year = {2025},
}

@article{CHO24,
    author = {Carles, R. and Hayashi, M. and Ozawa, T.},
    title = {Low regularity solutions to the logarithmic {S}chrödinger equation},
    journal = {Pure Appl. Anal.},
    volume = {6},
    number = {3},
    pages = {859--871},
    year = {2024},
}

@article{carles2025log,
  title={Propagation of coherent states in the logarithmic {S}chr\"odinger equation},
  author={Carles, R{\'e}mi and Dong, Fangyuan},
  journal={Commun. Math. Sci.},
  volume={23},
  number={8},
  pages={2151--2175},
  year={2025}
}

@article{CarlesFerriere21,
  title={Logarithmic {S}chr{\"o}dinger equation with quadratic potential},
  author={Carles, R{\'e}mi and Ferriere, Guillaume},
  journal = {Nonlinearity},
  volume  = {34},
  number  = {12},
  year    = {2021},
  pages   = {8283--8310}
}

@inproceedings{cazenave80,
  title={{\'E}quations d'{\'e}volution avec non lin{\'e}arit{\'e} logarithmique},
  author={Cazenave, Thierry and Haraux, Alain},
  booktitle={Annales de la Facult{\'e} des sciences de Toulouse: Math{\'e}matiques},
  volume={2},
  pages={21--51},
  year={1980}
}

@article {cazenave83,
    AUTHOR = {Cazenave, Thierry},
     TITLE = {Stable solutions of the logarithmic {S}chr\"odinger equation},
   JOURNAL = {Nonlinear Anal.},
  FJOURNAL = {Nonlinear Analysis. Theory, Methods \& Applications. An
              International Multidisciplinary Journal},
    VOLUME = {7},
      YEAR = {1983},
    NUMBER = {10},
     PAGES = {1127--1140},
}

@article{ACS20,
 author = {Ardila, Alex H. and Cely, Liliana and Squassina, Marco},
 title = {Logarithmic {B}ose-{E}instein condensates with harmonic potential},
 fjournal = {Asymptotic Analysis},
 journal = {Asymptotic Anal.},
 issn = {0921-7134},
 volume = {116},
 number = {1},
 pages = {27--40},
 year = {2020},
}

@article{Ferriere22,
  title         = {{WKB} analysis of the logarithmic nonlinear {S}chr\"odinger equation in an analytic framework},
  author        = {Ferriere, Guillaume},
  journal       = {Commun. Contemporary Mathematics},
  volume        = {24},
  number        = {8},
  pages         = {Article ID 2150082, 34~p.},
  year          = {2022},
}

@article{Ferriere21multi,
  title         = {Existence of multi-solitons for the focusing logarithmic non-linear {S}chr\"odinger equation},
  author        = {Ferriere, Guillaume},
  journal       = {Annales de l'Institut Henri Poincar\'e C, Analyse non lin\'eaire},
    fjournal    = {Annales de l'Institut Henri Poincar\'e C, Analyse non lin\'eaire},
  volume        = {38},
  number        = {3},
  pages         = {841--875},
  year          = {2021},
}

@article{Ferriere21,
  title         = {Convergence rate in {W}asserstein distance and semiclassical limit for the defocusing logarithmic {S}chr\"odinger equation},
  author        = {Ferriere, Guillaume},
  journal       = {Analysis \& PDE},
  volume        = {14},
  number        = {2},
  pages         = {617--666},
  year          = {2021},
}

@phdthesis{PhDFerriere,
  title        = {Equation de Schr\"{o}dinger non-lin\'{e}aire avec non-linéarité logarithmique},
  author       = {Ferriere, Guillaume},
  year         = 2021,
  note         = {\url{https://theses.hal.science/tel-03396661/}},
  school       = {Universit\'{e} de Montpellier},
  type         = {Ph.{D} thesis}
}

@article {BBM76,
    AUTHOR = {Bia{\l}ynicki-Birula, Iwo and Mycielski, Jerzy},
     TITLE = {Nonlinear wave mechanics},
   JOURNAL = {Ann. Physics},
  FJOURNAL = {Annals of Physics},
    VOLUME = {100},
      YEAR = {1976},
    NUMBER = {1-2},
     PAGES = {62--93},
}

@article {AMS14log,
    AUTHOR = {d'Avenia, Pietro and Montefusco, Eugenio and Squassina, Marco},
     TITLE = {On the logarithmic {S}chr\"odinger equation},
   JOURNAL = {Commun. Contemp. Math.},
  FJOURNAL = {Communications in Contemporary Mathematics},
    VOLUME = {16},
      YEAR = {2014},
    NUMBER = {2},
     PAGES = {1350032, 15},
}

@article {CarlesGallagher18,
    AUTHOR = {Carles, R\'{e}mi and Gallagher, Isabelle},
     TITLE = {Universal dynamics for the defocusing logarithmic {S}chr\"odinger equation},
   JOURNAL = {Duke Math. J.},
  FJOURNAL = {Duke Mathematical Journal},
    VOLUME = {167},
      YEAR = {2018},
    NUMBER = {9},
     PAGES = {1761--1801},
}

@article {Ardila16,
    AUTHOR = {Ardila, Alex H.},
     TITLE = {Orbital stability of {G}ausson solutions to logarithmic {S}chr\"odinger equations},
   JOURNAL = {Electron. J. Differential Equations},
  FJOURNAL = {Electronic Journal of Differential Equations},
      YEAR = {2016},
     PAGES = {Paper No. 335, 9},
}

@article {CarlesNouri17,
    AUTHOR = {Carles, R\'{e}mi and Nouri, Anne},
     TITLE = {Monokinetic solutions to a singular {V}lasov equation from a semiclassical perspective},
   JOURNAL = {Asymptot. Anal.},
  FJOURNAL = {Asymptotic Analysis},
    VOLUME = {102},
      YEAR = {2017},
    NUMBER = {1-2},
     PAGES = {99--117},
}

@article {Ferriere20,
    AUTHOR = {Ferriere, Guillaume},
     TITLE = {The focusing logarithmic {S}chr\"odinger equation: analysis of breathers and nonlinear superposition},
   JOURNAL = {Discrete Contin. Dyn. Syst.},
    VOLUME = {40},
      YEAR = {2020},
    NUMBER = {11},
     PAGES = {6247--6274},
}

@misc{Bensaid26,
  title={The {C}auchy problem for logarithmic {S}chr{\"o}dinger equation with a perturbation of power law nonlinearity},
  author={Bensaid, Mohamed},
  year={2026},
  note={Preprint, arXiv:2609.19217v1},
  url={https://arxiv.org/abs/2609.19217v1}
}

@misc{Bensaid26-2,
  author={Bensaid, Mohamed},
  title = {Existence of multi-solitons for logarithmic {S}chrödinger equation with a perturbation of a smooth $L^2–$subcritical function nonlinearity},
  year          = {2026},
  note={Preprint, arXiv:2609.18625},
  url           = {https://arxiv.org/abs/2609.18625}
}

@article{Hayashi26,
  author  = {Hayashi, M.},
  title   = {Uniqueness of solutions for the logarithmic {S}chr{\"o}dinger equation},
  journal = {Mathematische Annalen},
  volume  = {395},
  number   = {29},
  year    = {2026},
}

@article {BCST19,
    AUTHOR = {Bao, Weizhu and Carles, R\'emi and Su, Chunmei and Tang,
              Qinglin},
     TITLE = {Regularized numerical methods for the logarithmic
              {S}chr\"odinger equation},
   JOURNAL = {Numer. Math.},
  FJOURNAL = {Numerische Mathematik},
    VOLUME = {143},
      YEAR = {2019},
    NUMBER = {2},
     PAGES = {461--487},
}

@article {CC22,
    AUTHOR = {Carles, R\'emi and Su, Chunmei},
     TITLE = {Non-uniqueness and nonlinear instability of {G}aussons under
              repulsive harmonic potential},
   JOURNAL = {Comm. Partial Differential Equations},
  FJOURNAL = {Communications in Partial Differential Equations},
    VOLUME = {47},
      YEAR = {2022},
    NUMBER = {6},
     PAGES = {1176--1192},
}

@article {BCST19-2,
    AUTHOR = {Bao, Weizhu and Carles, R\'emi and Su, Chunmei and Tang,
              Qinglin},
     TITLE = {Error estimates of a regularized finite difference method for
              the logarithmic {S}chr\"odinger equation},
   JOURNAL = {SIAM J. Numer. Anal.},
  FJOURNAL = {SIAM Journal on Numerical Analysis},
    VOLUME = {57},
      YEAR = {2019},
    NUMBER = {2},
     PAGES = {657--680},
}

@article {BCST22,
    AUTHOR = {Bao, Weizhu and Carles, R\'emi and Su, Chunmei and Tang,
              Qinglin},
     TITLE = {Error estimates of local energy regularization for the
              logarithmic {S}chr\"odinger equation},
   JOURNAL = {Math. Models Methods Appl. Sci.},
  FJOURNAL = {Mathematical Models and Methods in Applied Sciences},
    VOLUME = {32},
      YEAR = {2022},
    NUMBER = {1},
     PAGES = {101--136},
}

@article {CC23,
    AUTHOR = {Carles, R\'emi and Su, Chunmei},
     TITLE = {Numerical study of the logarithmic {S}chr\"odinger equation
              with repulsive harmonic potential},
   JOURNAL = {Discrete Contin. Dyn. Syst. Ser. B},
  FJOURNAL = {Discrete and Continuous Dynamical Systems. Series B. A Journal
              Bridging Mathematics and Sciences},
    VOLUME = {28},
      YEAR = {2023},
    NUMBER = {5},
     PAGES = {3136--3159},
}

@article{WYZ24,
    author = {Wang, Li-Lian and Yan, Jingye and Zhang, Xiaolong},
    title = {Error analysis of a first-order {IMEX} scheme for the logarithmic {S}chrödinger equation},
    journal = {SIAM J. Numer. Anal.},
    volume = {62},
    number = {1},
    pages = {119--137},
    year = {2024},
}

@article{LWZ26,
    author = {Li, Meng and Wang, Huifei and Zhao, Yanmin},
    title = {Error analysis of a second-order {IMEX} scheme for the logarithmic {S}chrödinger equation},
    journal = {Appl. Numer. Math.},
    volume = {224},
    pages = {241--261},  
    year = {2026},
}

@article{ZW25,
    author = {Zhang, Xiaolong and Wang, Li-Lian},
    title = {Low regularity estimates of the {L}ie-{T}rotter time-splitting {F}ourier spectral method for the logarithmic {S}chrödinger equation},
    journal = {Math. Comp.},
    volume = {95},
    number = {358},
    pages = {773--801},
    year = {2025},
}

@article{WY26,
    author = {Wang, Tingchun and Yan, Jingye},
    title = {Optimal error estimates of two linearized finite difference schemes for the logarithmic {S}chrödinger equation},
    journal = {Appl. Numer. Math.},
    volume = {226},
    pages = {106--118},
    year = {2026},
}

@article {GLN10,
    AUTHOR = {Guerrero, P. and L\'opez, J. L. and Nieto, J.},
     TITLE = {Global {$H^1-$}solvability of the 3{D} logarithmic
              {S}chr\"odinger equation},
   JOURNAL = {Nonlinear Anal. Real World Appl.},
  FJOURNAL = {Nonlinear Analysis. Real World Applications. },
    VOLUME = {11},
      YEAR = {2010},
    NUMBER = {1},
     PAGES = {79--87},
}

@article{lasserlubich2020,
  title={Computing quantum dynamics in the semiclassical regime},
  author={Lasser, Caroline and Lubich, Christian},
  journal={Acta Numer.},
  volume={29},
  pages={229--401},
  year={2020},
  publisher={Cambridge University Press}
}

@article {lasserChimistesMag,
    AUTHOR = {Burkhard, Selina and D\"orich, Benjamin and Hochbruck, Marlis
              and Lasser, Caroline},
     TITLE = {Variational {G}aussian approximation for the magnetic
              {S}chr\"odinger equation},
   JOURNAL = {J. Phys. A},
  FJOURNAL = {Journal of Physics. A. Mathematical and Theoretical},
    VOLUME = {57},
      YEAR = {2024},
    NUMBER = {29},
     PAGES = {Paper No. 295202, 36},
}

@book {LubichBook,
    AUTHOR = {Lubich, Christian},
     TITLE = {From quantum to classical molecular dynamics: reduced models
              and numerical analysis},
    SERIES = {Zurich Lectures in Advanced Mathematics},
 PUBLISHER = {European Mathematical Society (EMS), Z\"urich},
      YEAR = {2008},
     PAGES = {x+144},
      ISBN = {978-3-03719-067-8},
}

@article {FaouLubich,
    AUTHOR = {Faou, Erwan and Lubich, Christian},
     TITLE = {A {P}oisson integrator for {G}aussian wavepacket dynamics},
   JOURNAL = {Comput. Vis. Sci.},
  FJOURNAL = {Computing and Visualization in Science},
    VOLUME = {9},
      YEAR = {2006},
    NUMBER = {2},
     PAGES = {45--55},
}

@article {Heller76,
    AUTHOR = {Heller, Eric. J.},
     TITLE = {Time dependent variational approach to semiclassical dynamics},
   JOURNAL = {J. Chem. Phys},
  FJOURNAL = {Journal of Chemical Physics},
    VOLUME = {1},
      YEAR = {1976},
    NUMBER = {64 (1)},
     PAGES = {63-73},
}

@article{Heller81,
    author = {Heller, Eric J.},
    title = {Frozen Gaussians: A very simple semiclassical approximation},
    journal = {J. of Chem. Phys.},
    volume = {75},
    number = {6},
    pages = {2923-2931},
    year = {1981},
}

@article{CoalsonKarplus,
    author = {Coalson, Rob D. and Karplus, Martin},
    title = {Multidimensional variational {G}aussian wave packet dynamics with application to photodissociation spectroscopy},
    journal = {J. Chem. Phys.},
    volume = {93},
    number = {6},
    pages = {3919-3930},
    year = {1990},
}

@article{Raab,
title = {On the {D}irac–{F}renkel/{M}c{L}achlan variational principle},
journal = {Chem. Phys. Lett.},
volume = {319},
number = {5},
pages = {674-678},
year = {2000},
author = {A. Raab},
}

@book{KramerSaraceno,
  title     = {Geometry of the Time-Dependent Variational Principle in Quantum Mechanics},
  author    = {Kramer, Peter and Saraceno, Marcos},
  series    = {Lecture Notes in Physics},
  volume    = {140},
  year      = {1981},
  publisher = {Springer Berlin Heidelberg},
}

@article{MCTDH,
title = {The multi-configurational time-dependent {H}artree approach},
author = {H.-D. Meyer and U. Manthe and L.S. Cederbaum},
journal = {Chem. Phys. Lett.},
volume = {165},
number = {1},
pages = {73-78},
year = {1990},
}

@article{MCTDH99,
    author = {Burghardt, I. and Meyer, H.-D. and Cederbaum, L. S.},
    title = {Approaches to the approximate treatment of complex molecular systems by the multiconfiguration time-dependent {H}artree method},
    journal = {J. Chem. Phys.},
    volume = {111},
    number = {7},
    pages = {2927-2939},
    year = {1999},
}

@article{MCTDHreport,
title = {The multiconfiguration time-dependent {H}artree ({M}{C}{T}{D}{H}) method: a highly efficient algorithm for propagating wavepackets},
journal = {Physics Reports},
volume = {324},
number = {1},
pages = {1-105},
year = {2000},
author = {M.H. Beck and A. Jäckle and G.A. Worth and H.-D. Meyer},
}

@article{vMCG,
author = {G.W. Richings and I. Polyak and K.E. Spinlove and G.A. Worth and I. Burghardt and B. Lasorne},
title = {Quantum dynamics simulations using {G}aussian wavepackets: the v{M}{C}{G} method},
journal = {Int. Rev. Phys. Chem.},
volume = {34},
number = {2},
pages = {269--308},
year = {2015},
}
\bibliographystyle{alpha}
\end{document}